\documentclass[amscd, amstex]{amsart}
\usepackage{amscd}
\usepackage{amsmath}
\usepackage{amssymb}
\usepackage{float}
\theoremstyle{plain}
\newtheorem{theorem}{Theorem}[section]
\newtheorem{proposition}[theorem]{Proposition}
\newtheorem{lemma}[theorem]{Lemma}

\newtheorem{corollary}[theorem]{Corollary}
\theoremstyle{definition}
\newtheorem{definition}[theorem]{Definition}

\theoremstyle{remark}
\newtheorem{example}[theorem]{Example}

\newtheorem{question}[theorem]{Question}
\newenvironment{pf}{\begin{proof}}{\end{proof}}
\makeatletter

\@addtoreset{equation}{section}
\makeatother

\begin{document}

\title[]
{$K3$ surfaces with involution, equivariant analytic torsion, and automorphic forms on the moduli space II:
a structure theorem for $r(M)>10$}
\author[]{Ken-Ichi Yoshikawa}
\address{Department of Mathematics, 
Faculty of Science,
Kyoto University,
Kyoto 606-8502, JAPAN}
\email{yosikawa@math.kyoto-u.ac.jp}
\address{Korea Institute for Advanced Study,
Hoegiro 87, Dongdaemun-gu,
Seoul 130-722, KOREA}

\begin{abstract}
We study the structure of the invariant of $K3$ surfaces with involution, which we obtained
using equivariant analytic torsion. It was known before that the invariant is expressed as the Petersson norm
of an automorphic form on the moduli space. 
When the rank of the invariant sublattice of the $K3$-lattice with respect to the involution is strictly bigger 
than $10$, we prove that this automorphic form is expressed as the tensor product of an explicit Borcherds lift 
and Igusa's Siegel modular form.
\end{abstract}

\maketitle

\tableofcontents

\section*{Introduction}\label{sect:0}
\par
In this paper, we study the structure of the invariant of $K3$ surfaces with involution introduced 
in \cite{Yoshikawa04}. Let us recall briefly this invariant.
\par
A $K3$ surface with holomorphic involution $(X,\iota)$ is called a $2$-elementary $K3$ surface 
if $\iota$ acts non-trivially on the holomorphic $2$-forms on $X$. Let ${\Bbb L}_{K3}$ be the $K3$ lattice, i.e.,
an even unimodular lattice of signature $(3,19)$, which is isometric to $H^{2}(X,{\bf Z})$ endowed with
the cup-product pairing. Let $M$ be a sublattice of ${\Bbb L}_{K3}$ with rank $r(M)$. 
A $2$-elementary $K3$ surface $(X,\iota)$ is of type $M$ if the invariant sublattice of $H^{2}(X,{\bf Z})$ 
with respect to the $\iota$-action is isometric to $M$. 
By \cite{Nikulin80a}, $M\subset{\Bbb L}_{K3}$ must be a primitive $2$-elementary Lorentzian sublattice.
The parity of the $2$-elementary lattice $M$ is denoted by $\delta(M)\in\{0,1\}$ (cf. \cite{Nikulin83}
and Sect.\,\ref{subsect:1.2}).
\par
Let $M^{\perp}$ be the orthogonal complement of $M$ in ${\Bbb L}_{K3}$. 
Let $\Omega_{M^{\perp}}$ be the period domain for $2$-elementary $K3$ surfaces of type $M$, 
which is an open subset of a quadric hypersurface of ${\bf P}(M^{\perp}\otimes{\bf C})$.
We fix a connected component $\Omega_{M^{\perp}}^{+}$ of $\Omega_{M^{\perp}}$, which is isomorphic to 
a bounded symmetric domain of type IV of dimension $20-r(M)$. 
Let ${\mathcal D}_{M^{\perp}}$ be the discriminant locus of $\Omega_{M^{\perp}}^{+}$, 
which is a reduced divisor on $\Omega_{M^{\perp}}^{+}$.
Let $O(M^{\perp})$ be the group of isometries of $M^{\perp}$.
Then $O(M^{\perp})$ acts properly discontinuously on $\Omega_{M^{\perp}}$ and ${\mathcal D}_{M^{\perp}}$.
Let $O^{+}(M^{\perp})$ be the subgroup of $O(M^{\perp})$ with index $2$ that preserves 
$\Omega_{M^{\perp}}^{+}$.
The coarse moduli space of $2$-elementary $K3$ surfaces of type $M$ is isomorphic to the analytic space
${\mathcal M}_{M^{\perp}}^{o}=(\Omega_{M^{\perp}}^{+}\setminus{\mathcal D}_{M^{\perp}})/O^{+}(M^{\perp})$ 
via the period map by the global Torelli theorem \cite{PSS71}, \cite{BurnsRapoport75}.
The period of a $2$-elementary $K3$ surface $(X,\iota)$ of type $M$ is denoted by
$\overline{\varpi}_{M}(X,\iota)\in{\mathcal M}_{M^{\perp}}^{o}$.
\par
Let $(X,\iota)$ be a $2$-elementary $K3$ surface of type $M$.
Let $\kappa$ be an $\iota$-invariant K\"ahler form on $X$. 
Let $X^{\iota}$ be the set of fixed points of $\iota$ and let $X^{\iota}=\sum_{i}C_{i}$ be the decomposition 
into the connected components. Let $\eta\in H^{0}(X,\Omega^{2}_{X})\setminus\{0\}$.
In \cite{Yoshikawa04}, we introduced a real-valued invariant
$$
\begin{aligned}
\tau_{M}(X,\iota)
&=
{\rm vol}(X,(2\pi)^{-1}\kappa)^{\frac{14-r(M)}{4}}\tau_{{\bf Z}_{2}}(X,\kappa)(\iota)
\prod_{i}{\rm Vol}(C_{i},(2\pi)^{-1}\kappa|_{C_{i}})\tau(C_{i},\kappa|_{C_{i}})
\\
&\quad
\times
\exp\left[
\frac{1}{8}\int_{X^{\iota}}
\log\left.\left(\frac{\eta\wedge\bar{\eta}}
{\kappa^{2}/2!}\cdot\frac{{\rm Vol}(X,(2\pi)^{-1}\kappa)}{\|\eta\|_{L^{2}}^{2}}
\right)\right|_{X^{\iota}}
c_{1}(X^{\iota},\kappa|_{X^{\iota}})
\right],
\end{aligned}
$$
where $\tau_{{\bf Z}_{2}}(X,\kappa)(\iota)$ is the equivariant analytic torsion of $(X,\kappa)$ 
with respect to the ${\bf Z}_{2}$-action induced by $\iota$, $\tau(C_{i},\kappa|_{C_{i}})$ is the analytic torsion
of $(C_{i},\kappa|_{C_{i}})$, and $c_{1}(X^{\iota},\kappa|_{X^{\iota}})$ is the first Chern form of 
$(X^{\iota},\kappa|_{X^{\iota}})$. (See \cite{Bismut95}, \cite{BGS88}, \cite{RaySinger73} and Sect.\,\ref{sect:5}.)
Since $\tau_{M}(X,\iota)$ depends only on the isomorphism class of $(X,\iota)$, we get the function
$$
\tau_{M}\colon{\mathcal M}_{M^{\perp}}^{o}
\ni\overline{\varpi}_{M}(X,\iota)\to\tau_{M}(X,\iota)\in{\bf R}_{>0}.
$$
By \cite{Yoshikawa04}, \cite{Yoshikawa09b}, there exists an automorphic form 
$\Phi_{M}$ on $\Omega_{M^{\perp}}^{+}$ with values in a certain $O^{+}(M^{\perp})$-equivariant 
holomorphic line bundle on $\Omega_{M^{\perp}}^{+}$, such that 
$$
\tau_{M}=\|\Phi_{M}\|^{-\frac{1}{2\nu}},
\qquad
{\rm div}\,\Phi_{M}=\nu\,{\mathcal D}_{M^{\perp}},
\qquad
\nu\in{\bf Z}_{>0}.
$$
Here $\|\cdot\|$ denotes the Petersson norm.
By \cite{Yoshikawa04}, $\Phi_{M}$ is given by the Borcherds $\Phi$-function \cite{Borcherds95}, 
\cite{Borcherds96} when $M$ is one of the two exceptional lattices in Proposition~\ref{Proposition2.1}.
For an arithmetic counterpart of the invariant $\tau_{M}$, 
we refer the reader to \cite{MaillotRoessler08}.
\par
In this paper, we give an explicit formula for $\tau_{M}$ for a class of non-exceptional $M$. 
We use two kinds of automorphic forms to express $\tau_{M}$, i.e., the Borcherds lift 
$\Psi_{M^{\perp}}(\cdot,F_{M^{\perp}})$ and Igusa's Siegel modular form $\chi_{g}$, which we explain briefly.
\par 
In \cite{Borcherds95}, \cite{Borcherds98}, Borcherds developed the theory of automorphic forms
with infinite product expansion over domains of type IV.
For an even $2$-elementary lattice $\Lambda$ of signature $(2,r(\Lambda)-2)$, we define the Borcherds lift 
$\Psi_{\Lambda}(\cdot,F_{\Lambda})$ as follows.
\par
Let $A_{\Lambda}$ be the discriminant group of $\Lambda$, which is a vector space over ${\bf Z}/2{\bf Z}$. 
Let ${\bf C}[A_{\Lambda}]$ be the group ring of $A_{\Lambda}$ and 
let $\rho_{\Lambda}\colon Mp_{2}({\bf Z})\to GL({\bf C}[A_{\Lambda}])$ be the Weil representation,
where $Mp_{2}({\bf Z})$ is the metaplectic double cover of $SL_{2}({\bf Z})$. 
Let $\{{\bf e}_{\gamma}\}_{\gamma\in A_{\Lambda}}$ be the standard basis of ${\bf C}[A_{\Lambda}]$.
Let $\eta(\tau)$ be the Dedekind $\eta$-function and set 
$\eta_{1^{-8}2^{8}4^{-8}}(\tau)=\eta(\tau)^{-8}\eta(2\tau)^{8}\eta(4\tau)^{-8}$.
Let $\theta_{{\Bbb A}_{1}^{+}}(\tau)$ be the theta function of the (positive-definite) $A_{1}$-lattice.
Then $\eta_{1^{-8}2^{8}4^{-8}}(\tau)$ and $\theta_{{\Bbb A}_{1}^{+}}(\tau)$ are modular forms 
for the subgroup $M\Gamma_{0}(4)\subset Mp_{2}({\bf Z})$ corresponding to the congruence subgroup
$\Gamma_{0}(4)\subset SL_{2}({\bf Z})$.
Following \cite{Borcherds00} and \cite{Scheithauer06}, we define a ${\bf C}[A_{\Lambda}]$-valued
holomorphic function $F_{\Lambda}(\tau)$ on the complex upper half-plane ${\frak H}$ as
\begin{equation}\label{eqn:(0.1)}
F_{\Lambda}(\tau)
=
\sum_{g\in M\Gamma_{0}(4)\backslash Mp_{2}({\bf Z})}
\left.\left\{\eta_{1^{-8}2^{8}4^{-8}}
\theta_{{\Bbb A}_{1}^{+}}^{12-r(\Lambda)}
\right\}\right|_{g}(\tau)\,
\rho_{\Lambda}(g^{-1})\,{\bf e}_{0}.
\end{equation}
Here we used the notation $\phi|_{g}(\tau)=\phi(\frac{a\tau+b}{c\tau+d})(c\tau+d)^{-k}$
for a modular form $\phi(\tau)$ for $M\Gamma_{0}(4)$ of weight $k$ with certain character and 
$g=(\binom{a\,b}{c\,d},\sqrt{c\tau+d})\in Mp_{2}({\bf Z})$.
By \cite{Borcherds00} and \cite{Scheithauer06}, $F_{\Lambda}(\tau)$ is an elliptic modular form for 
$Mp_{2}({\bf Z})$ of type $\rho_{\Lambda}$ with weight $2-\frac{r(\Lambda)}{2}$.
Then $\Psi_{\Lambda}(\cdot,F_{\Lambda})$ is defined as the Borcherds lift of $F_{\Lambda}(\tau)$, 
which is an automorphic form on $\Omega_{\Lambda}^{+}$ for $O^{+}(\Lambda)$ by \cite{Borcherds98}. 
(See \eqref{eqn:(8.3)} for an explicit infinite product expression of $\Psi_{\Lambda}(\cdot,F_{\Lambda})$.)
The Petersson norm $\|\Psi_{M^{\perp}}(\cdot,F_{M^{\perp}})\|^{2}$ is an $O^{+}(M^{\perp})$-invariant function 
on $\Omega_{M^{\perp}}^{+}$ and the value $\|\Psi_{M^{\perp}}(\overline{\varpi}_{M}(X,\iota),F_{M^{\perp}})\|$ 
makes sense.
\par
Recall that $\chi_{g}$ is the Siegel modular form on the Siegel upper half-space ${\frak S}_{g}$ 
of degree $g$ defined as the product of all even theta constants (cf. \cite{Igusa67})
\begin{equation}\label{eqn:(0.2)}
\chi_{g}(\varSigma)=\prod_{(a,b)\,{\rm even}}\theta_{a,b}(\varSigma),
\qquad
\varSigma\in{\frak S}_{g},
\qquad\quad
\chi_{0}=1.
\end{equation}
Then $\chi_{g}$ gives rise to another function on ${\mathcal M}_{M^{\perp}}^{o}$ as follows.
For a $2$-elementary $K3$ surface $(X,\iota)$ of type $M$, let $X^{\iota}$ denote the set of fixed points of $\iota$.
By \cite{Nikulin83}, $X^{\iota}$ is the disjoint union of (possibly empty) compact Riemann surfaces, 
whose topological type is determined by $M$. Let $g(M)\in{\bf Z}_{\geq0}$ denote the total genus of $X^{\iota}$. 
The period of $X^{\iota}$ is denoted by $\varOmega(X^{\iota})\in{\frak S}_{g(M)}/Sp_{2g(M)}({\bf Z})$.
By \cite{Yoshikawa04}, there exist a proper $O^{+}(M^{\perp})$-invariant Zariski closed subset 
$Z\subset{\mathcal D}_{M^{\perp}}$ and an $O^{+}(M^{\perp})$-equivariant holomorphic map 
$J_{M}\colon\Omega_{M^{\perp}}\setminus Z\to{\frak S}_{g(M)}/Sp_{2g(M)}({\bf Z})$ 
that induces the map of moduli spaces
$$
{\mathcal M}_{M^{\perp}}^{o}
\ni\overline{\varpi}_{M}(X,\iota)\to\varOmega(X^{\iota})\in
{\frak S}_{g(M)}/Sp_{2g(M)}({\bf Z}).
$$
Then $J_{M}^{*}\|\chi_{g(M)}\|^{2}$ is an $O^{+}(M^{\perp})$-invariant $C^{\infty}$ function 
on $\Omega_{M^{\perp}}^{o}$.
\par 
The following structure theorem for $\tau_{M}$ is the main result of this paper:

\begin{theorem}\label{Theorem0.1}
{\rm (cf. Theorem~\ref{Theorem9.1})}
Let $M$ be a primitive $2$-elementary Lorentzian sublattice of ${\Bbb L}_{K3}$. 
If $r(M)>10$ or $(r(M),\delta(M))=(10,1)$, then there exists a constant $C_{M}$ depending only on the lattice $M$ 
such that the following identity holds for every $2$-elementary $K3$ surface $(X,\iota)$ of type $M$
$$
\tau_{M}(X,\iota)^{-2^{g(M)+1}(2^{g(M)}+1)}
=
C_{M}\,
\|\Psi_{M^{\perp}}(\overline{\varpi}_{M}(X,\iota),
F_{M^{\perp}})\|^{2^{g(M)}}
\|\chi_{g(M)}(\varOmega(X^{\iota}))\|^{16}.
$$
\end{theorem}

It may be worth emphasizing that the structure of $\tau_{M}$ becomes transparent 
by considering elliptic modular forms for $M\Gamma_{0}(4)$ rather than $Mp_{2}({\bf Z})$.
After Bruinier \cite{Bruinier02}, Theorem~\ref{Theorem0.1} may not be surprizing. 
Indeed, if $M^{\perp}$ contains an even unimodular lattice of signature $(2,2)$ as a direct summand 
and if there is a Siegel modular form $S$ such that ${\rm div}(J_{M}^{*}S)$ is a Heegner divisor on 
$\Omega_{M^{\perp}}^{+}$, then $\Phi_{M}$ must be the product of a Borcherds lift and $J_{M}^{*}S$ 
by \cite[Th.\,0.8]{Bruinier02}. When $g(M)=2$, this explains the existence of a factorization of $\tau_{M}$ like 
Theorem~\ref{Theorem0.1}. 
(However, this does not seem to explain the common structures \eqref{eqn:(0.1)}, \eqref{eqn:(0.2)} 
of the elliptic modular forms $F_{\Lambda}$ and the Siegel modular forms $\chi_{g}$ appearing
in the expression of $\tau_{M}$.) It seems to be an interesting problem to understand the geometric origin of 
these common structures of the modular forms $F_{\Lambda}$ and $\chi_{g}$.
\par
There are $43$ isometry classes of primitive $2$-elementary Lorentzian sublattices $M\subset{\Bbb L}_{K3}$ 
such that $r(M)>10$ or $(r(M),\delta(M))=(10,1)$ (cf. Table 1 in Sect.\,\ref{sect:9}).
In fact, Theorem~\ref{Theorem0.1} remains valid for a certain primitive $2$-elementary Lorentzian sublattice 
$M\subset{\Bbb L}_{K3}$ with $r(M)=9$. (See Theorem~\ref{Theorem9.3}.)
By Theorems~\ref{Theorem0.1} and \ref{Theorem9.3} and \cite[Ths.\,8.2 and 8.7]{Yoshikawa04}, 
$\tau_{M}$ and $\Phi_{M}$ are determined for $46$ isometry classes of $M$. Since the total number of 
the isometry classes of primitive $2$-elementary Lorentzian sublattice of ${\Bbb L}_{K3}$ is $75$ by 
Nikulin \cite{Nikulin83}, the structures of $\tau_{M}$ and $\Phi_{M}$ are still open for the remaining $29$ lattices.
\par
Following \cite[Th.\,8.7]{Yoshikawa04}, 
we shall prove Theorem~\ref{Theorem0.1} by comparing the $O^{+}(M^{\perp})$-invariant currents 
$dd^{c}\log\tau_{M}$, $dd^{c}\log\|\Psi_{M^{\perp}}(\cdot,F_{M^{\perp}})\|^{2}$ and 
$dd^{c}\log J_{M}^{*}\|\chi_{g(M)}^{8}\|^{2}$. 
(See Sect.\,\ref{sect:9}.) The current $dd^{c}\log\tau_{M}$ was computed in \cite{Yoshikawa04}.
In Sect.\,\ref{sect:8}, the weight and the zero divisor of $\Psi_{M^{\perp}}(\cdot,F_{M^{\perp}})$ 
shall be determined, from which a formula for $dd^{c}\log\|\Psi_{M^{\perp}}(\cdot,F_{M^{\perp}})\|^{2}$ follows. 
In Sect.\,\ref{sect:4}, the current $dd^{c}\log J_{M}^{*}\|\chi_{g(M)}^{8}\|^{2}$ shall be computed, 
where the irreducibility of certain component of the divisor ${\mathcal D}_{M^{\perp}}/O^{+}(M^{\perp})$
on $\Omega^{+}_{M^{\perp}}/O^{+}(M^{\perp})$ plays a crucial role. (See Appendix Sect.\,\ref{subsect:11.3}.)
\par
In Proposition~\ref{Proposition9.2}, we shall prove that $\chi_{g(M)}$ vanishes identically on the locus
$J_{M}(\Omega_{M^{\perp}}^{+}\setminus{\mathcal D}_{M^{\perp}})$
when $(r(M),\delta(M))=(10,0)$ and $M$ is not exceptional. 
Hence Theorem~\ref{Theorem0.1} does {\em not} hold in these four cases.
This is similar to the exceptional case $(r(M),l(M),\delta(M))=(10,8,0)$, where $\chi_{g}$ should be replaced by
the product of two Jacobi's $\Delta$-functions \cite[Th.\,8.7]{Yoshikawa04}. 
\par
There is an application of the Borcherds lift $\Psi_{\Lambda}(\cdot,F_{\Lambda})$ to the moduli space
${\mathcal M}_{M^{\perp}}^{o}$.

\begin{theorem}\label{Theorem0.2}
If $r(M)\geq10$, then ${\mathcal M}_{M^{\perp}}^{o}$ is quasi-affine.
\end{theorem}

When ${\mathcal M}_{M^{\perp}}^{o}$ is the coarse moduli space of Enriques surfaces, 
this was proved by Borcherds \cite{Borcherds96}.
Since the coarse moduli space of ample $M$-polarized $K3$ surfaces 
(cf. \cite{AlexeevNikulin06}, \cite{Dolgachev96}, \cite{Nikulin80b}) is a finite covering of 
${\mathcal M}_{M^{\perp}}^{o}$, its quasi-affiness follows from that of ${\mathcal M}_{M^{\perp}}^{o}$.
The quasi-affiness of ${\mathcal M}_{M^{\perp}}^{o}$ is a consequence of the fact that
$\Psi_{M^{\perp}}(\cdot,F_{M^{\perp}})$ vanishes only on the discriminant locus ${\mathcal D}_{M^{\perp}}$ 
when $r(M^{\perp})\leq12$.
By \cite{Nikulin83}, there are $49$ isometry classes of primitive $2$-elementary Lorentzian sublattices 
$M\subset {\Bbb L}_{K3}$ with $r(M)\geq10$. 
In general, it is not easy to find a primitive sublattice $\Lambda\subset{\Bbb L}_{K3}$ of signature 
$(2,r(\Lambda)-2)$ admitting an automorphic form on $\Omega_{\Lambda}^{+}$ 
vanishing only on ${\mathcal D}_{\Lambda}$.
For example, there is {\it no} automorphic form on the coarse moduli space of polarized $K3$ surfaces of 
degree $2d$ vanishing only on the discriminant locus, 
if the discriminant locus is irreducible \cite[Sect.\,3.3]{Looijenga03}, \cite{Nikulin96}.
(See \cite{Borcherds95}, \cite{Borcherds96}, \cite{Borcherds00}, \cite{BKPSB98}, 
\cite[II]{GritsenkoNikulin98}, \cite{Kondo07a}, \cite{Scheithauer06}, \cite{Gritsenko10}... for affirmative examples).
For another application of $\Psi_{\Lambda}(\cdot,F_{\Lambda})$ to the negativity of the Kodaira dimension of
${\mathcal M}_{M^{\perp}}=\Omega_{M^{\perp}}^{+}/O^{+}(M^{\perp})$, see \cite{MaS10}. 
In fact, ${\mathcal M}_{M^{\perp}}$ is always unirational and hence $\kappa({\mathcal M}_{M^{\perp}})=-\infty$ 
by S. Ma \cite{MaS10}.
\par
This paper is organized as follows.
In Sect.\,\ref{sect:1}, we recall lattices and orthogonal modular varieties.
In Sect.\,\ref{sect:2}, we recall $2$-elementary $K3$ surfaces and their moduli spaces, 
and we study the singular fiber of an ordinary singular family of $2$-elementary $K3$ surfaces.
In Sect.\,\ref{sect:3}, we recall log del Pezzo surfaces of index $\leq2$ and its relation with 
$2$-elementary $K3$ surfaces.
In Sect.\,\ref{sect:4}, we study the current $dd^{c}J_{M}^{*}\|\chi_{g(M)}^{8}\|^{2}$ and we recall 
the notion of automorphic forms on $\Omega_{M^{\perp}}^{+}$.
In Sect.\,\ref{sect:5}, we recall the invariant $\tau_{M}$.
In Sect.\,\ref{sect:6}, we recall Borcherds products.
In Sect.\,\ref{sect:7}, we construct the elliptic modular form $F_{\Lambda}(\tau)$. 
In Sect.\,\ref{sect:8}, we study the Borcherds lift $\Psi_{\Lambda}(\cdot,F_{\Lambda})$.
In Sect.\,\ref{sect:9}, we prove Theorem~\ref{Theorem0.1}.
In Sect.\,\ref{sect:10}, we interpret Theorem~\ref{Theorem0.1} into a statement about 
the equivariant determinant of the Laplacian on real $K3$ surfaces.
In Sect.\,\ref{sect:11} Appendix, we prove some technical results about lattices.
\par
{\it Warning: }
In \cite{Yoshikawa04}, we used the notation $\Omega_{M}$, ${\mathcal M}_{M}$, ${\mathcal D}_{M}$ etc.
in stead of $\Omega_{M^{\perp}}$, ${\mathcal M}_{M^{\perp}}$, ${\mathcal D}_{M^{\perp}}$ etc.
\par
{\bf Acknowledgements }
We thank Professors Jean-Michel Bismut, Jean-Beno\^\i{}t Bost, Kai K\"ohler, Vincent Maillot, 
Damian Roessler for helpful discussions about the invariant $\tau_{M}$ and the automorphic form 
$\Phi_{M}$, in particular their arithmetic nature.
We thank Professor Siye Wu for helpful discussions about the elliptic modular form $F_{\Lambda}(\tau)$ 
and $S$-duality.
The main result of this paper is a refinement of \cite[Th.\,8.2]{Yoshikawa98}.
We thank Professor Richard Borcherds, from whom we learned how to construct modular forms of 
type $\rho_{\Lambda}$ from scalar-valued modular forms for $\Gamma_{0}(4)$ 
(cf. Proposition~\ref{Proposition7.1}, \cite[Prop.\,8.1]{Yoshikawa98}); 
the construction of $F_{\Lambda}(\tau)$ is an obvious extension of that of $F_{k}(\tau)$ 
in \cite[Sect.\,8.2]{Yoshikawa98}.
We also thank Professor Shigeyuki Kond\=o, who suggested us the proofs of some cases of 
Propositions~\ref{Proposition3.1} and \ref{Proposition9.2} in the earlier version of this paper.\footnote{
http://kyokan.ms.u-tokyo.ac.jp/users/preprint/pdf/2007-12.pdf}
We are partially supported by the Grants-in-Aid for Scientific Research (B) 19340016 and (S) 17104001, JSPS.

\section{Lattices}
\label{sect:1}
\par
A free ${\bf Z}$-module of finite rank endowed with a non-degenerate, integral, symmetric bilinear form 
is called a lattice.
The rank of a lattice $L$ is denoted by $r(L)$.
The signature of $L$ is denoted by ${\rm sign}(L)=(b^{+}(L),b^{-}(L))$.
We define $\sigma(L):=b^{+}(L)-b^{-}(L)$. 
A lattice $L$ is {\it Lorentzian} if ${\rm sign}(L)=(1,r(L)-1)$. 
For a lattice $L=({\bf Z}^{r},\langle\cdot,\cdot\rangle)$, we define $L(k):=({\bf Z}^{r},k\langle\cdot,\cdot\rangle)$.
The dual lattice of $L$ is denoted by $L^{\lor}\subset L\otimes{\bf Q}$.
The finite abelian group $A_{L}:=L^{\lor}/L$ is called the {\it discriminant group} of $L$.
For $\lambda\in L^{\lor}$, we write $\bar{\lambda}:=\lambda+L\in A_{L}$.
A lattice $L$ is {\it even} if $\langle x,x\rangle\in2{\bf Z}$ for all $x\in L$.
A sublattice $M\subset L$ is {\it primitive} if $L/M$ has no torsion elements.
The {\it level} of an even lattice $L$ is the smallest positive integer $l$ such that
$l\,\lambda^{2}/2\in{\bf Z}$ for all $\lambda\in L^{\lor}$.
The group of isometries of $L$ is denoted by $O(L)$.
We set $\Delta_{L}:=\{d\in L;\,\langle d,d\rangle=-2\}$ and define
$$
\Delta'_{L}:=\{d\in\Delta_{L},\,d/2\not\in L^{\lor}\},
\qquad
\Delta''_{L}:=\{d\in\Delta_{L},\,d/2\in L^{\lor}\}.
$$
Then $\Delta_{L}$, $\Delta_{L}'$, $\Delta_{L}''$ are preserved by $O(L)$. 
For $d\in\Delta_{L}$, the corresponding reflection $s_{d}\in O(L)$ is defined as 
$s_{d}(x):=x+\langle x,d\rangle\,d$ for $x\in L$.

\subsection{Discriminant forms}
\label{subsect:1.1}
\par
For an even lattice $L$, the {\it discriminant form} $q_{L}$ of $A_{L}$ is the quadratic form on $A_{L}$ 
with values in ${\bf Q}/2{\bf Z}$ defined as $q_{L}(\bar{l}):=l^{2}+2{\bf Z}$ for $\bar{l}\in A_{L}$. 
The corresponding bilinear form on $A_{L}$ with values in ${\bf Q}/{\bf Z}$ is denoted by $b_{L}$. 
Then $b_{L}(\bar{l},\bar{l'})=\langle l,l'\rangle+{\bf Z}$ for $\bar{l},\bar{l'}\in A_{L}$.
Since $\lambda\in L^{\lor}$ lies in $L$ if and only if $\langle\lambda,l\rangle\in{\bf Z}$ for all $l\in L^{\lor}$,
the bilinear form $b_{L}$ is non-degenerate, i.e.,
if $b_{L}(\gamma,x)\equiv0\mod{\bf Z}$ for all $x\in A_{L}$, then $\gamma=0$ in $A_{L}$.
We often write $\gamma^{2}$ (resp. $\langle\gamma,\delta\rangle$)
for $q_{L}(\gamma)$ (resp. $b_{L}(\gamma,\delta)$).
The automorphism group of $(A_{L},q_{L})$ is denoted by $O(q_{L})$.
See \cite{Nikulin80a} for more details.

\subsection{$2$-elementary lattices}
\label{subsect:1.2}
\par
Set ${\bf Z}_{2}:={\bf Z}/2{\bf Z}$.
An even lattice $L$ is {\it $2$-elementary} if there is an integer $l\in{\bf Z}_{\geq0}$ with 
$A_{L}\cong{\bf Z}_{2}^{l}$. 
For a $2$-elementary lattice $L$, we set $l(L):=\dim_{{\bf Z}_{2}}A_{L}$.
Then $r(L)\geq l(L)$ and $r(L)\equiv l(L)\mod 2$ by \cite[Th.\,3.6.2 (2)]{Nikulin80a}. We define
$$
\delta(L)
:=
\begin{cases}
\begin{array}{lll}
0&\mbox{if}&x^{2}\in{\bf Z}\,
\mbox{ for all }x\in L^{\lor},
\\
1&\mbox{if}&x^{2}\not\in{\bf Z}\,
\mbox{ for some }x\in L^{\lor}.
\end{array}
\end{cases}
$$
The triplet $({\rm sign}(L),l(L),\delta(L))$ determines the isometry class of an {\it indefinite} 
even $2$-elementary lattice $L$ by \cite[Th.\,3.6.2]{Nikulin80a}. 
\par
Since the mapping $A_{L}\ni\gamma\to\gamma^{2}\in\frac{1}{2}{\bf Z}/{\bf Z}\cong{\bf Z}_{2}$
is ${\bf Z}_{2}$-linear and since $b_{L}$ is non-degenerate, 
there exists a unique element ${\bf 1}_{L}\in A_{L}$ such that
$\langle\gamma,{\bf 1}_{L}\rangle\equiv\gamma^{2}\mod{\bf Z}$ for all $\gamma\in A_{L}$. 
By the uniqueness of ${\bf 1}_{L}$, $g({\bf 1}_{L})={\bf 1}_{L}$ for all $g\in O(q_{L})$.
By definition, ${\bf 1}_{L}=0$ if and only if $\delta(L)=0$.
If $L=L'\oplus L''$, then ${\bf 1}_{L}={\bf 1}_{L'}\oplus{\bf 1}_{L''}$.
\par
Let ${\Bbb U}=\binom{0\,1}{1\,0}$ and let ${\Bbb A}_{1}$, ${\Bbb D}_{2k}$, ${\Bbb E}_{7}$,
${\Bbb E}_{8}$ be the {\it negative-definite} Cartan matrix of type $A_{1}$, $D_{2k}$, $E_{7}$, 
$E_{8}$ respectively, which are identified with the corresponding even lattices. 
Then ${\Bbb U}$ and ${\Bbb E}_{8}$ are unimodular,
and ${\Bbb A}_{1}$, ${\Bbb D}_{2k}$ and ${\Bbb E}_{7}$ are $2$-elementary. Set
$$
{\Bbb L}_{K3}:={\Bbb U}\oplus{\Bbb U}\oplus{\Bbb U}\oplus{\Bbb E}_{8}\oplus{\Bbb E}_{8}.
$$
For a sublattice $\Lambda\subset{\Bbb L}_{K3}$,
we define $\Lambda^{\perp}:=\{l\in{\Bbb L}_{K3};\,\langle l,\Lambda\rangle=0\}$.
When $\Lambda\subset{\Bbb L}_{K3}$ is primitive,
then $(A_{\Lambda},-q_{\Lambda})\cong(A_{\Lambda^{\perp}},q_{\Lambda^{\perp}})$
by \cite[Cor.\,1.6.2]{Nikulin80a}.

\subsection{Lorentzian lattices}
\label{subsect:1.3}
\par
Let $L$ be a Lorentzian lattice. The set ${\mathcal C}_{L}:=\{v\in L\otimes{\bf R};\,v^{2}>0\}$
is called the positive cone of $L$. 
Since $L$ is Lorentzian, ${\mathcal C}_{L}$ consists of two connected components,
which are written as ${\mathcal C}^{+}_{L}$, ${\mathcal C}^{-}_{L}$. 
For $l\in L\otimes{\bf R}$, we set $h_{l}:=\{v\in{\mathcal C}_{L}^{+};\,\langle v,l\rangle=0\}$. 
Then $h_{l}\not=\emptyset$ if and only if $l^{2}<0$. 
We define $({\mathcal C}_{L}^{+})^{o}:={\mathcal C}_{L}^{+}\setminus\bigcup_{d\in\Delta_{L}}h_{d}$. 
Any connected component of $({\mathcal C}_{L}^{+})^{o}$ is called a {\it Weyl chamber} of $L$.
\par
Let $M\subset{\Bbb L}_{K3}$ be a primitive $2$-elementary Lorentzian sublattice. 
Let $I_{M}$ be the involution on $M\oplus M^{\perp}$ defined as 
$$
I_{M}(x,y)=(x,-y),
\qquad
(x,y)\in M\oplus M^{\perp}.
$$
Then $I_{M}$ extends uniquely to an involution on ${\Bbb L}_{K3}$ by \cite[Cor.\,1.5.2]{Nikulin80a}. 
We define
$$
g(M):=\{22-r(M)-l(M)\}/2,
\qquad
k(M):=\{r(M)-l(M)\}/2.
$$
For $d\in\Delta_{M^{\perp}}$, the smallest sublattice of ${\Bbb L}_{K3}$ containing $M$ and ${\bf Z}d$ 
is given by
$$
[M\perp d]:=(M^{\perp}\cap d^{\perp})^{\perp}.
$$ 
By Lemma~\ref{Lemma11.3} below, $[M\perp d]$ is again a $2$-elementary Lorentzian lattice such that
\begin{equation}\label{eqn:(1.1)}
I_{[M\perp d]}=s_{d}\circ I_{M}=I_{M}\circ s_{d},
\qquad
[M\perp d]^{\perp}=M^{\perp}\cap d^{\perp}.
\end{equation}
By e.g. \cite[Appendix, Tables 1,2,3]{FinashinKharlamov08},
$M$ and $M^{\perp}$ are expressed as direct sums of the $2$-elementary lattices
${\Bbb A}_{1}^{+}:={\Bbb A}_{1}(-1)$, ${\Bbb A}_{1}$, ${\Bbb U}$,
${\Bbb U}(2)$, ${\Bbb D}_{2k}$, ${\Bbb E}_{7}$,
${\Bbb E}_{8}$, ${\Bbb E}_{8}(2)$.

\subsection{Lattices of signature $(2,n)$ and orthogonal modular varieties}
\label{subsect:1.4}
\par
Let $\Lambda$ be a lattice with ${\rm sign}(\Lambda)=(2,n)$.
Define
$$
\Omega_{\Lambda}
:=
\{[x]\in{\bf P}(\Lambda\otimes{\bf C});\,
\langle x,x\rangle=0,\,
\langle x,\bar{x}\rangle>0\},
$$
which has two connected components $\Omega_{\Lambda}^{\pm}$.
Each of $\Omega_{\Lambda}^{\pm}$ is isomorphic to a bounded symmetric domain of type IV of dimension $n$. 
On $\Omega_{\Lambda}$, acts $O(\Lambda)$ projectively.
Set 
$$
O^{+}(\Lambda):=\{g\in O(\Lambda);\,g(\Omega_{\Lambda}^{\pm})=\Omega_{\Lambda}^{\pm}\},
$$
which is a subgroup of $O(\Lambda)$ of index $2$ such that 
$\Omega_{\Lambda}/O(\Lambda)=\Omega_{\Lambda}^{+}/O^{+}(\Lambda)$.
Since $O^{+}(\Lambda)$ is an arithmetic subgroup of ${\rm Aut}(\Omega_{\Lambda}^{+})$, 
$O^{+}(\Lambda)$ acts properly discontinuously on $\Omega_{\Lambda}^{+}$.
In particular, the stabilizer $O^{+}(\Lambda)_{[\eta]}:=\{g\in O^{+}(\Lambda);\,g\cdot[\eta]=[\eta]\}$
is finite for all $[\eta]\in\Omega_{\Lambda}^{+}$, and the quotient 
$$
{\mathcal M}_{\Lambda}
:=
\Omega_{\Lambda}/O(\Lambda)
=
\Omega_{\Lambda}^{+}/O^{+}(\Lambda)
$$
is an analytic space. 
There exists a compactification ${\mathcal M}_{\Lambda}^{*}$ of ${\mathcal M}_{\Lambda}$, 
called the Baily--Borel--Satake compactification \cite{BailyBorel66}, 
such that ${\mathcal M}_{\Lambda}^{*}$ is an irreducible normal projective variety of dimension $n$ with
$\dim({\mathcal M}_{\Lambda}^{*}\setminus{\mathcal M}_{\Lambda})\leq1$.
\par
For $\lambda\in\Lambda\otimes{\bf R}$, set 
$H_{\lambda}:=\{[x]\in\Omega_{\Lambda};\,\langle x,\lambda\rangle=0\}$.
Then $H_{\lambda}\not=\emptyset$ if and only if $\lambda^{2}<0$. 
We define the discriminant locus of $\Omega_{\Lambda}$ by
$$
{\mathcal D}_{\Lambda}:=\sum_{d\in\Delta_{\Lambda}/\pm1}H_{d},
$$
which is a reduced divisor on $\Omega_{\Lambda}$.
We define the reduced divisors ${\mathcal D}'_{\Lambda}$ and ${\mathcal D}''_{\Lambda}$ by
$$
{\mathcal D}'_{\Lambda}
=
\sum_{d\in\Delta'_{\Lambda}/\pm1}H_{d},
\qquad\qquad
{\mathcal D}''_{\Lambda}
=
\sum_{d\in\Delta''_{\Lambda}/\pm1}H_{d}.
$$
Since $\Delta_{\Lambda}=\Delta'_{\Lambda}\amalg\Delta''_{\Lambda}$,
we have ${\mathcal D}_{\Lambda}={\mathcal D}'_{\Lambda}+{\mathcal D}''_{\Lambda}$. 
\par
Assume that $\Lambda$ is a primitive $2$-elementary sublattice of ${\Bbb L}_{K3}$. We set
$$
\Omega_{\Lambda}^{o}
:=
\Omega_{\Lambda}\setminus{\mathcal D}_{\Lambda},
\qquad
{\mathcal M}_{\Lambda}^{o}:=\Omega_{\Lambda}^{o}/O(\Lambda).
$$ 
For $d\in\Delta_{\Lambda}$, we have the relation
$$
H_{d}\cap\Omega_{\Lambda}
=
\Omega_{\Lambda\cap d^{\perp}}
=
\Omega_{[\Lambda^{\perp}\perp d]^{\perp}}.
$$
We define the subsets $H_{d}^{o}\subset H_{d}$ ($d\in\Delta_{\Lambda}$) and
${\mathcal D}_{\Lambda}^{o}\subset{\mathcal D}_{\Lambda}$ by
$$
H_{d}^{o}:=\{[\eta]\in\Omega_{\Lambda}^{+};\,O^{+}(\Lambda)_{[\eta]}=\{\pm1,\,\pm s_{d}\}\},
\qquad
{\mathcal D}_{\Lambda}^{o}:=\sum_{d\in\Delta_{\Lambda}/\pm1}H^{o}_{d}.
$$
If $H_{d}\not=\emptyset$ (resp. ${\mathcal D}_{\Lambda}\not=\emptyset$),
then $H_{d}^{o}$ (resp. ${\mathcal D}_{\Lambda}^{o}$)
is a non-empty Zariski open subset of $\Omega_{\Lambda\cap d^{\perp}}$ (resp. ${\mathcal D}_{\Lambda}$) 
unless $M^{\perp}=({\Bbb A}_{1}^{+})^{\oplus2}\oplus{\Bbb A}_{1}$
(cf. \cite[Proof of Th.\,4.1 and Sect.\,5]{Yoshikawa09b}).
Since $O(\Lambda)$ preserves ${\mathcal D}_{\Lambda}$ and ${\mathcal D}_{\Lambda}^{o}$, we define
$$
\overline{\mathcal D}_{\Lambda}:={\mathcal D}_{\Lambda}/O(\Lambda),
\qquad
\overline{\mathcal D}_{\Lambda}^{o}:={\mathcal D}_{\Lambda}^{o}/O(\Lambda)
\subset
\overline{\mathcal D}_{\Lambda}.
$$
Then $\overline{\mathcal D}_{\Lambda}^{o}\cap{\rm Sing}\,{\mathcal M}_{\Lambda}=\emptyset$ 
by \cite[Prop.\,1.9 (5)]{Yoshikawa04}.
For the irreducibility of ${\mathcal D}'_{\Lambda}/O(\Lambda)$, 
see Proposition~\ref{Proposition11.6} (5) below.
\par
When $\Lambda={\Bbb U}(N)\oplus L$, a vector of $\Lambda\otimes{\bf C}$ is denoted by $(m,n,v)$, where 
$m,n\in{\bf C}$ and $v\in L\otimes{\bf C}$.
The tube domain $L\otimes{\bf R}+i\,{\mathcal C}_{L}$ is identified with $\Omega_{\Lambda}$ via the map 
\begin{equation}\label{eqn:(1.2)}
L\otimes{\bf R}+i\,{\mathcal C}_{L}
\ni z\to 
[(-z^{2}/2,1/N,z)]
\in\Omega_{\Lambda}\subset
{\bf P}(\Lambda\otimes{\bf C}),
\qquad
z\in L\otimes{\bf C}
\end{equation}
by \cite[p.542]{Borcherds98}.
The component of $\Omega_{\Lambda}$ corresponding to $L\otimes{\bf R}+i\,{\mathcal C}_{L}^{+}$ 
via the isomorphism \eqref{eqn:(1.2)} is written as $\Omega_{\Lambda}^{+}$.

\section{$K3$ surfaces with involution}
\label{sect:2}
\par

\subsection{$K3$ surfaces with involution and their moduli space}
\label{subsect:2.1}
\par
A compact, connected, smooth complex surface $X$ is called a \it{$K3$ surface }\rm 
if it is simply connected and has trivial canonical line bundle $\Omega^{2}_{X}$.
Let $X$ be a $K3$ surface. Then $H^{2}(X,{\bf Z})$ endowed with the cup-product pairing 
is isometric to the $K3$ lattice ${\Bbb L}_{K3}$.
An isometry of lattices $\alpha\colon H^{2}(X,{\bf Z})\cong{\Bbb L}_{K3}$ is called a {\it marking} of $X$.
The pair $(X,\alpha)$ is called a {\it marked $K3$ surface}, whose period is defined as
$$
\pi(X,\alpha):=[\alpha(\eta)]\in{\bf P}({\Bbb L}_{K3}\otimes{\bf C}),
\qquad
\eta\in H^{0}(X,\Omega^{2}_{X})\setminus\{0\}.
$$
\par
Let $M\subset{\Bbb L}_{K3}$ be a primitive $2$-elementary Lorentzian sublattice.
A $K3$ surface equipped with a holomorphic involution $\iota\colon X\to X$ is called 
a {\it $2$-elementary $K3$ surface of type $M$} if there exists a marking $\alpha$ of $X$ satisfying
$$
\iota^{*}|_{H^{0}(X,\Omega^{2}_{X})}=-1,
\qquad
\iota^{*}=\alpha^{-1}\circ I_{M}\circ\alpha.
$$
Then $\alpha(H^{2}_{+}(X,{\bf Z}))=M$, where $H^{2}_{\pm}(X,{\bf Z}):=\{l\in H^{2}(X,{\bf Z});\,\iota^{*}l=\pm l\}$.
\par
Let $(X,\iota)$ be a $2$-elementary $K3$ surface of type $M$ and let $\alpha$ be a marking with
$\theta^{*}=\alpha^{-1}\circ I_{M}\circ\alpha$. 
Since $H^{2,0}(X,{\bf C})\subset H^{2}_{-}(X,{\bf Z})\otimes{\bf C}$, 
we have $\pi(X,\alpha)\in\Omega_{M^{\perp}}^{o}$ by \cite[Th.\,3.10]{Namikawa85}. 
By \cite[Th.\,1.8]{Yoshikawa04} and Proposition~\ref{Proposition11.2} below, 
the $O(M^{\perp})$-orbit of $\pi(X,\iota)$ is independent of the choice of a marking $\alpha$
with $\iota^{*}=\alpha^{-1}I_{M}\alpha$.
The Griffiths period of $(X,\iota)$ is defined as the $O(M^{\perp})$-orbit
$$
\overline{\varpi}_{M}(X,\iota):=
O(M^{\perp})\cdot\pi(X,\alpha)\in
{\mathcal M}_{M^{\perp}}^{o}.
$$
By \cite{PSS71}, \cite{BurnsRapoport75}, \cite{Nikulin83}, \cite{Dolgachev96}, 
\cite[Th.\,1.8]{Yoshikawa04} and Proposition~\ref{Proposition11.2} below, 
the coarse moduli space of $2$-elementary $K3$ surfaces of type $M$ is isomorphic to 
${\mathcal M}_{M^{\perp}}^{o}$ via the map $\overline{\varpi}_{M}$.
In the rest of this paper, we identify the point $\overline{\varpi}_{M}(X,\iota)\in {\mathcal M}_{M^{\perp}}^{o}$ 
with the isomorphism class of $(X,\iota)$.
\par
For a $2$-elementary $K3$ surface $(X,\iota)$, set $X^{\iota}:=\{x\in X;\,\iota(x)=x\}$.

\begin{proposition}\label{Proposition2.1}
Let $(X,\iota)$ be a $2$-elementary $K3$ surface of type $M$. 
\begin{itemize}
\item[(1)]
If $M\cong{\Bbb U}(2)\oplus{\Bbb E}_{8}(2)$, then $X^{\iota}=\emptyset$.
\item[(2)]
If $M\cong{\Bbb U}\oplus{\Bbb E}_{8}(2)$, then $X^{\iota}$ is the disjoint union of two elliptic curves.
\item[(3)]
If $M\not\cong{\Bbb U}(2)\oplus{\Bbb E}_{8}(2),\,{\Bbb U}\oplus{\Bbb E}_{8}(2)$,
then there exist a smooth irreducible curve $C$ of genus $g(M)$ and smooth rational curves 
$E_{1},\ldots,E_{k(M)}$ such that $X^{\iota}=C\amalg E_{1}\amalg\cdots\amalg E_{k(M)}$.
\end{itemize}
\end{proposition}

\begin{pf}
See \cite[Th.\,4.2.2]{Nikulin83}.
\end{pf}

After Proposition~\ref{Proposition2.1}, a primitive $2$-elementary Lorentzian sublattice $M\subset{\Bbb L}_{K3}$ 
is said to be {\it non-exceptional} if $M\not\cong{\Bbb U}(2)\oplus{\Bbb E}_{8}(2),{\Bbb U}\oplus{\Bbb E}_{8}(2)$.
Let $(X,\iota)$ be a $2$-elementary $K3$ surface of type $M$.
When $M$ is non-exceptional and when $g(M)>0$,
the component of $X^{\iota}$ with genus $g(M)$ is called {\it the main component of $X^{\iota}$}.
\par
For $g\geq0$,  let ${\frak S}_{g}$ be the Siegel upper half-space of degree $g$. 
When $g=1$, ${\frak S}_{1}$ is the complex upper half-plane. We write ${\frak H}$ for ${\frak S}_{1}$.
Let $Sp_{2g}({\bf Z})$ be the symplectic group of degree $2g$ over ${\bf Z}$ and let 
${\mathcal A}_{g}:={\frak S}_{g}/Sp_{2g}({\bf Z})$ be the Siegel modular variety of degree $g$. 
Then ${\mathcal A}_{g}$ is a coarse moduli space of principally polarized Abelian varieties of dimension $g$.
The Satake compactification of ${\mathcal A}_{g}$ is denoted by ${\mathcal A}_{g}^{*}$. 
Then ${\mathcal A}_{g}^{*}$ has the stratification
${\mathcal A}_{g}^{*}={\mathcal A}_{g}\amalg{\mathcal A}_{g-1}\amalg\cdots\amalg{\mathcal A}_{0}$.
\par 
For a $2$-elementary $K3$ surface $(X,\iota)$ of type $M$, the period of $X^{\iota}$, i.e., the period of
${\rm Jac}(X^{\iota}):=H^{1}(X^{\iota},{\mathcal O}_{X^{\iota}})/H^{1}(X^{\iota},{\bf Z})$,
is denoted by $\varOmega(X^{\iota})\in{\mathcal A}_{g(M)}$.
For a $2$-elementary $K3$ surface $(X,\iota)$ of type $M$, we define
$$
\overline{J}_{M}^{o}(X,\iota)
=
\overline{J}_{M}^{o}(\overline{\varpi}_{M}(X,\iota))
:=
\varOmega(X^{\iota})
\in{\mathcal A}_{g(M)}.
$$
Let $\varPi_{M^{\perp}}\colon\Omega_{M^{\perp}}\to{\mathcal M}_{M^{\perp}}$ be the projection and set
$J_{M}^{o}:=\overline{J}_{M}^{o}\circ\varPi_{M^{\perp}}|_{\Omega_{M^{\perp}}^{o}}$.
Then $J_{M}^{o}$ is an $O(M^{\perp})$-equivariant holomorphic map from $\Omega_{M^{\perp}}^{o}$ to 
${\mathcal A}_{g(M)}$ with respect to the trivial $O(M^{\perp})$-action on ${\mathcal A}_{g(M)}$. 
By \cite[Th.\,3.3]{Yoshikawa04}, $J_{M}^{o}$ extends to an $O(M^{\perp})$-equivariant holomorphic map 
$J_{M}\colon\Omega_{M^{\perp}}^{o}\cup{\mathcal D}_{M^{\perp}}^{o}\to{\mathcal A}_{g(M)}^{*}$.
The corresponding holomorphic extension of $\overline{J}_{M}^{o}$ is denoted by 
$\overline{J}_{M}\colon{\mathcal M}_{M^{\perp}}^{o}\cup\overline{\mathcal D}_{M^{\perp}}^{o}
\to{\mathcal A}_{g(M)}^{*}$.

\begin{proposition}\label{Proposition2.2}
The map $\overline{J}_{M}$ extends to a meromorphic map from ${\mathcal M}_{M^{\perp}}^{*}$ to 
${\mathcal A}_{g(M)}^{*}$.
When $r(M)\geq19$, $\overline{J}_{M}$ extends to 
a holomorphic map from ${\mathcal M}_{M^{\perp}}^{*}$ to ${\mathcal A}_{g(M)}^{*}$.
\end{proposition}

\begin{pf}
By \cite{Borel72}, $\overline{J}_{M}$ extends to a holomorphic map from 
${\mathcal M}_{M^{\perp}}^{*}\setminus({\rm Sing}\,{\mathcal M}_{M^{\perp}}^{*}\cup
{\rm Sing}\,\overline{\mathcal D}_{M^{\perp}})$ 
to ${\mathcal A}_{g(M)}^{*}$. 
Since ${\mathcal M}_{M^{\perp}}^{*}$ is normal, we get
$\dim({\rm Sing}\,{\mathcal M}_{M^{\perp}}^{*}\cup{\rm Sing}\,\overline{\mathcal D}_{M^{\perp}})\leq
\dim{\mathcal M}_{M^{\perp}}^{*}-2$ when $r(M)\leq 18$,
so that $\overline{J}_{M}$ extends to a meromorphic map from ${\mathcal M}_{M^{\perp}}^{*}$ 
to ${\mathcal A}_{g(M)}^{*}$ by \cite{Siu75} in this case.
If $r(M)=19$, the result follows from \cite{Borel72}
because ${\mathcal M}_{M^{\perp}}^{*}$ is a compact Riemann surface and 
${\mathcal M}_{M^{\perp}}^{*}\setminus{\mathcal M}_{M^{\perp}}^{o}$ is a finite subset of 
${\mathcal M}_{M^{\perp}}^{*}$. 
If $r(M)=20$, the result is trivial because ${\mathcal M}_{M^{\perp}}^{*}$ is a point.
\end{pf}

\subsection{Degenerations of $2$-elementary $K3$ surfaces}
\label{sect:2.2}
\par
Let $\varDelta\subset{\bf C}$ be the unit disc and set $\varDelta^{*}:=\varDelta\setminus\{0\}$.
Let ${\mathcal Z}$ be a smooth complex threefold.
Let $p\colon{\mathcal Z}\to\varDelta$ be a proper surjective holomorphic function without critical points on
 ${\mathcal Z}\setminus p^{-1}(0)$.
Let $\iota\colon{\mathcal Z}\to{\mathcal Z}$ be a holomorphic involution preserving the fibers of $p$.
We set $Z_{t}=p^{-1}(t)$ and $\iota_{t}=\iota|_{Z_{t}}$ for $t\in\varDelta$.
Then $p\colon({\mathcal Z},\iota)\to\varDelta$ is called an {\it ordinary singular family}
of $2$-elementary $K3$ surfaces of type $M$ if $p$ has a unique, non-degenerate critical point on $Z_{0}$ 
and if $(Z_{t},\iota_{t})$ is a $2$-elementary $K3$ surface of type $M$ for all $t\in\varDelta^{*}$. 
Since $Z_{0}$ is a singular $K3$ surface,
$\iota_{0}\in{\rm Aut}(Z_{0})$ extends to an anti-symplectic holomorphic involution $\tilde{\iota}_{0}$ 
on the minimal resolution $\widetilde{Z}_{0}$ of $Z_{0}$, 
i.e., $(\widetilde{\iota}_{0})^{*}=-1$ on $H^{0}(\widetilde{Z}_{0},\Omega^{2}_{\widetilde{Z}_{0}})$.
Let $o\in{\mathcal Z}$ be the unique critical point of $p$.
There exists a system of coordinates $({\mathcal U},(z_{1},z_{2},z_{3}))$ centered at $o$ such that
$$
\iota(z)=(-z_{1},-z_{2},-z_{3})
\quad\hbox{ or }\quad
(z_{1},z_{2},-z_{3}),
\qquad
z\in{\mathcal U}.
$$
If $\iota(z)=(-z_{1},-z_{2},-z_{3})$ on ${\mathcal U}$, $\iota$ is said to be of {\it type $(0,3)$}. 
If $\iota(z)=(z_{1},z_{2},-z_{3})$ on ${\mathcal U}$, $\iota$ is said to be of {\it type $(2,1)$}.

\begin{theorem}\label{Theorem2.3}
Let $d\in\Delta_{M^{\perp}}$ and let $\overline{H}_{d}^{o}:=\varPi_{M^{\perp}}(H_{d}^{o})$
be the image of $H_{d}^{o}$ by the natural projection
$\varPi_{M^{\perp}}\colon\Omega_{M^{\perp}}\to{\mathcal M}_{M^{\perp}}$.
Let $\gamma\colon\varDelta\to{\mathcal M}_{M^{\perp}}$ be a holomorphic curve intersecting 
$\overline{H}_{d}^{o}$ transversally at $\gamma(0)$. 
Then there exists an ordinary singular family of $2$-elementary $K3$ surfaces
$p_{\mathcal Z}\colon({\mathcal Z},\iota)\to\varDelta$ 
of type $M$ with Griffiths period map $\gamma$ satisfying the following properties:
\begin{itemize}
\item[(1)]
$p_{\mathcal Z}$ is a projective morphism and the minimal resolution 
$(\widetilde{Z}_{0},\widetilde{\iota}_{0})$ is a $2$-elementary $K3$ surface of type $[M\perp d]$
 with Griffiths period $\gamma(0)$.
\item[(2)]
If $d\in\Delta'_{M^{\perp}}$, then $\iota$ is of type $(2,1)$ and $(\widetilde{Z}_{0})^{\widetilde{\iota}_{0}}$ 
is the normalization of $Z_{0}^{\iota_{0}}$ with total genus $g(M)-1$.
\end{itemize}
\end{theorem}

\begin{pf}
By \cite[Th.\,2.6]{Yoshikawa04}, there exists an ordinary singular family of $2$-elementary $K3$ surfaces
$p_{\mathcal Z}\colon({\mathcal Z},\iota)\to\varDelta$ of type $M$ with Griffiths period map $\gamma$ 
such that $p_{\mathcal Z}$ is projective. 
We prove that $(\widetilde{Z}_{0},\widetilde{\iota}_{0})$ is a $2$-elementary $K3$ surface of type $[M\perp d]$.
\par
Let $o_{\mathcal Z}\in Z_{0}$ be the unique critical point of $p_{\mathcal Z}$. 
Let $p_{\mathcal Y}\colon({\mathcal Y},\iota_{\mathcal Y})\to\varDelta$ be the family induced from
$p_{\mathcal Z}\colon({\mathcal Z},\iota)\to\varDelta$ by the map $\varDelta\ni t\to t^{2}\in\varDelta$.
Then ${\mathcal Y}={\mathcal Z}\times_{\varDelta}\varDelta$ and $p_{\mathcal Y}={\rm pr}_{2}$. 
The projection ${\rm pr}_{1}$ induces an identification between
$(Y_{t},\iota_{\mathcal Y}|_{Y_{t}})$ and $(Z_{t^{2}},\iota_{t^{2}})$ for all $t\in\varDelta$. 
Since the Picard-Lefschetz transformation for the family of $K3$ surfaces
$p_{\mathcal Y}|_{\varDelta^{*}}\colon{\mathcal Y}|_{\varDelta^{*}}\to\varDelta^{*}$
is trivial, there exists a marking
$\beta\colon R^{2}(p_{\mathcal Y}|_{\varDelta^{*}})_{*}{\bf Z}\cong{\Bbb L}_{K3,\varDelta^{*}}$.
Let $o_{\mathcal Y}$ be the unique singular point of ${\mathcal Y}$ 
with ${\rm pr}_{2}(o_{\mathcal Y})=o_{\mathcal Z}$. 
Since $({\mathcal Y},o_{\mathcal Y})$ is a three-dimensional ordinary double point,
there exist two different resolutions $\pi\colon({\mathcal X},E)\to({\mathcal Y},o_{\mathcal Y})$ and
$\pi'\colon({\mathcal X}',E')\to({\mathcal Y},o_{\mathcal Y})$,
which satisfy the following (i), (ii), (iii), (iv) (cf. \cite[Th.\,2.1 and Proof of Th.\,2.6]{Yoshikawa04}
and the references therein):
\begin{itemize}
\item[(i)]
Set $p:=p_{\mathcal Y}\circ\pi$ and $p':=p_{\mathcal Y}\circ\pi'$. 
Then $p\colon{\mathcal X}\to\varDelta$ and $p'\colon{\mathcal X}'\to\varDelta$
are simultaneous resolutions of $p_{\mathcal Y}\colon{\mathcal Y}\to\varDelta$,
and they are smooth families of $K3$ surfaces.
The marking $\beta$ induces a marking $\alpha$ for $p\colon{\mathcal X}\to\varDelta$
and a marking $\alpha'$ for $p'\colon{\mathcal X}'\to\varDelta$.
\item[(ii)]
$E=\pi^{-1}(o_{\mathcal Y})$ is a smooth rational curve on $X_{0}$, and
$E'=(\pi')^{-1}(o_{\mathcal Y})$ is a smooth rational curve on $X'_{0}$. 
The marked family $(p'\colon{\mathcal X}'\to\varDelta,\alpha')$ is the elementary modification of
$(p\colon{\mathcal X}\to\varDelta,\alpha)$ with center $E$ (cf. \cite[Sect.\,2.1]{Yoshikawa04}). 
Replacing $\beta$ by $g\circ\beta$, $g\in\Gamma(M):=\{g\in O({\Bbb L}_{K3});\,gI_{M}=I_{M}g\}$
if necessary, one has $d=\alpha(c_{1}([E]))$. 
\item[(iii)]
Let $e\colon{\mathcal X}\setminus E\to{\mathcal X}'\setminus E'$ be the isomorphism defined as 
$$
e:=(\pi'|_{X'\setminus E'})^{-1}\circ(\pi|_{X\setminus E}).
$$ 
Then $e$ is an isomorphism of fiber spaces over $\varDelta^{*}$ and the isomorphism
$e|_{X_{0}\setminus E}\colon X_{0}\setminus E\to X'_{0}\setminus E'$
extends to an isomorphism $\widetilde{e}_{0}\colon X_{0}\to X'_{0}$ with
\begin{equation}\label{eqn:(2.1)}
\alpha_{0}\circ(\widetilde{e}_{0})^{*}\circ
(\alpha'_{0})^{-1}=s_{d}.
\end{equation}
\item[(iv)]
There exists an isomorphism $\varphi_{K3}(I_{M})\colon{\mathcal X}\to{\mathcal X}'$ of fiber spaces 
over $\varDelta$ such that the following diagrams are commutative (cf. \cite[Eqs.\,(1.6), (2.8)]{Yoshikawa04}):
\begin{equation}\label{eqn:(2.2)}
\begin{CD}
({\mathcal X},E)
@>\pi>>
({\mathcal Y},o)
@>{\rm pr}_{1}>>
({\mathcal Z},o)
\\
@V \varphi_{K3}(I_{M}) VV @V \iota_{\mathcal Y}VV
@V \iota VV
\\
({\mathcal X}',E')
@>\pi'>>
({\mathcal Y},o)
@>{\rm pr}_{1}>>
({\mathcal Z},o)
\end{CD}
\qquad\qquad\qquad
\begin{CD}
R^{2}p'_{*}{\bf Z}
@>\varphi_{K3}(I_{M})^{*}>>
R^{2}p_{*}{\bf Z}
\\
@V \alpha' VV @VV \alpha V
\\
{\Bbb L}_{K3,\varDelta}
@>I_{M}>>
{\Bbb L}_{K3,\varDelta}
\end{CD}
\end{equation}
\end{itemize}
\par
We define $\theta:=(\widetilde{e}_{0})^{-1}\circ\varphi_{K3}(I_{M})|_{X_{0}}\in{\rm Aut}(X_{0})$. 
Since $\pi'\circ\widetilde{e}_{0}=\pi|_{X_{0}}$ by (iii) and hence 
$\pi'|_{X'_{0}\setminus E'}=(\pi|_{X_{0}\setminus E})\circ(\widetilde{e}_{0})^{-1}|_{X'_{0}\setminus E'}$,
we get by the first diagram of \eqref{eqn:(2.2)}
$$
\begin{aligned}
(\pi|_{X_{0}\setminus E})\circ
(\theta|_{X_{0}\setminus E})
&=
(\pi|_{X_{0}\setminus E})\circ
(\widetilde{e}_{0})^{-1}|_{X'_{0}\setminus E'}\circ
\varphi_{K3}(I_{M})|_{X_{0}\setminus E}
\\
&=
(\pi'|_{X'_{0}\setminus E'})\circ
\varphi_{K3}(I_{M})|_{X_{0}\setminus E}
\\
&=
(\iota_{\mathcal Y}|_{Y_{0}\setminus\{o\}})\circ
\pi|_{X_{0}\setminus E},
\end{aligned}
$$
which implies $(\pi|_{X_{0}})\circ\theta=(\iota_{\mathcal Y}|_{Y_{0}})\circ(\pi|_{X_{0}})$.
Since $X_{0}$ is the minimal resolution of $Z_{0}$, i.e., $X_{0}=\widetilde{Z}_{0}$ and since
$(Y_{0},\iota_{\mathcal Y}|_{Y_{0}})=(Z_{0},\iota_{0})$, 
the equality $(\pi|_{X_{0}})\circ\theta=(\iota_{\mathcal Y}|_{Y_{0}})\circ(\pi|_{X_{0}})$ implies that 
$\theta$ is the involution on $X_{0}$ induced from $\iota_{0}$. Thus $\theta=\widetilde{\iota}_{0}$. 
\par
By \eqref{eqn:(1.1)}, \eqref{eqn:(2.1)} and the second diagram of \eqref{eqn:(2.2)}, we get
\begin{equation}
\label{eqn:(2.3)}
\alpha_{0}\,\theta^{*}\,\alpha_{0}^{-1}
=
\alpha_{0}\varphi_{K3}(I_{M})^{*}(\alpha'_{0})^{-1}
\circ\alpha'_{0}(\widetilde{e}_{0}^{-1})^{*}\alpha_{0}^{-1}
=I_{M}\circ s_{d}=I_{[M\perp d]}.
\end{equation}
By \eqref{eqn:(2.3)}, $\theta=\widetilde{\iota}_{0}$ is an anti-symplectic involution of type $[M\perp d]$. 
This proves (1).
\par
Let $d\in\Delta'_{M^{\perp}}$. If $\iota$ is of type $(0,3)$,
then $g([M\perp d])=g(M)$ by \cite[Prop.\,2.5]{Yoshikawa04}.
Since $d\in\Delta'_{M^{\perp}}$ implies $g([M\perp d])=g(M)-1$ by Lemma~\ref{Lemma11.5} below, 
we get a contradiction. Hence $\iota$ must be of type $(2,1)$.
Since $(\widetilde{Z}_{0},\widetilde{\iota}_{0})$ is a $2$-elementary $K3$ surface of type $[M\perp d]$,
$(\widetilde{Z}_{0})^{\widetilde{\iota}_{0}}$ has total genus $g([M\perp d])=g(M)-1$ by Lemma~\ref{Lemma11.5}. 
Since $\widetilde{Z}_{0}\to Z_{0}$ is the blow-up at the ordinary double point $o_{\mathcal Z}$,
it follows from the local description $\iota(z)=(z_{1},z_{2},-z_{3})$ near $o_{\mathcal Z}$
that the set of fixed points $(\widetilde{Z}_{0})^{\widetilde{\iota}_{0}}$ is the normalization of $Z_{0}^{\iota_{0}}$.
This proves (2).
\end{pf}

Let ${\mathcal C}$ be a (possibly disconnected) smooth complex surface. 
Let $p\colon{\mathcal C}\to\varDelta$ be a proper, surjective holomorphic function without critical points 
on ${\mathcal C}\setminus p^{-1}(0)$.
Then $p\colon{\mathcal C}\to\varDelta$ is called an ordinary singular family of curves if $p$ has a unique, 
non-degenerate critical point on $p^{-1}(0)$. We set $C_{t}:=p^{-1}(t)$ for $t\in\varDelta$.

\begin{lemma}\label{Lemma2.4}
Let $p\colon{\mathcal C}\to\varDelta$ be an ordinary singular family of curves and let
$g=\dim H^{0}(C_{t},\Omega^{1}_{C_{t}})$ for $t\not=0$. 
Let $J\colon\varDelta^{*}\to{\mathcal A}_{g}$ be the holomorphic map defined as
$J(t):=\varOmega(C_{t})$ for $t\in\varDelta^{*}$. 
Then $J$ extends to a holomorphic map from $\varDelta$ to ${\mathcal A}_{g}^{*}$ by setting 
$J(0):=\varOmega(\widehat{C_{0}})$, where $\widehat{C}_{0}$ is the normalization of $C_{0}$.
\end{lemma}

\begin{pf}
Since the result is obvious when $g=0$, we assume $g>0$.
The extendability of $J$ follows from e.g. \cite[Chap.\,III Th.\,16.1]{BPV84}.
Assume that $p$ has connected fibers.
Either $C_{0}$ is the join of two smooth curves $A$ and $B$ intersecting transversally at 
${\rm Sing}\,C_{0}$ or $C_{0}$ is irreducible. The result follows from e.g. \cite[Cors.\,3.2 and 3.8]{Fay73}.
\par
Assume that ${\mathcal C}$ is not connected. 
Let ${\mathcal C}={\mathcal C}_{0}\amalg\ldots\amalg{\mathcal C}_{k}$ be the decomposition 
into the connected components and set $p_{i}:=p|_{{\mathcal C}_{i}}$. 
Since the period matrix of $C_{t}$ is the direct sum of those of the curves $p_{i}^{-1}(t)$, 
the result follows from the case where $p$ has connected fibers and \cite[Chap.\,III Th.\,(16.1)]{BPV84}.
\end{pf}

\begin{theorem}\label{Theorem2.5}
For $d\in\Delta_{M^{\perp}}$, the following equality holds
$$
J_{M}|_{H_{d}^{o}}=J_{[M\perp d]}^{o}|_{H_{d}^{o}}.
$$
\end{theorem}

\begin{pf}
Let ${\frak p}\in\overline{H}_{d}^{o}$ and let $\gamma\colon\varDelta\to{\mathcal M}_{M^{\perp}}$ 
be a holomorphic curve intersecting $\overline{H}_{d}^{o}$ transversally at ${\frak p}=\gamma(0)$.
Let $p_{\mathcal Z}\colon({\mathcal Z},\iota)\to\varDelta$ be an ordinary singular family of 
$2$-elementary $K3$ surfaces of type $M$ with Griffiths period map $\gamma$, 
such that $p_{\mathcal Z}$ is projective and such that $(\widetilde{Z}_{0},\widetilde{\iota}_{0})$ 
is a $2$-elementary $K3$ surface of type $[M\perp d]$ with Griffiths period $\gamma(0)$ 
(cf. Theorem~\ref{Theorem2.3}).
Let $o\in{\mathcal Z}$ be the unique critical point of $p_{\mathcal Z}$.
Since $\overline{J}_{M}({\frak p})=\overline{J}_{M}(\gamma(0))=\lim_{t\to0}\overline{J}_{M}(\gamma(t))$ 
by the continuity of $\overline{J}_{M}$ and since
$\overline{J}_{[M\perp d]}^{o}({\frak p})=\overline{J}_{[M\perp d]}^{o}(\widetilde{Z}_{0},\widetilde{\iota}_{0})=
\varOmega(({\widetilde Z}_{0})^{\widetilde{\iota}_{0}})$ 
by Theorem~\ref{Theorem2.3}, it suffices to prove
\begin{equation}\label{eqn:(2.4)}
\overline{J}_{M}({\frak p})
=
\lim_{t\to0}\overline{J}_{M}(\gamma(t))
=
\varOmega(({\widetilde Z}_{0})^{\widetilde{\iota}_{0}})
=
\overline{J}_{[M\perp d]}^{o}({\frak p}).
\end{equation}
Set ${\mathcal Z}^{\iota}:=\{z\in{\mathcal Z};\,\iota(z)=z\}$.
\par
Assume that $\iota$ is of type $(0,3)$. 
By \cite[Prop.\,2.5 (1)]{Yoshikawa04}, ${\mathcal C}:={\mathcal Z}^{\iota}\setminus\{o\}$ 
is a smooth complex surface and $p|_{\mathcal C}\colon{\mathcal C}\to\varDelta$ 
is a proper holomorphic submersion. Then
\begin{equation}\label{eqn:(2.5)}
\lim_{t\to0}\overline{J}_{M}^{o}(Z_{t},\iota_{t})
=
\lim_{t\to0}\varOmega(C_{t})
=
\varOmega(C_{0}).
\end{equation}
Since $Z_{0}^{\iota_{0}}=C_{0}\amalg\{o\}$, we get
$(\widetilde{Z}_{0})^{\widetilde{\iota}_{0}}=C_{0}\amalg{\bf P}^{1}$,
which yields that
\begin{equation}\label{eqn:(2.6)}
\varOmega(C_{0})
=
\varOmega(({\widetilde Z}_{0})^{\widetilde{\iota}_{0}}).
\end{equation}
Eq.\,\eqref{eqn:(2.4)} follows from \eqref{eqn:(2.5)} and \eqref{eqn:(2.6)} in this case.
\par
Assume that $\iota$ is of type $(2,1)$. By \cite[Prop.\,2.5 (2)]{Yoshikawa04},  
$p|_{{\mathcal Z}^{\iota}}\colon{\mathcal Z}^{\iota}\to\varDelta$ is an ordinary singular family of curves.
Let $W\to Z_{0}^{\iota_{0}}$ be the normalization. We get
\begin{equation}\label{eqn:(2.7)}
\lim_{t\to0}\overline{J}_{M}^{o}(Z_{t},\iota_{t})
=
\lim_{t\to0}
\varOmega(Z_{t}^{\iota_{t}})
=
\varOmega(W)\in{\mathcal A}_{g(M)}^{*}
\end{equation}
by Lemma~\ref{Lemma2.4}. In the same manner as in the proof of Theorem~\ref{Theorem2.3} (2), 
we get $W=(\widetilde{Z}_{0})^{\widetilde{\iota}_{0}}$, which together with \eqref{eqn:(2.7)}, 
yields \eqref{eqn:(2.4)} in this case.
Since ${\frak p}$ is an arbitrary point of $\overline{H}_{d}^{o}$, we get the result.
\end{pf}

The following propositions shall be used in the proof of Proposition~\ref{Proposition4.2} (3).

\begin{proposition}\label{Proposition2.6}
If $g(M)=1$ and $d\in\Delta'_{M^{\perp}}$, then
$J_{M}(H_{d}^{o})={\mathcal A}_{0}={\mathcal A}_{1}^{*}\setminus{\mathcal A}_{1}$.
\end{proposition}

\begin{pf}
By Lemma~\ref{Lemma11.5} below, $g([M\perp d])=g(M)-1=0$.
By Theorem~\ref{Theorem2.5}, we get
$J_{M}(H_{d}^{o})=J_{[M\perp d]}^{o}(H_{d}^{o})={\mathcal A}_{0}
={\mathcal A}_{1}^{*}\setminus{\mathcal A}_{1}$.
\end{pf}

\begin{proposition}\label{Proposition2.7}
If $g(M)=1$, then $\overline{J_{M}^{o}(\Omega_{M^{\perp}}^{o})}={\mathcal A}_{1}^{*}$.
\end{proposition}

\begin{pf}
By Proposition~\ref{Proposition2.2}, $\overline{J}_{M}$ extends to a meromorphic map from 
${\mathcal M}_{M^{\perp}}^{*}$ to ${\mathcal A}_{1}^{*}$. 
Since
$J_{M}^{o}(\Omega_{M^{\perp}}^{o})=\overline{J}_{M}({\mathcal M}_{M^{\perp}}^{o})$ and 
since $\dim{\mathcal A}_{1}^{*}=1$, we have
$\overline{J_{M}^{o}(\Omega_{M^{\perp}}^{o})}={\mathcal A}_{1}^{*}$
if $\overline{J}_{M}^{o}$ is non-constant. 
We see that $\overline{J}^{o}_{M}$ is non-constant.
\par
Since $g(M)=1$, we get by \cite[p.1434, Table 1]{Nikulin83}
or \cite[Appendix, Table 2]{FinashinKharlamov08}
\begin{equation}\label{eqn:(2.8)}
M^{\perp}
\cong
{\Bbb U}\oplus{\Bbb A}_{1}^{+}\oplus{\Bbb A}_{1}^{\oplus m-1}
\quad
(1\leq m\leq 10),
\quad
{\Bbb U}(2)\oplus{\Bbb U}(2)\oplus{\Bbb D}_{4},
\quad
{\Bbb U}\oplus{\Bbb U}(2).
\end{equation}
By \eqref{eqn:(2.8)}, $\Delta'_{M^{\perp}}\not=\emptyset$. Let $d\in\Delta'_{M^{\perp}}$.
By Proposition~\ref{Proposition2.6}, we get
$J_{M}(H_{d}^{o})={\mathcal A}_{0}={\mathcal A}_{1}^{*}\setminus{\mathcal A}_{1}$.
Since $J_{M}(\Omega_{M^{\perp}}^{o})\subset{\mathcal A}_{1}$,
this implies that $\overline{J}_{M}$ is non-constant.
\end{pf}

\begin{proposition}\label{Proposition2.8}
If $g(M)=1$ and $d\in\Delta''_{M^{\perp}}$, then $J_{M}(H_{d}^{o})\subset{\mathcal A}_{1}$.
\end{proposition}

\begin{pf}
Since $d\in\Delta''_{M^{\perp}}$, we get $g([M\perp d])=g(M)=1$ by Lemma~\ref{Lemma11.5} below.
By Theorem~\ref{Theorem2.5}, we get
$J_{M}(H_{d}^{o})=J_{[M\perp d]}^{o}(H_{d}^{o})
\subset 
J_{[M\perp d]}^{o}(\Omega_{[M\perp d]^{\perp}}^{o})
\subset
{\mathcal A}_{1}$.
\end{pf}

\begin{proposition}\label{Proposition2.9}
If $g(M)=2$ and $d\in\Delta'_{M^{\perp}}$, then
$\overline{J_{M}(H_{d}^{o})}={\mathcal A}_{2}^{*}\setminus{\mathcal A}_{2}$.
\end{proposition}

\begin{pf}
By Lemma~\ref{Lemma11.5} below, $g([M\perp d])=1$.
By Theorem~\ref{Theorem2.5}, we get
$$
\overline{J_{M}(H_{d}^{o})}
=
\overline{J_{[M\perp d]}^{o}(H_{d}^{o})}
=
\overline{J_{[M\perp d]}(\Omega_{[M\perp d]^{\perp}}^{o})}
=
{\mathcal A}_{1}^{*}
=
{\mathcal A}_{2}^{*}\setminus{\mathcal A}_{2},
$$
where the third equality follows from Proposition~\ref{Proposition2.7}.
\end{pf}

We define the divisor 
${\mathcal N}_{2}\subset{\mathcal A}_{2}$ as
$$
{\mathcal N}_{2}
:=
\{\varOmega(E_{1}\times E_{2})\in{\mathcal A}_{2};\,E_{1},E_{2}\hbox{ are elliptic curves}\}.
$$

\begin{proposition}\label{Proposition2.10}
Let $g(M)=2$ and $d\in\Delta''_{M^{\perp}}$. 
Then $J_{M}(H_{d}^{o})\cap{\mathcal N}_{2}\not=\emptyset$ if and only if 
$M\cong{\Bbb A}_{1}^{+}\oplus{\Bbb A}_{1}^{\oplus8}$ and $d/2\equiv{\bf 1}_{M^{\perp}}\mod M^{\perp}$.
In particular, if either $M\not\cong{\Bbb A}_{1}^{+}\oplus{\Bbb A}_{1}^{\oplus8}$ 
or $d/2\not\equiv{\bf 1}_{M^{\perp}}\mod M^{\perp}$, then
$$
J_{M}(H_{d}^{o})\subset{\mathcal A}_{2}\setminus{\mathcal N}_{2}.
$$
\end{proposition}

\begin{pf}
Assume $J_{M}(H_{d}^{o})\cap{\mathcal N}_{2}\not=\emptyset$. By Theorem~\ref{Theorem2.5},
$$
J_{[M\perp d]}(\Omega_{[M\perp d]^{\perp}}^{o})\cap{\mathcal N}_{2}
\supset
J_{[M\perp d]}(H_{d}^{o})\cap{\mathcal N}_{2}
=
J_{M}(H_{d}^{o})\cap{\mathcal N}_{2}\not=\emptyset.
$$
Let $(X,\iota)$ be a $2$-elementary $K3$ surface of type $[M\perp d]$ such that
$J_{[M\perp d]}(X,\iota)\in{\mathcal N}_{2}$.
If $[M\perp d]\not\cong{\Bbb U}\oplus{\Bbb E}_{8}(2),{\Bbb U}(2)\oplus{\Bbb E}_{8}(2)$, 
there exists an irreducible smooth curve $C$ of genus $g([M\perp d])$ with 
$J_{[M\perp d]}(X,\iota)=\varOmega(C)$ by Proposition~\ref{Proposition2.1}.
By $d\in\Delta''_{M^{\perp}}$ and Lemma~\ref{Lemma11.5} below,
we get $g([M\perp d])=2$. However, the period of an irreducible smooth curve of genus $2$ lies in 
${\mathcal A}_{2}\setminus{\mathcal N}_{2}$. 
This contradicts the condition $\varOmega(C)\in{\mathcal N}_{2}$.
Thus $[M\perp d]\cong{\Bbb U}\oplus{\Bbb E}_{8}(2)$ or
$[M\perp d]\cong{\Bbb U}(2)\oplus{\Bbb E}_{8}(2)$.
If $[M\perp d]\cong{\Bbb U}(2)\oplus{\Bbb E}_{8}(2)$, then 
$C=\emptyset$ by Proposition~\ref{Proposition2.1} (1), which contradicts
the condition $\varOmega(C)\in{\mathcal N}_{2}$.
We get $[M\perp d]\cong{\Bbb U}\oplus{\Bbb E}_{8}(2)$ and hence 
$M^{\perp}\cap d^{\perp}=[M\perp d]^{\perp}\cong{\Bbb U}^{\oplus2}\oplus{\Bbb E}_{8}(2)$.
Since $d\in\Delta''_{M^{\perp}}$, we get $(r(M),l(M))=(r([M\perp d])-1,l([M\perp d])+1)=(9,9)$
by Proposition~\ref{Proposition11.6} below. Since $r(M)=9$, we get $\delta(M)=1$. 
All together, we get $(r(M),l(M),\delta(M))=(9,9,1)$ and hence
$M\cong{\Bbb A}_{1}^{+}\oplus{\Bbb A}_{1}^{\oplus8}$.
Set $L:={\bf Z}\,d\subset M^{\perp}$. Then $L\cong{\Bbb A}_{1}$.
Since $(r(M^{\perp}),l(M^{\perp}),\delta(M^{\perp}))=(13,9,1)$, we get the decomposition
$M^{\perp}\cong(M^{\perp}\cap d^{\perp})\oplus L\cong{\Bbb U}^{\oplus2}\oplus{\Bbb E}_{8}(2)\oplus L$
by comparing the triplets $(r,l,\delta)$, which implies
${\bf 1}_{M^{\perp}}={\bf 1}_{M^{\perp}\cap d^{\perp}}\oplus{\bf 1}_{L}={\bf 1}_{L}=d/2$ in $A_{M^{\perp}}$.
\par
Conversely, assume $M\cong{\Bbb A}_{1}^{+}\oplus{\Bbb A}_{1}^{\oplus8}$ and
$d/2\equiv{\bf 1}_{M^{\perp}}\mod M^{\perp}$. 
Since $d/2\in(M^{\perp})^{\lor}$, we get $d\in\Delta''_{M^{\perp}}$.
Since $(r(M^{\perp}),l(M^{\perp}))=(r(M^{\perp}\cap d^{\perp})+1,l(M^{\perp}\cap d^{\perp})+1)$
by Proposition~\ref{Proposition11.6} (3), (4) below, 
we get $r(M^{\perp}\cap d^{\perp})=r(M^{\perp})-1=12$ and $l(M^{\perp}\cap d^{\perp})=l(M^{\perp})-1=8$. 
Since $\delta(M^{\perp})=1$, we get $M^{\perp}=(M^{\perp}\cap d^{\perp})\oplus L$ by comparing $(r,l,\delta)$.
Let us see that $\delta(M^{\perp}\cap d^{\perp})=0$.
Let $x\in(M^{\perp}\cap d^{\perp})^{\lor}$ be an arbitrary element and let
$k\in{\bf Z}$. Set $y=x+k(d/2)\in(M^{\perp})^{\lor}=(M^{\perp}\cap d^{\perp})^{\lor}\oplus L^{\lor}$.
Since ${\bf 1}_{M^{\perp}}\equiv d/2\mod M^{\perp}$, we get by the definition of ${\bf 1}_{M^{\perp}}$
$$
-k/2
=
\langle y,d/2\rangle
\equiv
\langle y,{\bf 1}_{M^{\perp}}\rangle
\equiv
\langle y,y\rangle
\equiv
\langle x,x\rangle-k^{2}/2
\mod{\bf Z}.
$$
Hence $x^{2}\equiv k(k-1)/2\equiv0\mod{\bf Z}$, which implies $\delta(M^{\perp}\cap d^{\perp})=0$.
Since $(r,l,\delta)=(12,8,0)$ for $M^{\perp}\cap d^{\perp}$, 
we get $M^{\perp}\cap d^{\perp}\cong{\Bbb U}^{\oplus2}\oplus{\Bbb E}_{8}(2)$ and hence
$[M\perp d]\cong{\Bbb U}\oplus{\Bbb E}_{8}(2)$.
By Theorem~\ref{Theorem2.5},
$J_{M}(H_{d}^{o})\subset J_{[M\perp d]}^{o}(\Omega_{[M\perp d]}^{o})\subset{\mathcal N}_{2}$, 
where the last inclusion follows from Proposition~\ref{Proposition2.1} (2). 
Hence $J_{M}(H_{d}^{o})\cap{\mathcal N}_{2}\not=\emptyset$.
\end{pf}

\section{Log del Pezzo surfaces and $2$-elementary $K3$ surfaces}
\label{sect:3}
\par
In this section, we recall the notion of log del Pezzo surfaces of index $\leq2$ and DNP surfaces, 
for which we refer the reader to \cite{AlexeevNikulin06} and \cite{Nakayama07}.
In Sect.\,\ref{sect:3}, the canonical divisor of a normal complex surface $S$ is denoted by $K_{S}$.
Hence the canonical line bundle of $S$ is denoted by ${\mathcal O}_{S}(K_{S})$ in stead of $\Omega^{2}_{S}$ 
in this section.

\subsection{Log del Pezzo surfaces of index $2$ and DNP surfaces}
\label{subsect:3.1}
\par
A normal projective surface $S$ is a {\it log del Pezzo surface} 
if it has only log terminal singularities and 
if its anticanonical divisor $-K_{S}$ is an ample ${\bf Q}$-Cartier divisor.
The {\it index} of $S$ is the smallest integer $\nu\in{\bf Z}_{>0}$ such that $-\nu K_{S}$ is Cartier
\cite[Sect.\,1]{AlexeevNikulin06}.
\par
A smooth projective surface $Y$ is a {\it DNP surface}
if $h^{1}(Y)=0$, $K_{Y}\not\sim0$
and if there exists an effective divisor $C\in|-2K_{Y}|$
with only simple singularities
\cite[Sect.\,2.1]{AlexeevNikulin06}.
A DNP surface $Y$ is rational if $|-2K_{Y}|\not=\emptyset$. 
If $Y$ is a DNP surface and if $C\in|-2K_{Y}|$
is a smooth divisor, the pair $(Y,C)$ is called
a right DNP pair.
\par
Let $S$ be a log del Pezzo surface of index $2$.
By \cite[Th.\,1.5]{AlexeevNikulin06}, $|-2K_{S}|$ contains
a smooth curve. Let $C\in|-2K_{S}|$ be smooth. 
To the pair $(S,C)$, one can associate a right DNP pair and
a $2$-elementary $K3$ surface as follows
\cite[Sect.\,2.1]{AlexeevNikulin06},
\cite[Sect.\,6.6]{Nakayama07}.
\par
Let $\alpha\colon\widetilde{S}\to S$ be the minimal resolution. 
Since $S$ has only log terminal singularities of index $2$, 
we deduce from \cite[Sect.\,1.2]{AlexeevNikulin06} 
the existence of a non-zero $\alpha$-exceptional simple 
normal crossing divisor $E$ on $\widetilde{S}$ such that
$-2K_{\widetilde{S}}\sim\alpha^{*}(-2K_{S})+E$.
If $D$ is a connected component of $E$, the germ
$(S,\alpha(D))\in{\rm Sing}\,S$ is isomorphic to
one of the singularities $K_{n}$ 
in \cite[Sect.\,1.2]{AlexeevNikulin06},
\cite[Example 4.17]{Nakayama07}.
\par
Let $\beta\colon Y\to\widetilde{S}$ be the blowing-up
at the nodes of $E$. By \cite[Sect.\,1.2]{AlexeevNikulin06},
the proper transform $E_{Y}$ of $E$ is the disjoint union
of $(-4)$-curves on $Y$ and the total transform
$\beta^{*}E$ is the disjoint union of the configurations
in \cite[Sect.\,1.5 (9)]{AlexeevNikulin06}. 
Set $p:=\alpha\circ\beta$. The birational morphism 
$p\colon Y\to S$ is called the {\it right resolution} of $S$.
\par
Let $C_{Y}:=p^{-1}(C)\subset Y$ be the total transform of $C$
with respect to the birational morphism $p\colon Y\to S$.
Since $C\in|-2K_{S}|$ is smooth and hence
$C\cap{\rm Sing}\,S=\emptyset$,
$p|_{C_{Y}}\colon C_{Y}\to C$ is an isomorphism.
By \cite[Sect.\,2.1]{AlexeevNikulin06}, 
\cite[p.415 Eq.\,(6-1)]{Nakayama07}, $Y$ is a DNP surface
and the pair $(Y,C_{Y}+E_{Y})$ is a right DNP pair.
We call $(Y,C_{Y}+E_{Y})$ 
{\it the right DNP pair associated to $(S,C)$}.
\par
Since $C_{Y}+E_{Y}\in|-2K_{Y}|$, there exists a double
covering $\pi\colon X\to Y$ with branch divisor $C_{Y}+E_{Y}$.
Let $\iota\colon X\to X$ be the non-trivial covering
transformation of $\pi\colon X\to Y$. 
By \cite[Sect.\,2.1]{AlexeevNikulin06}, 
\cite[Sect.\,6.6]{Nakayama07},
$(X,\iota)$ is a $2$-elementary $K3$ surface such that
$X^{\iota}\cong C_{Y}+E_{Y}$. We call $(X,\iota)$
{\it the $2$-elementary $K3$ surface associated to $(S,C)$}.
In this case, $g(C)=g(X^{\iota})\geq2$
by \cite[Th.\,4.1]{AlexeevNikulin06}.
\par
Conversely, if $(X,\iota)$ is a $2$-elementary $K3$ surface
with $g(X^{\iota})\geq2$, then $(X/\iota,X^{\iota})$
is a right DNP pair. By \cite[Th.\,4.1]{AlexeevNikulin06}, 
there exists a unique pair $(S,C)$, 
where $S$ is a log del Pezzo surface of index $\leq2$ and 
$C\in|-2K_{S}|$ is a smooth member, 
such that $(X,\iota)$ is associated to $(S,C)$.

\subsection{Some properties of the main component of $X^{\iota}$}
\label{subsect:3.2}
\par
Let $M\subset{\Bbb L}_{K3}$ be a primitive $2$-elementary
Lorentzian sublattice.
Assume that $M$ is non-exceptional and that $g(M)\geq1$. 
Recall that if $(X,\iota)$ is a $2$-elementary $K3$ surface
of type $M$ and if
$X^{\iota}=C\amalg E_{1}\amalg\ldots\amalg E_{k}$
denotes the decomposition into the connected components
with $g(C)=g(M)\geq1$ and $E_{i}\cong{\bf P}^{1}$
for $1\leq i\leq k$, then $C$ is called the {\it main component} 
of $X^{\iota}$.

\begin{proposition}\label{Proposition3.1}
Assume that $r(M)>10$ or $(r(M),\delta(M))=(10,1)$.
Then $0\leq g(M)\leq 5$ and the following hold:
\begin{itemize}
\item[(1)]
If $g(M)=3$, then there exists a $2$-elementary $K3$ surface 
$(X,\iota)$ of type $M$ such that the main component 
of $X^{\iota}$ is non-hyperelliptic.
\item[(2)]
If $g(M)=4$, then there exists a $2$-elementary $K3$ surface 
$(X,\iota)$ of type $M$ such that the main component 
of $X^{\iota}$ is isomorphic to the complete intersection 
of a smooth quadric and a (possibly singular) cubic 
in ${\bf P}^{3}$.
\item[(3)]
If $g(M)=5$, then there exists a $2$-elementary $K3$ surface 
$(X,\iota)$ of type $M$ such that the main component 
of $X^{\iota}$ is the normalization of an irreducible 
plane quintic with one node.
\end{itemize}
\end{proposition}

\begin{pf}
By \cite[p.1434, Table 1]{Nikulin83} and 
the assumption on $M$, we get $0\leq g(M)\leq5$.
\par
In what follows, we use Nakayama's notation
\cite[p.410, Table 6]{Nakayama07} and 
\cite[p.494-495, Table 14]{Nakayama07}
for the type of log del Pezzo surfaces of index $\leq2$.
See \cite[p.410, Tables 9]{Nakayama07} 
for the relation between the type of 
a log del Pezzo surface and the type of the associated
$2$-elementary $K3$ surface.
\par
Let $S$ be a log del Pezzo surface and
let $\Gamma\in|-2K_{S}|$ be a non-singular member.
Let $M\subset{\Bbb L}_{K3}$ be the type of 
the $2$-elementary $K3$ surface $(X,\iota)$ 
associated to $(S,\Gamma)$.
The main component of $X^{\iota}$ is isomorphic to $\Gamma$ 
by construction. 
\newline{\bf (1) }
Since $g(M)=3$, $r(M)\geq10$ and 
$(r(M),\delta(M))\not=(10,0)$,
the type of $S$ is one of $[2]_{+}(b)$ $(0\leq b\leq 4)$ 
by \cite[p.410, Table 6 and p.444 Table 9]{Nakayama07}. 
By \cite[pp.494-495, Table 14]{Nakayama07}, 
$S$ is a hypersurface of ${\bf P}(1,1,1,2)$ 
defined by the equation:
$$
\begin{array}{lll}
\hbox{Lattice }M^{\perp}
&
\hbox{Type of }S
&
\hbox{Equation defining }S
\\
{\Bbb U}^{\oplus2}\oplus{\Bbb D}_{4}\oplus{\Bbb A}_{1}^{\oplus 4}
&
[2]_{+}(0)
&
xyu=z^{4}+F_{3}(x,z)x+G_{3}(y,z)y,
\\
{\Bbb U}^{\oplus2}\oplus{\Bbb D}_{4}\oplus{\Bbb A}_{1}^{\oplus 4-b}
&
[2]_{+}(b)\,\,(1\leq b\leq4)
&
xyu=F_{4-b}(x,z)x^{b}+G_{3}(y,z)y,
\end{array}
$$
where 
${\rm wt}(x)={\rm wt}(y)={\rm wt}(z)=1$, ${\rm wt}(u)=2$
and $F_{k}(x,y)\in{\bf C}[x,y]$ is a homogeneous 
polynomial of degree $k$ and $G_{3}(x,y)\in{\bf C}[x,y]$ 
is a homogeneous polynomial of degree $3$.
Since 
${\mathcal O}_{S}(-2K_{S})\cong
{\mathcal O}_{{\bf P}(1,1,1,2)}(2)|_{S}$
by the adjunction formula and hence 
$u-U(x,y,z)\in H^{0}(S,{\mathcal O}_{S}(-2K_{S}))$
for a homogeneous polynomial $U(x,y,z)\in{\bf C}[x,y,z]$
of degree $2$,
a general member of the linear system $|-2K_{S}|$ 
is a hypersurface of ${\bf P}^{2}$ 
defined by the equation:
$$
\begin{array}{ll}
\hbox{Type of }S
&
\hbox{Equation defining }\Gamma
\\
\,
[2]_{+}(0)
&
xy\,U(x,y,z)=z^{4}+F_{3}(x,z)x+G_{3}(y,z)y
\\
\,
[2]_{+}(b)\,\,(1\leq b\leq4)
&
xy\,U(x,y,z)=F_{4-b}(x,z)x^{b}+G_{3}(y,z)y.
\end{array}
$$
In particular,
$\Gamma\in|-2K_{S}|$ is a smooth plane quartic
if $F_{k}(x,y)$, $G_{3}(x,y)$ and $U(x,y,z)$ are 
sufficiently general. Since a smooth plane quartic curve
is non-hyperelliptic, we get (1).
\newline{\bf (2) }
Since $g(M)=4$, $r(M)\geq10$ and 
$(r(M),\delta(M))\not=(10,0)$,
the type of $S$ is one of $[0;1,1]_{+}(b)$ 
$(1\leq b\leq 3)$ 
by \cite[p.410, Table 6 and p.444, Table 9]{Nakayama07}. 
By \cite[pp.494-495, Table 14]{Nakayama07}, 
$S$ is a complete intersection of ${\bf P}(1,1,1,1,2)$ 
defined by the equations:
$$
\begin{array}{lll}
\hbox{Lattice }M^{\perp}
&
\hbox{Type of }S
&
\hbox{Equations defining }S
\\
({\Bbb A}_{1}^{+})^{\oplus2}\oplus{\Bbb E}_{8}\oplus
{\Bbb A}_{1}^{\oplus 2}
&
[0;1,1]_{+}(1)
&
\begin{cases}
xw=yz,
\\
xu=(w+c\,z)zw+(w+c'\,y)yw
\end{cases}
\\
({\Bbb A}_{1}^{+})^{\oplus2}\oplus{\Bbb E}_{8}\oplus
{\Bbb A}_{1}^{\oplus 3-b}
&
[0;1,1]_{+}(b)\,\,(b\geq2)
&
\begin{cases}
xw=yz,
\\
xu=(w+c\,z)zw+w^{3-b}y^{b}
\end{cases}
\end{array}
$$
where 
${\rm wt}(x)={\rm wt}(y)={\rm wt}(z)={\rm wt}(w)=1$,
${\rm wt}(u)=2$ and $c,c'\in{\bf C}$ are constants.
Since 
${\mathcal O}_{S}(-2K_{S})\cong
{\mathcal O}_{{\bf P}(1,1,1,1,2)}(2)|_{S}$ 
by the adjunction formula or by Lemma~\ref{Lemma3.2} below
and hence 
$u-U(x,y,z,w)\in H^{0}(S,{\mathcal O}_{S}(-2K_{S}))$
for a homogeneous polynomial 
$U(x,y,z,w)\in{\bf C}[x,y,z,w]$ of degree $2$,
$\Gamma\in|-2K_{S}|$ is a complete intersection of
${\bf P}^{3}$ defined by the equations:
$$
\begin{array}{ll}
\hbox{Type of }S
&
\hbox{Equations defining }\Gamma
\\
\,
[0;1,1]_{+}(1)
&
xw=yz,
\quad 
x\,U(x,y,z,w)=(w+c\,z)zw+(w+c'\,y)yw
\\
\,
[0;1,1]_{+}(b)\,\,(2\leq b\leq3)
&
xw=yz,
\quad 
x\,U(x,y,z,w)=(w+c\,z)zw+w^{3-b}y^{b}.
\end{array}
$$
By choosing $c$, $c'$, $U(x,y,z,w)$ sufficiently general, 
$\Gamma$ is a complete intersection 
in ${\bf P}^{3}$ of a smooth quadric and 
a (possibly singular) cubic. This proves (2).
\newline{\bf (3) }
{\bf Step 1 }
Since $g(M)=5$, $r(M)\geq10$ and $(r(M),\delta(M))\not=(10,0)$,
the type of $S$ is one of $[1;1,1]_{+}(1,b)$ 
$(2\leq b\leq 3)$ 
by \cite[p.410, Table 6 and p.444, Table 9]{Nakayama07}. 
By \cite[pp.494-495, Table 14]{Nakayama07}, 
$S$ is a subvariety of the weighted projective space
${\bf P}(1,1,2,2,4)$ defined by the equations:
$$
\begin{array}{lll}
\hbox{Lattice }M^{\perp}
&
\hbox{Type of }S
&
\hbox{Equations defining }S
\\
{\Bbb U}\oplus{\Bbb A}_{1}^{+}\oplus{\Bbb E}_{8}\oplus
{\Bbb A}_{1}^{\oplus 3-b}
&
[1;1,1]_{+}(b)\,\,(2\leq b\leq3)
&
\begin{cases}
xw=yz,
\\
z^{2}w=(xu-y^{2b-1}w^{3-b})\,x,
\\
zw^{2}=(xu-y^{2b-1}w^{3-b})\,y,
\end{cases}
\end{array}
$$
where
${\rm wt}(x)={\rm wt}(y)=1$, 
${\rm wt}(z)={\rm wt}(w)=2$, ${\rm wt}(u)=4$.
Notice that $S$ is {\it not} a complete intersection 
in ${\bf P}(1,1,2,2,4)$.
\par
Since ${\mathcal O}_{S}(-2K_{S})\cong{\mathcal O}_{{\bf P}(1,1,2,2,4)}(4)|_{S}$ 
by Lemma~\ref{Lemma3.2} below and hence 
$u-U(x,y,z,w)\in H^{0}(S,{\mathcal O}_{S}(-2K_{S}))$
for a weighted homogeneous polynomial 
$U(x,y,z,w)\in{\bf C}[x,y,z,w]$ of degree $4$,
$\Gamma$ is a subvariety of ${\bf P}(1,1,2,2)$ 
defined by the following equations:
$$
\begin{array}{ll}
\hbox{Type of }S
&
\hbox{Equations defining }\Gamma
\\

[1;1,1]_{+}(b)\,\,(2\leq b\leq3)
&
\begin{cases}
xw=yz,
\\
z^{2}w=(xU(x,y,z,w)-y^{2b-1}w^{3-b})\,x,
\\
zw^{2}=(xU(x,y,z,w)-y^{2b-1}w^{3-b})\,y.
\end{cases}
\end{array}
$$
{\bf Step 2 }
Set 
$\Sigma_{1}:=\{(x:y:z:w)\in{\bf P}(1,1,2,2);\,xw=yz\}$. 
By \cite[Lemma 7.6]{Nakayama07},
$\Sigma_{1}\cong{\bf P}({\mathcal O}_{{\bf P}^{1}}
\oplus{\mathcal O}_{{\bf P}^{1}}(1))$
is a Hirzebruch surface, which contains
$\Gamma$ as an irreducible divisor. 
The projection
$p\colon\Sigma_{1}\to{\bf P}^{1}$ is given by 
the formula:
$p\colon\Sigma_{1}\ni(x:y:z:w)\to(x:y)=(z:w)
\in{\bf P}^{1}$. By \cite[Lemma 7.6]{Nakayama07},
the negative section $\sigma$ of
$p\colon\Sigma_{1}\to{\bf P}^{1}$ is given by
$\sigma\colon{\bf P}^{1}\ni(z:w)\to(0:0:z:w)
\in\Sigma_{1}$. 
\par
Let $\ell$ and $C$ be the divisors on $\Sigma_{1}$
defined as
$$
\begin{array}{ll}
\ell
&:=
\{(0:y:0:w)\in\Sigma_{1};\,
(y:w)\in{\bf P}(1,2)\}
=
p^{-1}(0:1)
\subset
\Sigma_{1},
\\
C
&:=\{(x:y:z:w)\in\Sigma_{1};\,z^{2}w=(xU(x,y,z,w)-y^{2b-1}w^{3-b})x\}
\subset\Sigma_{1}.
\end{array}
$$
We have the equation of divisors $C=\Gamma+\ell$.
Since 
${\mathcal O}_{{\bf P}(1,1,2,2)}(2)|_{\Sigma_{1}}
\cong{\mathcal O}_{\Sigma_{1}}(\sigma+2\ell)$
by \cite[Lemma 7.6]{Nakayama07} and since
$z^{2}w-\{x\,U(x,y,z,w)-y^{2b-1}w^{3-b}\}x$
is an element of
$H^{0}({\bf P}(1,1,2,2),
{\mathcal O}_{{\bf P}(1,1,2,2)}(6))$, 
we get
$$
\Gamma+\ell
=
C
=
{\rm div}(z^{2}w-\{x\,U(x,y,z,w)-y^{2b-1}w^{3-b}\}x)
\in
|{\mathcal O}_{\Sigma_{1}}(3(\sigma+2\ell))|.
$$
Hence 
$\Gamma=C-\ell\in|{\mathcal O}_{\Sigma_{1}}(3\sigma+5\ell)|$.
Regard $H^{2}(\Sigma_{1},{\bf Z})$ as the Neron--Severi lattice
of $\Sigma_{1}$. Then we have the equations
$\sigma\cdot\sigma=-1$, $\ell\cdot\ell=0$, $\sigma\cdot\ell=1$.
Since $\Gamma$ is linearly equivalent to
$3\sigma+5\ell$, we get
$\Gamma\cdot\Gamma=21$, $\Gamma\cdot\sigma=2$.
\newline{\bf Step 3 }
Let $\pi\colon\Sigma_{1}\to{\bf P}^{2}$ 
be the blowing-down of the $(-1)$-curve $\sigma$
and set $\overline{\Gamma}:=\pi(\Gamma)$.
Let $\mu:={\rm mult}_{\pi(\sigma)}\overline{\Gamma}$
be the multiplicity of $\overline{\Gamma}$ at
$\pi(\sigma)$. Then $\mu=\Gamma\cdot\sigma=2$, 
so that $\pi(\sigma)$ is a double point of 
$\overline{\Gamma}$. Since $\Gamma$ is smooth and
since $\pi\colon\Gamma\setminus\sigma\to
\overline{\Gamma}\setminus\{\pi(\sigma)\}$ 
is an isomorphism, $\overline{\Gamma}$ has a unique 
singular point at $\pi(\sigma)$ and
$\pi|_{\Gamma}\colon\Gamma\to\overline{\Gamma}$
is the normalization.
Since $(\deg\overline{\Gamma})^{2}=
\overline{\Gamma}\cdot\overline{\Gamma}=
\Gamma\cdot\Gamma+\mu^{2}=25$,
we get $\deg\overline{\Gamma}=5$.
\par
Since $\pi\colon\Sigma_{1}\to{\bf P}^{2}$
is the blowing-down of $\sigma$ and since
$\overline{\Gamma}=\pi(\Gamma)$,
$\pi(\sigma)$ is a node of $\overline{\Gamma}$
if and only if $\Gamma$ intersects $\sigma$ transversally
at two different points. 
Since $\Gamma\cdot\sigma=2$, 
$\pi(\sigma)$ is a node of $\overline{\Gamma}$
if and only if $\#(\Gamma\cap\sigma)=2$.
By the definitions of $\Gamma$ and $\sigma$,
$$
\begin{aligned}
\#(\Gamma\cap\sigma)
&=
\#\{(x:y:z:w)\in\Gamma;\,x=y=0\}
\\
&=
\#\{(0:0:z:w)\in{\bf P}(1,1,2,2);\,zw=0\}=2.
\end{aligned}
$$
Hence $\overline{\Gamma}\subset{\bf P}^{2}$ is a quintic 
with one node, and $\Gamma$ is the normalization 
of $\overline{\Gamma}$.
\end{pf}

\subsection{Some log del Pezzo surfaces of index $2$}
\label{subsection:3.3}
\par
Let $n$ be an integer with $0\leq n\leq3$.
Let $(x:y:z:w:u)$ be the system of homogeneous coordinates 
of the weighted projective space 
${\bf P}(1,1,n+1,n+1,2(n+1))$ with weights 
${\rm wt}(x)={\rm wt}(y)=1$, ${\rm wt}(z)={\rm wt}(w)=n+1$,
${\rm wt}(u)=2(n+1)$.
Set
$$
W:=\{(x:y:z:w:u)\in{\bf P}(1,1,n+1,n+1,2(n+1));\,
xw=yz\}.
$$
In \cite[Prop.\,7.13]{Nakayama07}, Nakayama gave
a system of homogeneous polynomials that defines 
a log del Pezzo surface of index $2$ as a subvariety 
of $W$.

\begin{lemma}\label{Lemma3.2}
Let $S\subset W$ be a log del Pezzo surface of index $2$ as in \cite[Prop.\,7.13]{Nakayama07}.
Then the following isomorphism of holomorphic line bundles on $S$ holds:
$$
{\mathcal O}_{{\bf P}(1,1,n+1,n+1,2(n+1))}(2(n+1))|_{S}\cong{\mathcal O}_{S}(-2K_{S}).
$$
\end{lemma}

\begin{pf}
Let ${\mathcal F}$ be the vector bundle of rank $2$ 
over ${\bf P}(1,1,n+1,n+1)$ defined as
$$
{\mathcal F}
:=
{\mathcal O}_{{\bf P}(1,1,n+1,n+1)}
\oplus
{\mathcal O}_{{\bf P}(1,1,n+1,n+1)}(2(n+1)).
$$
Let 
${\bf P}({\mathcal F})\to{\bf P}(1,1,n+1,n+1)$
be the ${\bf P}^{1}$-bundle associated with ${\mathcal F}$
and let 
${\mathcal O}_{{\bf P}({\mathcal F})}(1)\to
{\bf P}({\mathcal F})$ be the tautological quotient line bundle.
Let 
$$
\Psi\colon{\bf P}({\mathcal F})\to
{\bf P}(1,1,n+1,n+1,2(n+1))
$$
be the birational morphism as in 
\cite[Lemma 7.5]{Nakayama07}.
\par
Set 
$\Sigma_{n}:=\{(x:y:z:w)\in{\bf P}(1,1,n+1,n+1);\,
xw=yz\}$.
By \cite[Lemma 7.6]{Nakayama07},
$\Sigma_{n}\cong{\bf P}({\mathcal O}_{{\bf P}^{1}}
\oplus{\mathcal O}_{{\bf P}^{1}}(n+1))$
is a Hirzebruch surface.
We set
$$
{\mathcal E}
:=
{\mathcal F}|_{\Sigma_{n}},
\qquad
{\bf P}({\mathcal E})
:=
{\bf P}({\mathcal F})|_{\Sigma_{n}},
\qquad
{\mathcal O}_{{\bf P}({\mathcal E})}(1)
:=
{\mathcal O}_{{\bf P}({\mathcal F})}(1)
|_{{\bf P}({\mathcal E})}.
$$
Then ${\bf P}({\mathcal E})\to\Sigma_{n}$ is 
the ${\bf P}^{1}$-bundle associated with ${\mathcal E}$,
and 
${\mathcal O}_{{\bf P}({\mathcal E})}(1)\to
{\bf P}({\mathcal E})$ is the tautological quotient line bundle.
Set $\Phi:=\Psi|_{{\bf P}({\mathcal E})}$.
Then $\Phi({\bf P}({\mathcal E}))=W$ 
by \cite[Prop.\,7.13]{Nakayama07}.
By \cite[Prop.\,7.8]{Nakayama07},
$\Phi\colon{\bf P}({\mathcal E})\to W$
is a birational morphism.
By \cite[p.461 l.10]{Nakayama07}, we have
$\Psi^{*}
{\mathcal O}_{{\bf P}(1,1,n+1,n+1,2(n+1))}(2(n+1))
\cong
{\mathcal O}_{{\bf P}({\mathcal F})}(1)$
and hence
\begin{equation}\label{eqn:(3.1)}
\Phi^{*}{\mathcal O}_{{\bf P}(1,1,n+1,n+1,2(n+1))}(2(n+1))|_{W}
\cong
{\mathcal O}_{{\bf P}({\mathcal F})}(1)|_{{\bf P}({\mathcal E})}
=
{\mathcal O}_{{\bf P}({\mathcal E})}(1).
\end{equation}
\par
Let $V\subset{\bf P}({\mathcal E})$ be the proper transform 
of $S$ with respect to the birational morphism
$\Phi\colon{\bf P}({\mathcal E})\to W$.
We set $\varphi:=\Phi|_{V}$
(cf. \cite[p.465, l.15]{Nakayama07}). Then
$\varphi\colon V\to S$ is a birational morphism.
By \cite[p.464 l.1-l.11]{Nakayama07}, we have
\begin{equation}\label{eqn:(3.2)}
\varphi^{*}{\mathcal O}_{S}(-2K_{S})\cong{\mathcal O}_{{\bf P}({\mathcal E})}(1)|_{V}.
\end{equation}
By \eqref{eqn:(3.1)} and \eqref{eqn:(3.2)}, we have an isomorphism of holomorphic line bundles on $V$:
\begin{equation}\label{eqn:(3.3)}
\varphi^{*}{\mathcal O}_{{\bf P}(1,1,n+1,n+1,2(n+1))}(2(n+1))|_{S}
\cong
\varphi^{*}{\mathcal O}_{S}(-2K_{S}).
\end{equation}
Since 
$\varphi|_{V\setminus\varphi^{-1}({\rm Sing}\,S)}\colon V\setminus\varphi^{-1}({\rm Sing}\,S)\to 
S\setminus{\rm Sing}\,S$ 
is an isomorphism by \cite[p.464, l.9-l.10]{Nakayama07} and since $S$ is normal, 
the desired isomorphism follows from \eqref{eqn:(3.3)}.
\end{pf}

\subsection{Even theta-characteristics on the main component of $X^{\iota}$}
\label{subsection:3.4}
\par

Recall that a {\it theta-characteristic} on a compact Riemann surface $C$ is a half canonical line bundle
on $C$, i.e., a holomorphic line bundle on $C$ whose square is the canonical line bundle of $C$.
A theta-characteristic $L$ is {\it even} if $h^{0}(L)\equiv0\mod2$. 
A theta-characteristic $L$ is effective if $h^{0}(L)>0$. If $g(C)$ denotes the genus of $C$, 
there are exactly $2^{g(C)-1}(2^{g(C)}+1)$ even theta-characteristics on $C$.

\begin{proposition}\label{Proposition3.3}
Let $C$ be a compact Riemann surface of genus $g(C)$.
\begin{itemize}
\item[(1)]
If $g(C)\leq2$, $C$ has no effective even 
theta-characteristics.
\item[(2)]
When $g(C)=3$, $C$ has no effective even 
theta-characteristics if and only if 
$C$ is non-hyperelliptic.
\item[(3)]
When $g(C)=4$, $C$ has no effective even 
theta-characteristics if and only if 
$C$ is a complete intersection of a smooth quadric
and a cubic in ${\bf P}^{3}$.
\item[(4)]
If $C$ is the normalization of an irreducible plane
quintic with one node, then $g(C)=5$ and
$C$ has no effective even theta-characteristics.
\end{itemize}
\end{proposition}

\begin{pf}
The assertions (1), (2), (3) are classical.
The assertion (4) follows from 
\cite[Lemma 0.18 (i), (ii)]{SmithVarley85}.
\end{pf}

\begin{proposition}\label{Proposition3.4}
Let $M\subset{\Bbb L}_{K3}$ be a primitive $2$-elementary
Lorentzian sublattice. 
If $r(M)>10$ or $(r(M),\delta(M))=(10,1)$, 
then there exists a $2$-elementary $K3$ surface 
$(X,\iota)$ of type $M$ such that 
$X^{\iota}$ has no effective even theta-characteristics.
\end{proposition}

\begin{pf}
The result follows from Propositions~\ref{Proposition3.1} and \ref{Proposition3.3}.
\end{pf}

\section{Automorphic forms on the period domain}
\label{sect:4}
\par

\subsection{Igusa's Siegel modular form and its pull-back on $\Omega_{M^{\perp}}$}
\label{subsect:4.1}
\par
Let ${\mathcal F}_{g}$ be the Hodge line bundle on 
${\mathcal A}_{g}$.
Then ${\mathcal F}_{g}$ is an ample line bundle on 
${\mathcal A}_{g}$ in the sense of orbifolds.
There is an integer $\nu\in{\bf N}$ such that 
${\mathcal F}_{g}^{\nu}$ is a line bundle on 
${\mathcal A}_{g}$ in the ordinary sense and 
such that ${\mathcal F}_{g}^{m\nu}$ extends to a very ample 
line bundle on ${\mathcal A}_{g}^{*}$ for $m\gg0$. In this case, 
let $\overline{\mathcal F}_{g}^{m\nu}$ denote the holomorphic 
extension of ${\mathcal F}_{g}^{m\nu}$ to 
${\mathcal A}_{g}^{*}$. An element of 
$H^{0}({\mathcal A}_{g},{\mathcal F}_{g}^{k})$
is identified with a Siegel modular form on ${\frak S}_{g}$ 
for $Sp_{2g}({\bf Z})$ of weight $k$. For $g>0$, we define
$$
\chi_{g}(\varSigma):=
\prod_{(a,b)\,{\rm even}}
\theta_{a,b}(\varSigma),
\qquad
\varSigma\in{\frak S}_{g},
$$
where $a,b\in\{0,\frac{1}{2}\}^{g}$ and
$\theta_{a,b}(\varSigma):=\sum_{n\in{\bf Z}^{g}}
\exp\{\pi i{}^{t}(n+a)\varSigma(n+a)+2\pi i{}^{t}(n+a)b\}$ 
is the corresponding theta 
constant. Here $(a,b)$ is {\it even} 
if $4{}^{t}ab\equiv0\mod2$.
When $g=0$, we define $\chi_{0}:=1$. 
By \cite[Lemma 10]{Igusa67}, $\chi_{g}^{8}$ 
is a Siegel modular form of weight $2^{g+1}(2^{g}+1)$. 
Let $\theta_{{\rm null},g}$ be the reduced divisor on
${\mathcal A}_{g}$ defined as
$$
\theta_{{\rm null},g}:=
\{[\varSigma]\in{\mathcal A}_{g};\,
\chi_{g}(\varSigma)=0\}.
$$
It is classical that $\theta_{{\rm null},2}={\mathcal N}_{2}$.
In Sect.\,\ref{sect:9}, $\chi_{g}^{8}$ shall play a crucial role.
\par
Define the {\it Petersson metric} on ${\mathcal F}_{g}$ by
\begin{equation}\label{eqn:(4.1)}
\|\xi\|^{2}(\varSigma):=(\det{\rm Im}\,\varSigma)|\xi|^{2},
\qquad
(\varSigma,\xi)\in{\frak S}_{g}\times{\bf C}.
\end{equation}
Since $\chi_{g}^{8}$ is a Siegel modular form,
$\|\chi_{g}^{8}\|^{2}=(\det{\rm Im}\,\varSigma)^{w(g)}|\chi_{g}(\varSigma)^{8}|^{2}$,
$w(g)=2^{g+1}(2^{g}+1)$, is a $C^{\infty}$ function on ${\mathcal A}_{g}$ in the sense of orbifolds.

\begin{lemma}\label{Lemma4.1}
Let $p\colon{\mathcal C}\to\varDelta$ be an ordinary
singular family of curves of genus $g>0$ such that
$C_{0}$ is irreducible. Let ${\frak o}:={\rm Sing}\,C_{0}$.
\begin{itemize}
\item[(1)]
There exists a holomorphic function
$h(t)\in{\mathcal O}(\varDelta)$ such that
$$
\log\|\chi_{g}(\varOmega(C_{t}))^{8}\|^{2}=
2^{2g-2}\log|t|^{2}+\log|h(t)|^{2}+O(\log\log|t|^{-1})
\qquad
(t\to0).
$$
\item[(2)]
If $g=1$ or $g=2$, then $h(0)\not=0$.
\end{itemize}
\end{lemma}

\begin{pf}
We follow \cite[p.370, Sect.\,3]{Mumford83}.
For $\varSigma\in{\frak S}_{g}$, 
we write $\varSigma=\binom{z\,{}^{t}\omega}{\omega\,Z}$,
where $z\in{\frak H}$, $\omega\in{\bf C}^{g-1}$,
$Z\in{\frak S}_{g-1}$.
\newline{\bf (1) }
Since $C_{0}$ is an irreducible curve of arithmetic genus
$g>0$ with one node, the normalization of $C_{0}$ is
a smooth curve of genus $g-1$.
By \cite[Cor.3.8]{Fay73}, there exists a holomorphic function 
$\psi(t)$ on $\varDelta$ with values in complex symmetric 
$g\times g$-matrices such that
\begin{equation}\label{eqn:(4.2)}
\varOmega(C_{t})=\left[\frac{\log t}{2\pi i}A+\psi(t)\right]\in{\mathcal A}_{g},
\qquad
A=
\begin{pmatrix}
1&{}^{t}{\bf 0}_{g-1}\\
{\bf 0}_{g-1}&O_{g-1}
\end{pmatrix}.
\end{equation}
Write 
$\psi(0)=\binom{\psi_{0}\,{}^{t}\omega_{0}}{\omega_{0}\,\,Z_{0}}$.
Then $Z_{0}\in{\frak S}_{g-1}$ and
$\lim_{t\to0}\varOmega(C_{t})=[Z_{0}]\in{\mathcal A}_{g-1}
\subset{\mathcal A}_{g}^{*}$.
\par 
For $a,b\in\{0,\frac{1}{2}\}^{g}$, write
$a=(a_{1},a')$, $b=(b_{1},b')$, where 
$a_{1},b_{1}\in\{0,\frac{1}{2}\}$,
$a',b'\in\{0,\frac{1}{2}\}^{g-1}$. Let $a_{1}=\frac{1}{2}$. 
There is a holomorphic function $f_{a',b'}(\zeta,\omega,Z)$ 
such that
\begin{equation}\label{eqn:(4.3)}
\begin{aligned}
\theta_{a,b}(\varSigma)
&=
\sum_{n=(n_{1},n')\in{\bf Z}\times{\bf Z}^{g-1}}
e^{
\pi i(n_{1}+\frac{1}{2})^{2}z
+2\pi i(n_{1}+\frac{1}{2}){}^{t}\omega(n'+a')
+\pi i{}^{t}(n'+a')Z(n'+a')+2\pi i{}^{t}(n+a)b}
\\
&=
e^{\frac{\pi iz}{4}}\{e^{-\pi ib_{1}}\theta_{a',b'}(-\omega/2,Z)
+
e^{\pi ib_{1}}\theta_{a',b'}(\omega/2,Z)
+
e^{2\pi iz}f_{a',b'}(e^{2\pi iz},\omega,Z)\}
\\
&=
e^{\frac{\pi iz}{4}}\{2i^{2b_{1}}\,\theta_{a',b'}(\omega/2,Z)
+
e^{2\pi iz}f_{a',b'}(e^{2\pi iz},\omega,Z)\},
\end{aligned}
\end{equation}
where we used $4{}^{t}ab\in2{\bf Z}$ and the identity
$\theta_{a',b'}(-\omega/2,Z)=
(-1)^{4{}^{t}a'b'}\theta_{a',b'}(\omega/2,Z)$
\cite[p.167, Prop.\,3.14]{Mumford83b} to get the third equality.
The number of even $(a,b)$ with $a_{1}=1/2$ is given by
$2^{2(g-1)}$. 
\par
Similarly, let $a_{1}=0$. Then the pair $(a',b')$ must be even.
There is a holomorphic function $g_{a',b'}(\zeta,\omega,Z)$ 
such that
\begin{equation}\label{eqn:(4.4)}
\begin{aligned}
\theta_{a,b}(\varSigma)
&=
\sum_{n=(n_{1},n')\in{\bf Z}\times{\bf Z}^{g-1}}
e^{
\pi in_{1}^{2}z
+2\pi in_{1}{}^{t}\omega(n'+a')
+\pi i{}^{t}(n'+a')Z(n'+a')+2\pi i{}^{t}(n+a)b}
\\
&=
(-1)^{2{}^{t}a'b'}\theta_{a',b'}(Z)
+
e^{\pi iz}g_{a',b'}(e^{\pi iz},\omega,Z).
\end{aligned}
\end{equation}
By \eqref{eqn:(4.3)}, \eqref{eqn:(4.4)}, there is a holomorphic function 
$F(\zeta,\omega,Z)$ such that
\begin{equation}\label{eqn:(4.5)}
\chi_{g}(\varSigma)^{8}
=
\prod_{(a,b)\,{\rm even}}\theta_{a,b}(\varSigma)^{8}
=
(e^{\frac{\pi iz}{4}})^{8\cdot2^{2(g-1)}}F(e^{\pi iz},\omega,Z)
=
(e^{2\pi iz})^{2^{2g-2}}F(e^{\pi iz},\omega,Z).
\end{equation}
Since $\chi_{g}^{8}$ is a Siegel modular form and hence 
$\chi_{g}(\varOmega+A)^{8}=\chi_{g}(\varOmega)^{8}$,
$F(\zeta,\omega,Z)$ is an even function in $\zeta$.
By \eqref{eqn:(4.2)}, $z=(\log t)/2\pi i+\psi_{11}(t)$ for some $\psi_{11}(t)\in{\mathcal O}(\varDelta)$. 
Hence $\exp(2\pi iz)=t\exp(2\pi i\psi_{11}(t))$.
By \eqref{eqn:(4.5)}, there exists $h(t)\in{\mathcal O}(\varDelta)$ such that
\begin{equation}\label{eqn:(4.6)}
\chi_{g}\left(\frac{\log t}{2\pi i}A+\psi(t)\right)^{8}=t^{2^{2g-2}}h(t).
\end{equation}
Since 
${\rm Im}(\frac{\log t}{2\pi i}A+\psi(t))=
(-\frac{1}{2\pi}\log|t|)A+{\rm Im}\,\psi(0)+O(|t|)$
with
$\psi(0)=\binom{\psi_{0}\,{}^{t}\omega_{0}}{\omega_{0}\,\,Z_{0}}$,
$Z_{0}\in{\frak S}_{g-1}$, we get
\begin{equation}\label{eqn:(4.7)}
\det{\rm Im}\left(\frac{\log t}{2\pi i}A+\psi(t)\right)=-\frac{\det{\rm Im}\,Z_{0}}{2\pi}\log|t|+O(1).
\end{equation}
By \eqref{eqn:(4.2)}, \eqref{eqn:(4.6)}, \eqref{eqn:(4.7)}, we get (1).
\newline{\bf (2) }
Let $g=1$. Since $p\colon{\mathcal C}\to\varDelta$ is 
an ordinary singular family of elliptic curves, 
$(\varDelta,0)$ is regarded as a local coordinate of 
${\mathcal A}_{1}^{*}$ centered at the cusp $+i\infty$.
Since $\chi_{1}^{8}(\tau)=\eta(\tau)^{24}$ vanishes of order $1$
at the cusp of ${\mathcal A}_{1}^{*}$, we get (2) in this case.
\par
Let $g=2$. Then $\omega_{0}\in{\bf C}$, $Z_{0}\in{\frak H}$ and $\theta_{a',b'}(Z_{0})\not=0$ in \eqref{eqn:(4.4)}.
Set $\Lambda_{0}:={\bf Z}+Z_{0}{\bf Z}$.
By \eqref{eqn:(4.3)}, the assertion (2) follows if $\theta_{a,b}(\omega_{0}/2,Z_{0})\not=0$ for all $(a,b)\in\{0,1/2\}$. 
Since 
${\rm div}\,\theta_{a,b}(\cdot,Z_{0})=
[(a+\frac{1}{2})Z_{0}+(b+\frac{1}{2})]\in{\bf C}/\Lambda_{0}$
by \cite[Lemma 4.1]{Mumford83b}, it suffices to prove
$\frac{\omega_{0}}{2}\not\in
(\frac{1}{2}{\bf Z})Z_{0}+\frac{1}{2}{\bf Z}$, i.e.,
$\omega_{0}\not\in\Lambda_{0}$. 
Let $\iota\colon\widehat{C}_{0}\to C_{0}$ be the normalization.
Since ${\frak o}$ is the node of $C_{0}$, we can write
$\iota^{-1}({\frak o})=\{\widehat{\frak o}_{1},
\widehat{\frak o}_{2}\}$ with 
$\widehat{\frak o}_{1}\not=\widehat{\frak o}_{2}$. 
By \cite[p.53 Cor.\,3.8]{Fay73}, there exist a symplectic basis 
$\{\alpha,\beta\}$ of $H_{1}(\widehat{C}_{0},{\bf Z})$ and 
a holomorphic $1$-form $v$ on $\widehat{C}_{0}$ such that
$\int_{\alpha}v=1$, $\int_{\beta}v=Z_{0}$ and
$\int_{\widehat{\frak o}_{1}}^{\widehat{\frak o}_{2}}v=\omega_{0}$. 
Since $\widehat{\frak o}_{1}\not=\widehat{\frak o}_{2}$,
we get $\omega_{0}\not\in\Lambda_{0}$. This proves (2).
\end{pf}

\par
Let $\omega_{{\frak S}_{g}}$ be the 
$Sp_{2g}({\bf Z})$-invariant K\"ahler form 
on ${\frak S}_{g}$ defined as
$$
\omega_{{\frak S}_{g}}(\varSigma)
:=-dd^{c}\log\det{\rm Im}\,\varSigma,
\qquad
\varSigma\in{\frak S}_{g}.
$$ 
Let $\omega_{{\mathcal A}_{g}}$ be the K\"ahler form 
on ${\mathcal A}_{g}$ in the sense of orbifolds induced 
from $\omega_{{\frak S}_{g}}$. Then 
$$
\omega_{{\mathcal A}_{g}}=c_{1}({\mathcal F}_{g},\|\cdot\|).
$$
\par
Let ${\mathcal I}(M)\subset{\bf Z}$ be the ideal defined
as follows: $q\in{\mathcal I}(M)$ if and only if
there exists 
$\overline{\mathcal F}_{g(M)}^{q}\in
H^{1}({\mathcal A}_{g(M)}^{*},
{\mathcal O}^{*}_{{\mathcal A}_{g(M)}^{*}})$ with
$\overline{\mathcal F}_{g(M)}^{q}|_{{\mathcal A}_{g(M)}}=
{\mathcal F}_{g(M)}^{q}$.
\par
Let 
$i\colon\Omega_{M^{\perp}}^{o}
\cup{\mathcal D}_{M^{\perp}}^{o}
\hookrightarrow \Omega_{M^{\perp}}$ 
be the inclusion. For $q\in{\mathcal I}(M)$, we set
$$
\lambda_{M}^{q}:=
i_{*}{\mathcal O}_{\Omega_{M^{\perp}}^{o}
\cup{\mathcal D}_{M^{\perp}}^{o}}
(J_{M}^{*}\overline{\mathcal F}_{g(M)}^{q}).
$$
By \cite[Lemma 3.6]{Yoshikawa04} and Proposition~\ref{Proposition2.2},
the ${\mathcal O}_{\Omega_{M^{\perp}}}$-module $\lambda_{M}^{q}$ 
is an invertible sheaf on $\Omega_{M^{\perp}}$.
We identify $\lambda_{M}^{q}$ with the corresponding 
holomorphic line bundle on $\Omega_{M^{\perp}}$. 
By \cite[Lemma 3.7]{Yoshikawa04} and Proposition~\ref{Proposition2.2}, 
the $O(M^{\perp})$-action on 
$\lambda_{M}^{q}|_{\Omega_{M^{\perp}}^{o}
\cup{\mathcal D}_{M^{\perp}}^{o}}$
induced from the $O(M^{\perp})$-equivariant map $J_{M}$,
extends to the one on $\lambda_{M}^{q}$. Hence
$\lambda_{M}^{q}$ is equipped with the structure of 
an $O(M^{\perp})$-equivariant line bundle on 
$\lambda_{M}^{q}$.
\par
Let $\|\cdot\|_{\lambda_{M}^{q}}$ be the 
$O(M^{\perp})$-invariant Hermitian metric on 
$\lambda_{M}^{q}|_{\Omega_{M^{\perp}}^{o}}$ defined as
$$
\|\cdot\|_{\lambda_{M}^{q}}:=(J_{M}^{o})^{*}\|\cdot\|.
$$
By \eqref{eqn:(4.1)}, $(J_{M}^{o})^{*}\omega_{{\mathcal A}_{g(M)}}$ 
is a $C^{\infty}$ closed semi-positive $(1,1)$-form on 
$\Omega_{M^{\perp}}^{o}$ such that
$q\,(J_{M}^{o})^{*}\omega_{{\mathcal A}_{g(M)}}=
c_{1}(\lambda_{M}^{q}|_{\Omega_{M^{\perp}}^{o}},
\|\cdot\|_{\lambda_{M}^{q}})$.
Since $\dim\Omega_{M^{\perp}}\setminus
(\Omega_{M^{\perp}}^{o}\cup{\mathcal D}_{M^{\perp}}^{o})
\leq\dim\Omega_{M^{\perp}}-2$
when $r(M)\leq18$, 
we can define the closed positive 
$(1,1)$-current $J_{M}^{*}\omega_{{\mathcal A}_{g(M)}}$
on $\Omega_{M^{\perp}}$ as the trivial extension of 
$(J_{M}^{o})^{*}\omega_{{\mathcal A}_{g(M)}}$
from $\Omega_{M^{\perp}}^{o}$ to $\Omega_{M^{\perp}}$
by \cite[p.\,53 Th.\,1]{Siu74} and \cite[Th.\,3.9]{Yoshikawa04}. 
When $r(M)=19$,
$(J_{M}^{o})^{*}\omega_{{\mathcal A}_{g(M)}}$
extends trivially to a closed positive 
$(1,1)$-current $J_{M}^{*}\omega_{{\mathcal A}_{g(M)}}$
on $\Omega_{M^{\perp}}$, because 
$(J_{M}^{o})^{*}\omega_{{\mathcal A}_{g(M)}}$
has Poincar\'e growth along ${\mathcal D}_{M^{\perp}}$
by \cite[Prop.\,3.8]{Yoshikawa04}.
By \cite[p.\,53 Th.\,1]{Siu74} and \cite[Th.\,3.13]{Yoshikawa04},
the Hermitian metric $\|\cdot\|_{\lambda_{M}^{q}}$ 
on $\lambda_{M}^{q}|_{\Omega_{M^{\perp}}^{o}}$ extends to 
a singular Hermitian metric on $\lambda_{M}^{q}$ 
with curvature current
\begin{equation}\label{eqn:(4.8)}
c_{1}(\lambda_{M}^{q},\|\cdot\|_{\lambda_{M}^{q}})=q\,J_{M}^{*}\omega_{{\mathcal A}_{g(M)}}.
\end{equation}
\par
Let $\ell\in{\bf Z}_{>0}$ be such that
$2^{g(M)+1}(2^{g(M)}+1)\ell\in{\mathcal I}(M)$.
Then
${\mathcal F}_{g(M)}^{2^{g(M)+1}(2^{g(M)}+1)\ell}$
extends to a holomorphic line bundle
on ${\mathcal A}_{g(M)}^{*}$.
Since $\chi_{g(M)}^{8\ell}$ is a holomorphic section of 
${\mathcal F}_{g(M)}^{2^{g(M)+1}(2^{g(M)}+1)\ell}$,
$J_{M}^{*}\chi_{g(M)}^{8\ell}$ is an 
$O(M^{\perp})$-invariant holomorphic section of 
$\lambda_{M}^{2^{g(M)+1}(2^{g(M)}+1)\ell}$.
If
$J_{M}^{o}(\Omega_{M^{\perp}}^{o})\not\subset
\theta_{{\rm null},g(M)}$,  we define 
$$
{\frak D}:={\rm div}(J_{M}^{*}\chi_{g(M)}^{8\ell}).
$$ 
Since $J_{M}$ is $O(M^{\perp})$-equivariant 
with respect to the trivial $O(M^{\perp})$-action on 
${\mathcal A}_{g(M)}^{*}$, ${\frak D}$ 
is an $O(M^{\perp})$-invariant effective divisor on
$\Omega_{M^{\perp}}$. 
By \cite[p.\,53 Th.\,1]{Siu74}, \cite[Th.\,3.13]{Yoshikawa04} and \eqref{eqn:(4.8)},
$\log\|J_{M}^{*}\chi_{g(M)}^{8}\|$ lies in $L^{1}_{\rm loc}(\Omega_{M^{\perp}})$ and 
satisfies the following equation of currents on $\Omega_{M^{\perp}}$
\begin{equation}\label{eqn:(4.9)}
-dd^{c}\log\|J_{M}^{*}\chi_{g(M)}^{8\ell}\|^{2}
=
2^{g(M)+1}(2^{g(M)}+1)\ell\,J_{M}^{*}\omega_{{\mathcal A}_{g(M)}}-\delta_{\frak D}.
\end{equation}
\par
Recall that the divisor ${\mathcal D}'_{M^{\perp}}$ was defined in Sect.\,\ref{subsect:1.4}.

\begin{proposition}\label{Proposition4.2}
Assume that $r(M)>10$ or $(r(M),\delta(M))=(10,1)$ and that $g(M)>0$. Hence $M$ is non-exceptional.
Let $\ell\in{\bf Z}_{>0}$ be such that $2^{g(M)+1}(2^{g(M)}+1)\ell\in{\mathcal I}(M)$.
Then the following hold:
\begin{itemize}
\item[(1)]
$J_{M}^{o}(\Omega_{M^{\perp}}^{o})\not\subset
\theta_{{\rm null},g(M)}$.
\item[(2)]
There exist an integer $a\in{\bf Z}_{\geq0}$ 
and an $O(M^{\perp})$-invariant (possibly empty)
effective divisor $E$ on $\Omega_{M^{\perp}}$ such that
$\dim(E\cap{\mathcal D}'_{M^{\perp}})<
\dim{\mathcal D}'_{M^{\perp}}$ and
$$
{\frak D}=2(2^{2g(M)-2}+a)\ell\,{\mathcal D}'_{M^{\perp}}+E.
$$
In particular, the following equation of currents
on $\Omega_{M^{\perp}}$ holds:
$$
-dd^{c}\log\|J_{M}^{*}\chi_{g(M)}^{8\ell}\|^{2}
=
2^{g(M)+1}(2^{g(M)}+1)\ell\,
J_{M}^{*}\omega_{{\mathcal A}_{g(M)}}-
2(2^{2g(M)-2}+a)\ell\,
\delta_{{\mathcal D}'_{M^{\perp}}}-\delta_{E}.
$$
\item[(3)]
If $g(M)=1$ or $g(M)=2$, then $a=0$ and $E=0$ in $(2)$.
\end{itemize}
\end{proposition}

\begin{pf}
{\bf (1) }
Let $(X,\iota)$ be a $2$-elementary $K3$ surface
of type $M$ and let $C$ be the main component of $X^{\iota}$.
By Riemann's theorem \cite[p.338]{GriffithsHarris78}
and Riemann's singularities theorem \cite{Mayer64},
\cite[p.348]{GriffithsHarris78}, 
$C$ has an effective even theta-characteristic
if and only if 
$J_{M}^{o}(X,\iota)\in\theta_{{\rm null},g(M)}$.
Since $r(M)>10$ or $(r(M),\delta(M))=(10,1)$,
there exists by Proposition~\ref{Proposition3.4}
a $2$-elementary $K3$ surface $(X,\iota)$ of type $M$ 
such that the main component of $X^{\iota}$
has {\it no} effective even theta-characteristics, i.e.,
$J_{M}^{o}(X,\iota)\not\in\theta_{{\rm null},g(M)}$.
This proves (1).
\par{\bf (2) }
Since ${\frak D}$ is an $O(M^{\perp})$-invariant 
effective divisor on $\Omega_{M^{\perp}}$, we can write 
${\frak D}=\sum_{d\in\Delta'_{M^{\perp}}}m(d)\,H_{d}+E$,
where $m(d)\in{\bf Z}_{\geq0}$ and $E$ is an effective 
divisor on $\Omega_{M^{\perp}}$ with
$\dim({\mathcal D}'_{M^{\perp}}\cap E)\leq
\dim{\mathcal D}'_{M^{\perp}}-1$.
Since $g(H_{d})=H_{g(d)}$ for all $g\in O(M^{\perp})$
and $d\in\Delta'_{M^{\perp}}$, 
the $O(M^{\perp})$-invariance of ${\frak D}$ implies 
that $m(g(d))=m(d)$ for all $g\in O(M^{\perp})$ and
$d\in\Delta'_{M^{\perp}}$.
Since $O(M^{\perp})$ acts transitively on 
$\Delta'_{M^{\perp}}$ by \cite[Prop.\,3.3]{FinashinKharlamov08}
and Proposition~\ref{Proposition11.6} (5) below, 
there exists $\alpha\in{\bf Z}_{\geq0}$ with
\begin{equation}\label{eqn:(4.10)}
{\frak D}=\alpha\,{\mathcal D}'_{M^{\perp}}+E.
\end{equation}
\par
Let $d\in\Delta'_{M^{\perp}}$
and ${\frak p}\in\overline{H}_{d}^{o}$. 
Let $\gamma\colon\varDelta\to{\mathcal M}_{M^{\perp}}$ 
be a holomorphic curve intersecting 
$\overline{H}_{d}^{o}$ transversally at 
$\gamma(0)={\frak p}$ such that
$\gamma(\varDelta\setminus\{0\})\subset
{\mathcal M}_{M}\setminus
(\overline{\mathcal D}_{M^{\perp}}\cup{\frak D})$. 
By Theorem~\ref{Theorem2.3} (1), there exists an 
ordinary singular family of $2$-elementary $K3$ surfaces
$p_{\mathcal Z}\colon({\mathcal Z},\iota)\to\varDelta$ 
of type $M$ with Griffiths period map $\gamma$,
such that $(\widetilde{Z}_{0},\widetilde{\iota}_{0})$ 
is a $2$-elementary $K3$ surface of type 
$[M\perp d]$ with Griffiths period $\gamma(0)$. 
\par
Since the natural projection 
$\varPi_{M^{\perp}}\colon\Omega_{M^{\perp}}
\to{\mathcal M}_{M^{\perp}}$
is doubly ramified along $H_{d}^{o}$ by
\cite[Prop.\,1.9 (4)]{Yoshikawa04}, there exists
a holomorphic curve $c\colon\varDelta\to\Omega_{M^{\perp}}$ 
intersecting $H_{d}^{o}$ transversally at 
$c(0)\in H_{d}^{o}$ such that
$\varPi_{M^{\perp}}(c(t))=\gamma(t^{2})$. 
Hence we have
\begin{equation}\label{eqn:(4.11)}
J_{M}(c(t))=\varOmega(Z_{t^{2}}^{\iota_{t^{2}}}).
\end{equation}
\par
Since $d\in\Delta'_{M^{\perp}}$, $\iota$ is of type $(2,1)$ by Theorem~\ref{Theorem2.3} (2).
By \cite[Prop.\,2.5]{Yoshikawa04},
$p|_{{\mathcal Z}^{\iota}}\colon{\mathcal Z}^{\iota}\to\varDelta$ is an ordinary singular family of curves.
Let ${\mathcal C}\subset{\mathcal Z}^{\iota}$ be the connected component such that
$C_{t}:={\mathcal C}\cap Z_{t}^{\iota_{t}}$ is the main component of $Z_{t}^{\iota_{t}}$ for all
$t\in\varDelta\setminus\{0\}$. 
Since the normalization of $Z_{0}^{\iota_{0}}$ is given by $(\widetilde{Z}_{0})^{\widetilde{\iota}_{0}}$, 
the normalization of $C_{0}$ has genus $g(M)-1$ by Theorem~\ref{Theorem2.3} (2). 
Hence $C_{0}$ is singular and $p|_{\mathcal C}\colon{\mathcal C}\to\varDelta$
is an ordinary singular family of curves. 
Since the normalization of $C_{0}$ has genus $g(M)-1$ and since $C_{0}$ has a unique node
as its singular set, $C_{0}$ is irreducible.
\par
We apply Lemma~\ref{Lemma4.1} to the ordinary singular family
$p|_{\mathcal C}\colon{\mathcal C}\to\varDelta$ with irreducible $C_{0}$. 
Since $\varOmega(C_{t})=\varOmega(Z_{t}^{\iota_{t}})$ for all $t\in\varDelta\setminus\{0\}$, 
there exists $h(t)\in{\mathcal O}(\varDelta)$ by Lemma~\ref{Lemma4.1} (1) such that
\begin{equation}\label{eqn:(4.12)}
\log\|\chi_{g(M)}(\varOmega(Z_{t}^{\iota_{t}}))^{8}\|^{2}
=
2^{2g(M)-2}\,\log|t|^{2}+\log|h(t)|^{2}+O(\log\log|t|^{-1}).
\end{equation}
Since $\gamma(\varDelta\setminus\{0\})\cap{\frak D}=\emptyset$ by the choice of $\gamma$,
$h(t)$ does not vanish identically on $\varDelta$ by \eqref{eqn:(4.6)}.
Let $a\in{\bf Z}_{\geq0}$ be the multiplicity of $h(t)$ at $t=0$. 
By \eqref{eqn:(4.11)}, \eqref{eqn:(4.12)}, we get
\begin{equation}\label{eqn:(4.13)}
\log\|\chi_{g(M)}(J_{M}(c(t)))^{8\ell}\|^{2}=2(2^{2g(M)-2}+a)\ell\,\log|t|^{2}+O(\log\log|t|^{-1}),
\end{equation}
which yields that $H_{d}\subset{\rm supp}\,{\frak D}$ for $d\in\Delta'_{M^{\perp}}$.
Comparing \eqref{eqn:(4.9)}, \eqref{eqn:(4.10)} and \eqref{eqn:(4.13)}, 
we get $\alpha=2(2^{2g(M)-2}+a)\ell$ in \eqref{eqn:(4.10)}.
Since ${\frak D}$ and ${\mathcal D}'_{M^{\perp}}$ are $O(M^{\perp})$-invariant, so is $E$ by \eqref{eqn:(4.10)}. 
This proves (2).
\newline{\bf (3) }
Let $g(M)=1$ or $g(M)=2$. By Proposition~\ref{Proposition3.3} (1), we get 
the inclusion ${\frak D}\subset{\mathcal D}_{M^{\perp}}$.
This, together with \eqref{eqn:(4.10)}, implies the inclusion $E\subset{\mathcal D}''_{M^{\perp}}$. 
Since $r(M)\geq10$ and hence $M\not\cong{\Bbb A}_{1}^{+}\oplus{\Bbb A}_{1}^{\oplus8}$,
there exists by Propositions~\ref{Proposition2.8} and \ref{Proposition2.10} a dense Zariski open 
subset $U$ of ${\mathcal D}''_{M^{\perp}}$ with 
$J_{M}(U)\subset{\mathcal A}_{g(M)}\setminus\theta_{{\rm null},g(M)}$.
By the inclusion $E\subset{\mathcal D}''_{M^{\perp}}$, we get
$J_{M}(E\cap U)\subset{\mathcal A}_{g(M)}\setminus\theta_{{\rm null},g(M)}$. 
If $E\not=0$, $J_{M}^{*}\chi_{g(M)}^{8\ell}$ would not vanish on the non-empty dense Zariski open subset 
$E\cap U$ of $E$, which contradicts the fact that $E\subset{\frak D}={\rm div}(J_{M}^{*}\chi_{g(M)}^{8\ell})$.
This proves that $E=0$. The equality $a=0$ follows from \eqref{eqn:(4.12)}, \eqref{eqn:(4.13)} and
the nonvanishing $h(0)\not=0$ in Lemma~\ref{Lemma4.1} (2). This proves the proposition.
\end{pf}

\begin{lemma}\label{Lemma4.3}
Let $p\colon{\mathcal C}\to\varDelta$ be an ordinary singular family of curves of genus $2$ 
such that $C_{0}$ is the join of two elliptic curves intersecting at one point transversally. Then
$$
\log\|\chi_{2}(\varOmega(C_{t}))^{8}\|^{2}=8\log|t|^{2}+O(\log\log|t|^{-1})
\qquad
(t\to0).
$$
\end{lemma}

\begin{pf}
Since $g=2$ and $C_{0}$ is reducible, we deduce from
\cite[Cor.3.8]{Fay73} the existence of a holomorphic map
$\psi\colon\varDelta\to{\frak S}_{2}$ with
$$
{\varOmega}(C_{t})=[\psi(t)],
\quad
\psi(0)=
\begin{pmatrix}
\psi_{1}&0\\0&\psi_{2}
\end{pmatrix},
\quad
\psi'(0)=
\begin{pmatrix}
0&a\\a&0
\end{pmatrix},
\quad
\psi_{1},\psi_{2}\in{\frak H}, 
\quad a\not=0.
$$
The result follows from e.g. \cite[Eq.\,(A.24)]{Yoshikawa99}.
\end{pf}

\begin{proposition}\label{Proposition4.4}
Let $g(M)=2$ and $r(M)<10$, i.e., $M\cong{\Bbb A}_{1}^{+}\oplus{\Bbb A}_{1}^{\oplus8}$. 
Let ${\mathcal H}_{M^{\perp}}({\bf 1}_{M^{\perp}},-1/2)$ be the Heegner divisor defined as
$$
{\mathcal H}_{M^{\perp}}({\bf 1}_{M^{\perp}},-1/2)
:=
\sum_{\{\lambda\in{\bf 1}_{M^{\perp}}+M^{\perp};\,
\lambda^{2}=-1/2\}/\pm1}
H_{\lambda}
=
\sum_{d\in\Delta''_{M^{\perp}}/\pm1,\,
d/2\in{\bf 1}_{M^{\perp}}+M^{\perp}}
H_{d}.
$$
Then the following equation of divisors on $\Omega_{M^{\perp}}$ holds:
$$
{\rm div}(J_{M}^{*}\chi_{2}^{8\ell})
=
8\ell\,{\mathcal D}'_{M^{\perp}}
+
16\ell\,{\mathcal H}_{M^{\perp}}({\bf 1}_{M^{\perp}},-1/2).
$$
In particular, the following equations of currents
on $\Omega_{M^{\perp}}$ holds:
$$
-dd^{c}\log\|J_{M}^{*}\chi_{2}^{8\ell}\|^{2}
=
40\ell\,J_{M}^{*}\omega_{{\mathcal A}_{2}}
-
8\ell\,\delta_{{\mathcal D}'_{M^{\perp}}}
-
16\ell\,\delta_{{\mathcal H}_{M^{\perp}}
({\bf 1}_{M^{\perp}},-\frac{1}{2})}.
$$
\end{proposition}

\begin{pf}
Let $d\in\Delta''_{M^{\perp}}$ and $d/2={\bf 1}_{M^{\perp}}$.
Since $M\cong{\Bbb A}_{1}^{+}\oplus{\Bbb A}_{1}^{\oplus8}$,
we get $[M\perp d]\cong{\Bbb U}\oplus{\Bbb E}_{8}(2)$ by the proof of Proposition~\ref{Proposition2.10}.
Let ${\frak p}\in\overline{H}_{d}^{o}$. 
Let $\gamma\colon\varDelta\to{\mathcal M}_{M^{\perp}}$ be a holomorphic curve intersecting 
$\overline{H}_{d}^{o}$ transeversally at $\gamma(0)={\frak p}$ such that
$\gamma(\varDelta\setminus\{0\})\subset{\mathcal M}_{M}\setminus
(\overline{\mathcal D}_{M^{\perp}}\cup{\frak D})$. 
By Theorem~\ref{Theorem2.3} (1), there exists an ordinary singular family of $2$-elementary $K3$ surfaces
$p_{\mathcal Z}\colon({\mathcal Z},\iota)\to\varDelta$ of type $M$ with Griffiths period map $\gamma$,
such that $(\widetilde{Z}_{0},\widetilde{\iota}_{0})$ is a $2$-elementary $K3$ surface of type 
$[M\perp d]\cong{\Bbb U}\oplus{\Bbb E}_{8}(2)$ with Griffiths period $\gamma(0)$. 
As in the proof of Proposition~\ref{Proposition4.2} (2), there exists a holomorphic curve 
$c\colon\varDelta\to\Omega_{M^{\perp}}$ intersecting $H_{d}^{o}$ transversally at $c(0)\in H_{d}^{o}$ 
and satisfying \eqref{eqn:(4.11)}.
If $\iota$ is of type $(0,3)$, then $Z_{0}^{\iota_{0}}$ is the disjoint union of a smooth curve of genus $2$
and an isolated point by \cite[Prop.\,2.5]{Yoshikawa04}, which implies that
$J_{M}(c(0))\in{\mathcal A}_{2}\setminus{\mathcal N}_{2}$.
By Theorem~\ref{Theorem2.5}, this leads to the contradiction
$$
J_{M}(c(0))=J^{o}_{[M\perp d]}(c(0))=J^{o}_{{\Bbb U}\oplus{\Bbb E}_{8}(2)}(c(0))\in{\mathcal N}_{2},
$$
where the last inclusion follows from Proposition~\ref{Proposition2.1} (2). Hence $\iota$ is of type $(2,1)$.
\par
By \cite[Prop.\,2.5]{Yoshikawa04}, $p|_{{\mathcal Z}^{\iota}}\colon{\mathcal Z}^{\iota}\to\varDelta$ 
is an ordinary singular family of curves.
Since the normalization of $(Z_{0})^{\iota_{0}}$ is isomorphic to $(\widetilde{Z}_{0})^{\widetilde{\iota}_{0}}$
by Theorem~\ref{Theorem2.3} (2) and since $(\widetilde{Z}_{0},\widetilde{\iota}_{0})$ is of type
$[M\perp d]\cong{\Bbb U}\oplus{\Bbb E}_{8}(2)$, we deduce from Proposition~\ref{Proposition2.1} (2) 
that $(Z_{0})^{\iota_{0}}$ is the join of two elliptic curves intersecting at one point transversally. 
By Lemma~\ref{Lemma4.3}, we get
\begin{equation}\label{eqn:(4.14)}
\log\|\chi_{2}(\varOmega(Z_{t}^{\iota_{t}}))^{8}\|^{2}=8\,\log|t|^{2}+O(\log\log|t|^{-1})
\qquad
(t\to0).
\end{equation}
By \eqref{eqn:(4.11)} and \eqref{eqn:(4.14)}, we get
\begin{equation}\label{eqn:(4.15)}
\log\|\chi_{2}(J_{M}(c(t)))^{8}\|^{2}=16\,\log|t|^{2}+O(\log\log|t|^{-1})
\qquad
(t\to0).
\end{equation}
\par
By Proposition~\ref{Proposition2.1} (3), we get
$J_{M}(\Omega_{M^{\perp}}^{o})=
J_{M}^{o}(\Omega_{M^{\perp}}^{o})\subset{\mathcal A}_{2}\setminus\theta_{{\rm null},2}$.
By Proposition~\ref{Proposition2.10}, we get 
$J_{M}(\bigcup_{d\in\Delta''_{M^{\perp}},\,d/2\not\equiv{\bf 1}_{M^{\perp}}}H_{d}^{o})
\subset{\mathcal A}_{2}\setminus\theta_{{\rm null},2}$.
By these two inclusions,
$$
J_{M}(\Omega_{M^{\perp}}^{o}\cup\bigcup_{d\in\Delta''_{M^{\perp}},\,d/2\not\equiv{\bf 1}_{M^{\perp}}}H_{d}^{o})
\subset
{\mathcal A}_{2}\setminus\theta_{{\rm null},2},
$$
which implies that $J_{M}^{*}\chi_{2}^{8\ell}$ does not
vanish on
$\Omega_{M^{\perp}}^{o}\cup
\bigcup_{d\in\Delta''_{M^{\perp}},\,
d/2\not\equiv{\bf 1}_{M^{\perp}}}H_{d}^{o}$. 
Hence
$$
\begin{aligned}
(\Omega_{M^{\perp}}^{o}\cup{\mathcal D}_{M^{\perp}}^{o})
\cap
{\frak D}
&\subset
(\Omega_{M^{\perp}}^{o}\cup{\mathcal D}_{M^{\perp}}^{o})
\setminus
(
\Omega_{M^{\perp}}^{o}\cup
\bigcup_{d\in\Delta''_{M^{\perp}},\,
d/2\not\equiv{\bf 1}_{M^{\perp}}}H_{d}^{o}
)
\\
&=
{\mathcal D}_{M^{\perp}}^{o}\setminus
\bigcup_{d\in\Delta''_{M^{\perp}},\,
d/2\not\equiv{\bf 1}_{M^{\perp}}}H_{d}^{o}
\\
&\subset
{\mathcal D}'_{M^{\perp}}\cup 
{\mathcal H}_{M^{\perp}}({\bf 1}_{M^{\perp}},-1/2).
\end{aligned}
$$
Since $\Omega_{M^{\perp}}\setminus
(\Omega_{M^{\perp}}^{o}\cup{\mathcal D}_{M^{\perp}}^{o})$ 
is an analytic subset of codimension $2$ 
in $\Omega_{M^{\perp}}$,
we get 
\begin{equation}\label{eqn:(4.16)}
{\frak D}
\subset
{\mathcal D}'_{M^{\perp}}\cup{\mathcal H}_{M^{\perp}}({\bf 1}_{M^{\perp}},-1/2).
\end{equation}
Since the proof of Proposition~\ref{Proposition4.2} (2) works in the case
$M\cong{\Bbb A}_{1}^{+}\oplus{\Bbb A}_{1}^{\oplus8}$,
\eqref{eqn:(4.13)} remains valid. Moreover, we get $a=0$ in \eqref{eqn:(4.13)} by Lemma~\ref{Lemma4.1} (2). 
The desired formula follows from \eqref{eqn:(4.10)}, \eqref{eqn:(4.13)} with $a=0$, 
\eqref{eqn:(4.15)}, \eqref{eqn:(4.16)}.
\end{pf}

\subsection{Automorphic forms on $\Omega_{\Lambda}^{+}$}
\label{subsect:4.2}
\par
Let $\Lambda$ be a lattice of signature $(2,r(\Lambda)-2)$.
We fix a vector $l_{\Lambda}\in\Lambda\otimes{\bf R}$ with $\langle l_{\Lambda},l_{\Lambda}\rangle\geq0$,
and we set
$$
j_{\Lambda}(\gamma,[\eta])
:=
\frac{\langle\gamma(\eta),l_{\Lambda}\rangle}{\langle\eta,l_{\Lambda}\rangle}
\qquad
[\eta]\in\Omega_{\Lambda}^{+},
\quad
\gamma\in O^{+}(\Lambda).
$$
Since $H_{l_{\Lambda}}=\emptyset$, $j_{\Lambda}(\gamma,\cdot)$
is a nowhere vanishing holomorphic function on $\Omega_{\Lambda}^{+}$.
\par
Let $\Gamma\subset O^{+}(\Lambda)$ be a cofinite subgroup. 
A holomorphic function $f\in{\mathcal O}(\Omega_{\Lambda}^{+})$ is called an
{\it automorphic form on $\Omega_{\Lambda}^{+}$ for $\Gamma$ of weight $p$} if
$$
f(\gamma\cdot[\eta])=\chi(\gamma)\,j_{\Lambda}(\gamma,[\eta])^{p}\,f([\eta]),
\qquad
[\eta]\in\Omega_{\Lambda}^{+},
\quad
\gamma\in\Gamma,
$$
where $\chi\colon\Gamma\to{\bf C}^{*}$ is a unitary character.
For an automorphic form $f$ on $\Omega_{\Lambda}^{+}$ for $\Gamma$ of weight $p$, 
the Petersson norm $\|f\|$ is the function on $\Omega_{\Lambda}^{+}$ defined as
$$
\|f([\eta])\|^{2}:=K_{\Lambda}([\eta])^{p}\,|f([\eta])|^{2},
\qquad
K_{\Lambda}([\eta]):=\frac{\langle\eta,\bar{\eta}\rangle}{|\langle\eta,l_{\Lambda}\rangle|^{2}}.
$$ 
If $r(\Lambda)\geq 5$, then $\|f\|^{2}$ is a $\Gamma$-invariant $C^{\infty}$ function on $\Omega_{\Lambda}^{+}$,
because the group $\Gamma/[\Gamma,\Gamma]$ is finite and Abelian and hence $\chi$ is finite in this case.
\par
We also consider automorphic forms on $\Omega_{M^{\perp}}^{+}$ with values in the sheaf $\lambda_{M}^{q}$. 
Let $M\subset{\Bbb L}_{K3}$ be a primitive $2$-elementary Lorentzian sublattice.
Let $\chi$ be a character of $O^{+}({M}^{\perp})$. Let $p,q\in{\bf Z}$.
Then $\Psi\in H^{0}(\Omega_{M^{\perp}}^{+},\lambda_{M}^{q})$ 
is called an automorphic form on $\Omega_{M^{\perp}}^{+}$ for $O^{+}({M}^{\perp})$ of weight $(p,q)$ 
if for all $\gamma\in O^{+}(M^{\perp})$,
$$
\Psi(\gamma\cdot[\eta])=\chi(\gamma)\,j_{M^{\perp}}(\gamma,[\eta])^{p}\,\gamma(\Psi([\eta])),
\qquad
[\eta]\in\Omega_{M^{\perp}}^{+}.
$$
\par
For an automorphic form $\Psi$ on $\Omega_{M^{\perp}}^{+}$ for $O^{+}(M^{\perp})$ of weight $(p,q)$, 
the Petersson norm of $\Psi$ is a $C^{\infty}$ function on $\Omega_{M^{\perp}}^{+}$ defined as
\begin{equation}\label{eqn:(4.17)}
\|\Psi([\eta])\|^{2}:=K_{M^{\perp}}([\eta])^{p}\cdot\|\Psi([\eta])\|_{\lambda_{M}^{q}}^{2},
\qquad
[\eta]\in\Omega_{M^{\perp}}^{+}.
\end{equation}
%

\section{The invariant $\tau_{M}$ of $2$-elementary $K3$ surfaces of type $M$}
\label{sect:5}
\par
Let $(X,\iota)$ be a $2$-elementary $K3$ surface 
of type $M$. Identify ${\bf Z}_{2}$ with 
the subgroup of ${\rm Aut}(X)$ generated by $\iota$. 
Let $\kappa$ be a ${\bf Z}_{2}$-invariant 
K\"ahler form on $X$. Set 
${\rm vol}(X,\kappa):=(2\pi)^{-2}\int_{X}\kappa^{2}/2!$.
Let $\eta$ be a nowhere vanishing holomorphic
$2$-form on $X$. The $L^{2}$-norm of $\eta$ is
defined as $\|\eta\|_{L^{2}}^{2}:=
(2\pi)^{-2}\int_{X}\eta\wedge\bar{\eta}$.
\par
Let 
$\square_{0,q}=2(\bar{\partial}+\bar{\partial}^{*})^{2}$ 
be the $\bar{\partial}$-Laplacian acting on 
$C^{\infty}$ $(0,q)$-forms on $X$. 
Let $\sigma(\square_{0,q})$ be the spectrum of $\square_{0,q}$. 
For $\lambda\in\sigma(\square_{0,q})$, let $E_{0,q}(\lambda)$ be 
the eigenspace of $\square_{0,q}$ with respect to the eigenvalue 
$\lambda$. Since ${\bf Z}_{2}$ preserves $\kappa$, 
$E_{0,q}(\lambda)$ is a finite-dimensional representation of 
${\bf Z}_{2}$. For $s\in{\bf C}$, set
$$
\zeta_{0,q}(\iota)(s)
:=
\sum_{\lambda\in\sigma(\square_{0,q})\setminus\{0\}}
{\rm Tr}\,(\iota|_{E_{0,q}(\lambda)})\,\lambda^{-s}.
$$ 
Then $\zeta_{0,q}(\iota)(s)$ converges absolutely when 
${\rm Re}\,s>\dim X$, admits a meromorphic continuation to 
the complex plane ${\bf C}$, and is holomorphic at $s=0$. 
The {\it equivariant analytic torsion} of
the trivial Hermitian line bundle on 
$(X,\kappa)$ is defined as
$$
\tau_{{\bf Z}_{2}}(X,\kappa)(\iota)
:=
\exp[
-\sum_{q\geq0}(-1)^{q}q\,\zeta'_{0,q}(\iota)(0)].
$$
We refer to \cite{RaySinger73}, \cite{BGS88}, 
\cite{BismutLebeau91}, \cite{GilletSoule92},
\cite{Bismut95}, \cite{Ma00}, \cite{KohlerRoessler01}
for more about equivariant and non-equivariant
analytic torsion.
\par
Let $X^{\iota}=\sum_{i}C_{i}$ be the decomposition
of the fixed point set of $\iota$ into the connected
components. Let $c_{1}(C_{i},\kappa|_{C_{i}})$ 
be the Chern form of $(TC_{i},\kappa|_{C_{i}})$ and let
$\tau(C_{i},\kappa|_{C_{i}})$ be the analytic torsion
of the trivial Hermitian line bundle on 
$(C_{i},\kappa|_{C_{i}})$. We define 
$$
\begin{aligned}
\tau_{M}(X,\iota)
&:=
{\rm vol}(X,(2\pi)^{-1}\kappa)^{\frac{14-r(M)}{4}}
\tau_{{\bf Z}_{2}}(X,\kappa)(\iota)
\prod_{i}{\rm Vol}(C_{i},(2\pi)^{-1}\kappa|_{C_{i}})
\tau(C_{i},\kappa|_{C_{i}})
\\
&\quad
\times
\exp\left[
\frac{1}{8}\int_{X^{\iota}}
\log\left.\left(\frac{\eta\wedge\bar{\eta}}
{\kappa^{2}/2!}\cdot
\frac{{\rm Vol}(X,(2\pi)^{-1}\kappa)}{\|\eta\|_{L^{2}}^{2}}
\right)\right|_{X^{\iota}}
c_{1}(X^{\iota},\kappa|_{X^{\iota}})
\right],
\end{aligned}
$$
which is independent of the choice of $\kappa$
by \cite[Th.\,5.7]{Yoshikawa04}.
Hence $\tau_{M}(X,\iota)$ is an invariant of 
the pair $(X,\iota)$, 
so that $\tau_{M}$ descends to a function on 
${\mathcal M}_{M^{\perp}}^{o}$.

\begin{theorem}\label{Theorem5.1}
There exist an integer $\nu\in{\bf Z}_{>0}$ and an automorphic form $\Phi_{M}$ on $\Omega_{M^{\perp}}$ 
for $O^{+}(M^{\perp})$ of weight $(\nu(r(M)-6),4\nu)$ with zero divisor $\nu\,{\mathcal D}_{M^{\perp}}$ 
such that for every $2$-elementary $K3$ surface $(X,\iota)$ of type $M$, 
$$
\tau_{M}(X,\iota)=
\|\Phi_{M}(\overline{\varpi}_{M}(X,\iota))\|^{-\frac{1}{2\nu}}.
$$
\end{theorem}

\begin{pf}
See \cite[Main Th.]{Yoshikawa04}, \cite[Th.\,1.1]{Yoshikawa09b} and Proposition~\ref{Proposition11.2} below.
\end{pf}

\section{Borcherds products}
\label{sect:6}
\par

\subsection{Modular forms for $Mp_{2}({\bf Z})$}
\label{subsect:6.1}
\par
Recall that ${\frak H}\subset{\bf C}$ is the complex 
upper half-plane. 
Let $Mp_{2}({\bf Z})$ be the metaplectic double cover
of $SL_{2}({\bf Z})$ (cf. \cite[Sect.\,2]{Borcherds00}),
which is generated by the two elements
$S:=(\binom{0\,-1}{1\,\,\,\,0},\sqrt{\tau})$ and
$T:=(\binom{1\,1}{0\,1},1)$.
For 
$\gamma=(\binom{a\,b}{c\,d},\sqrt{c\tau+d})
\in Mp_{2}({\bf Z})$ and $\tau\in{\frak H}$, we define
$j(\gamma,\tau):=\sqrt{c\tau+d}$ and
$\gamma\cdot\tau:=(a\tau+b)/(c\tau+d)$.
\par
Let $M$ be an even lattice. 
Let ${\bf C}[A_{M}]$ be the group ring of 
the discriminant group $A_{M}$. 
Let $\{{\bf e}_{\gamma}\}_{\gamma\in A_{M}}$
be the standard basis of ${\bf C}[A_{M}]$.
The Weil representation 
$\rho_{M}\colon Mp_{2}({\bf Z})\to
GL({\bf C}[A_{M}])$ is defined as follows
\cite[Sect.\,2]{Borcherds00}:
\begin{equation}\label{eqn:(6.1)}
\rho_{M}(T)\,{\bf e}_{\gamma}:=e^{\pi i\gamma^{2}}{\bf e}_{\gamma},
\qquad
\rho_{M}(S)\,{\bf e}_{\gamma}
:=
\frac{i^{-\sigma(M)/2}}{|A_{M}|^{1/2}}\sum_{\delta\in A_{M}}e^{-2\pi i\langle\gamma,\delta\rangle}{\bf e}_{\delta}.
\end{equation}
\par
A ${\bf C}[A_{M}]$-valued holomorphic function 
$F(\tau)$ on ${\frak H}$ is a
{\it modular form of type $\rho_{M}$
with weight $w\in\frac{1}{2}{\bf Z}$}
if the following conditions (a), (b) are satisfied: 
\begin{itemize}
\item[(a)]
For $\gamma\in Mp_{2}({\bf Z})$ and $\tau\in{\frak H}$,
$F(\gamma\cdot\tau)=
j(\gamma,\tau)^{2w}\,\rho_{M}(\gamma)\cdot F(\tau)$.
\item[(b)] 
$F(\tau)=
\sum_{\gamma\in A_{M}}{\bf e}_{\gamma}
\sum_{k\in\frac{1}{l}{\bf Z}}c_{\gamma}(k)\,
e^{2\pi i k\tau}$,
where $l$ is the level of $M$,
$c_{\gamma}(k)\in{\bf Z}$ for all 
$k\in\frac{1}{l}{\bf Z}$ and $c_{\gamma}(k)=0$ for $k\ll0$.
\end{itemize}
By the first condition of \eqref{eqn:(6.1)}, 
\cite[Eq.\,(1.4)]{Bruinier02} and Condition (a), we get
\begin{equation}\label{eqn:(6.2)}
c_{\gamma}(k)
=
\begin{cases}
\begin{array}{lll}
0&{\rm if}&
k\not\in\gamma^{2}/2+{\bf Z},
\\
c_{-\gamma}(k)
&{\rm if}&
k\in\gamma^{2}/2+{\bf Z}.
\end{array}
\end{cases}
\end{equation}
\par
The group $O(M)$ acts on ${\bf C}[A_{M}]$ by
$g({\bf e}_{\gamma}):={\bf e}_{\bar{g}(\gamma)}$,
where $\bar{g}\in O(q_{M})$ is the element induced
by $g\in O(M)$.
For a modular form $F$ of type $\rho_{M}$, 
we define
${\rm Aut}(M,F):=\{g\in O(M);\,g(F)=F\}$.
Then ${\rm Aut}(M,F)$ is a cofinite subgroup of $O(M)$,
since ${\rm Aut}(M,F)\supset\ker\{O(M)\to O(q_{M})\}$.

\subsection{Borcherds products}
\label{subsect:6.2}
\par
Let $\Lambda$ be an even lattice of signature $(2,r(\Lambda)-2)$. 
Assume that $\Lambda={\Bbb U}(N)\oplus L$, for simplicity.
A vector of $\Lambda\otimes{\bf Q}$ is denoted by $(m,n,v)$, where 
$m,n\in{\bf Q}$ and $v\in L\otimes{\bf Q}$.
We write a vector of $A_{\Lambda}$ in the same manner.
If $F(\tau)=\sum_{\gamma\in A_{\Lambda}}
f_{\gamma}(\tau)\,{\bf e}_{\gamma}$ 
is a modular form of type $\rho_{\Lambda}$,
then $F(\tau)$ induces 
a modular form $F|_{L}(\tau)$ of type $\rho_{L}$ 
with the same weight as follows 
\cite[Th.\,5.3]{Borcherds98}:
\begin{equation}\label{eqn:(6.3)}
F|_{L}(\tau):=\sum_{\lambda\in A_{L}}f_{L+\lambda}(\tau)\,{\bf e}_{\lambda},
\qquad
f_{L+\lambda}(\tau):=\sum_{n=0}^{N-1}f_{(\frac{n}{N},0,\bar{\lambda})}(\tau).
\end{equation}
\par
Write
$F|_{L}(\tau)=
\sum_{\gamma\in A_{L}}{\bf e}_{\gamma}\sum_{k\in\frac{\gamma^{2}}{2}+{\bf Z}}c_{L,\gamma}(k)\,e^{2\pi i k\tau}$.
By \cite[Sect.\,6, p.517]{Borcherds98}, 
$F|_{L}(\tau)$ induces a chamber structure of 
${\mathcal C}_{L}^{+}$:
\begin{equation}\label{eqn:(6.4)}
({\mathcal C}_{L}^{+})^{0}_{F|_{L}}
:=
{\mathcal C}_{L}^{+}\setminus\bigcup_{\lambda\in L^{\lor},\,
\lambda^{2}<0,\,c_{L,\bar{\lambda}}(\lambda^{2}/2)\not=0}
h_{\lambda}
=
\amalg_{\alpha\in A}{\mathcal W}_{\alpha},
\end{equation}
where $h_{\lambda}=\lambda^{\perp}
=\{v\in L\otimes{\bf R};\,
\langle v,\lambda\rangle=0\}$ and
$\{{\mathcal W}_{\alpha}\}_{\alpha\in A}$ 
is the set of connected components of 
$({\mathcal C}_{L}^{+})^{0}_{F|_{L}}$.
Each component ${\mathcal W}_{\alpha}$ is called 
a {\it Weyl chamber} of $F|_{L}(\tau)$.
If $\lambda\in L\otimes{\bf R}$
satisfies $\langle\lambda,w\rangle>0$ 
for all $w\in{\mathcal W}_{\alpha}$, we write
$\lambda\cdot{\mathcal W}_{\alpha}>0$.

\begin{theorem}\label{Theorem6.1}
Let 
$F(\tau)=\sum_{\gamma\in A_{\Lambda}}{\bf e}_{\gamma}
\sum_{k\in\frac{\gamma^{2}}{2}+{\bf Z}}c_{\gamma}(k)\,e^{2\pi i k\tau}$ 
be a modular form of type $\rho_{\Lambda}$ with weight $\sigma(\Lambda)/2$.
Then there exists a meromorphic automorphic form $\Psi_{\Lambda}(z,F)$ on $\Omega_{\Lambda}^{+}$
for ${\rm Aut}(\Lambda,F)\cap O^{+}(\Lambda)$ of weight $c_{0}(0)/2$ such that
$$
{\rm div}(\Psi_{\Lambda}(\cdot,F))
=
\frac{1}{2}\sum_{\lambda\in\Lambda^{\lor},\,\lambda^2<0}
c_{\bar{\lambda}}(\lambda^{2}/2)\,H_{\lambda}
=
\sum_{\lambda\in\Lambda^{\lor}/\pm1,\,\lambda^2<0}
c_{\bar{\lambda}}(\lambda^{2}/2)\,H_{\lambda}.
$$
If ${\mathcal W}$ is a Weyl chamber of $F|_{L}$, then there exists a vector 
$\varrho(L,F|_{L},{\mathcal W})\in L\otimes{\bf Q}$ such that $\Psi_{\Lambda}(z,F)$ is expressed as 
the following infinite product near the cusp under the identification \eqref{eqn:(1.2)}: 
For $z\in L\otimes{\bf R}+i\,{\mathcal W}$ with $({\rm Im}\,z)^{2}\gg0$, 
$$
\Psi_{\Lambda}(z,F)
=
e^{2\pi i\langle\varrho(L,F|_{L},{\mathcal W}),z\rangle}
\prod_{\lambda\in L^{\lor},\,\lambda\cdot{\mathcal W}>0}
\prod_{n\in{\bf Z}/N{\bf Z}}
(1-e^{2\pi i(\langle\lambda,z\rangle+\frac{n}{N})}
)^{c_{(\frac{n}{N},0,\bar{\lambda})}
(\lambda^{2}/2)}.
$$
\end{theorem}

\begin{pf}
See \cite[Th.\,13.3]{Borcherds98}, \cite[Th.\,3.22]{Bruinier02}.
\end{pf}

The automorphic form $\Psi_{\Lambda}(z,F)$ is called
the {\it Borcherds product} or the {\it Borcherds lift} of $F(\tau)$, 
and the vector $\varrho(L,F|_{L},{\mathcal W})$ is 
called the {\it Weyl vector} of $\Psi_{\Lambda}(z,F)$.
See \cite[Th.\,10.4]{Borcherds98}, \cite[p.321 Correction]{Borcherds00} 
for an explicit formula for $\varrho(L,F|_{L},{\mathcal W})$.

\section{$2$-elementary lattices and elliptic modular forms}
\label{sect:7}
\par
Throughout Section 7, we assume that {\it $\Lambda$ is an even, $2$-elementary lattice.}
\par
Set
$M\Gamma_{0}(4):=
\{(\binom{a\,b}{c\,d},\sqrt{c\tau+d})
\in Mp_2({\bf Z});\,c\equiv0\mod 4\}$.
Let $w\in\frac{1}{2}{\bf Z}$ and
let $\chi\colon M\Gamma_{0}(4)\to{\bf C}^{*}$ 
be a character.
A holomorphic function $f(\tau)$ on ${\frak H}$ 
is a modular form for $M\Gamma_{0}(4)$ of weight 
$w$ with character $\chi$
if the following (a), (b) are satisfied: 
\begin{itemize}
\item[(a)]
$f(\gamma\cdot\tau)=
j(\gamma,\tau)^{2w}\chi(\gamma)\,f(\tau)$
for all $\gamma\in M\Gamma_{0}(4)$ and
$\tau\in{\frak H}$.
\item[(b)]
$f(\tau)=
\sum_{k\in\frac{1}{4}{\bf Z}}c(k)\,
e^{2\pi i k\tau}$ with $c(k)=0$ for $k\ll0$.
\end{itemize}
\par
Set $q=e^{2\pi i\tau}$ for $\tau\in{\frak H}$.
Let $\eta(\tau)=q^{1/24}\prod_{n=1}^{\infty}(1-q^{n})$
be the Dedekind $\eta$-function and let
$\vartheta_{2}(\tau)$, $\vartheta_{3}(\tau)$,
$\vartheta_{4}(\tau)$ be the Jacobi theta functions:
$$
\vartheta_{2}(\tau)=
\sum_{n\in\bf Z}q^{(n+\frac{1}{2})^{2}/2},
\qquad
\vartheta_{3}(\tau)=\sum_{n\in\bf Z}q^{n^{2}/2},
\qquad
\vartheta_{4}(\tau)=
\sum_{n\in\bf Z}(-1)^{n}q^{n^{2}/2}.
$$
Notice that we use the
notation $q=e^{2\pi i\tau}$ while $q=e^{\pi i\tau}$
in \cite[Chap.\,4]{ConwaySloan99}.
Recall that ${\Bbb A}_{1}$ is the negative-definite
one-dimensional $A_{1}$-lattice $\langle-2\rangle$
and ${\Bbb A}_{1}^{+}=\langle2\rangle$.
For $d\in\{0,1\}$,
let $\theta_{{\Bbb A}_{1}^{+}+d/2}(\tau)$ 
be the theta function of ${\Bbb A}_{1}^{+}$:
$$
\theta_{{\Bbb A}_{1}^{+}}(\tau):=
\vartheta_{3}(2\tau),
\qquad
\theta_{{\Bbb A}_{1}^{+}+1/2}(\tau):=
\vartheta_{2}(2\tau).
$$
By \cite[Lemma 5.2]{Borcherds00}, there exists
a character 
$\chi_{\theta}\colon M\Gamma_{0}(4)\to\{\pm1,\pm i\}$
such that $\theta_{{\Bbb A}_{1}^{+}}(\tau)$ is 
a modular form for $M\Gamma_{0}(4)$ of weight $1/2$
with character $\chi_{\theta}$.
\par
For $k\in{\bf Z}$, define 
$f_{k}^{(0)}(\tau),f_{k}^{(1)}(\tau)
\in{\mathcal O}({\frak H})$ and 
$\{c_{k}^{(0)}(l)\}_{l\in{\bf Z}}$,
$\{c_{k}^{(1)}(l)\}_{l\in{\bf Z}+k/4}$ by
$$
\begin{aligned}
f_{k}^{(0)}(\tau)
&:=
\frac{\eta(2\tau)^{8}\,
\theta_{{\Bbb A}_{1}^{+}}(\tau)^{k}}
{\eta(\tau)^8\eta(4\tau)^8}
=
\sum_{l\in{\bf Z}}c_{k}^{(0)}(l)\,q^{l}
=
q^{-1}+8+2k+O(q),
\\
f_{k}^{(1)}(\tau)
&:=
-16\frac{\eta(4\tau)^8\,
\theta_{{\Bbb A}_{1}^{+}+\frac{1}{2}}(\tau)^{k}}
{\eta(2\tau)^{16}}
=
\sum_{l\in\frac{k}{4}+{\bf Z}}2c_{k}^{(1)}(l)\,q^{l}
=-2^{k+4}\,q^{\frac{k}{4}}\{1+(k+16)q^{2}+O(q^{4})\}.
\end{aligned}
$$
We define holomorphic functions 
$g_{k}^{(i)}(\tau)\in{\mathcal O}({\frak H})$, 
$i\in{\bf Z}/4{\bf Z}$ by
$$
g_{k}^{(i)}(\tau):=
\sum_{l\equiv i\,\,{\rm mod}\,\,4}
c_{k}^{(0)}(l)\,q^{l/4}.
$$
By definition,
$$
\sum_{i\in{\bf Z}/4{\bf Z}}g^{(i)}_{k}(\tau)
=
\frac{\eta(\tau/2)^{8}\,
\theta_{{\Bbb A}_{1}^{+}}(\tau/4)^{k}}
{\eta(\tau)^{8}\eta(\tau/4)^{8}}
=
f^{(0)}_{k}(\tau/4).
$$
\par 
For a modular form $\phi(\tau)$ of weight $l$ 
for $M\Gamma_{0}(4)$ and
for $g\in Mp_2(\bf Z)$, we define
$$
\phi|_{g}(\tau):=\phi(g\cdot\tau)\,j(g,\tau)^{-2l}.
$$ 
The following key construction of modular forms 
of type $\rho_{\Lambda}$ is due to Borcherds.

\begin{proposition}\label{Proposition7.1}
Let $\phi(\tau)$ be a modular form for $M\Gamma_{0}(4)$ of weight $l$ with character 
$\chi_{\theta}^{\sigma(\Lambda)}$ and set
$$
{\mathcal B}_{\Lambda}[\phi](\tau)
:=
\sum_{g\in M\Gamma_{0}(4)\backslash Mp_{2}({\bf Z})}\phi|_{g}(\tau)\,\rho_{\Lambda}(g^{-1})\,{\bf e}_{0}.
$$ 
\begin{itemize}
\item[(1)]
$\rho_{\Lambda}(g){\bf e}_{0}=\chi_{\theta}(g)^{\sigma(\Lambda)}{\bf e}_{0}$
for all $g\in M\Gamma_{0}(4)$.
\item[(2)]
${\mathcal B}_{\Lambda}[\phi](\tau)$ is independent of the choice of representatives of 
$M\Gamma_{0}(4)\backslash Mp_{2}({\bf Z})$.
Moreover, ${\mathcal B}_{\Lambda}[\phi](\tau)$ is a modular form for $Mp_{2}({\bf Z})$ of type 
$\rho_{\Lambda}$ of weight $l$.
\end{itemize}
\end{proposition}

\begin{pf}
{\bf (1) }
Let $\chi_{A_{\Lambda}}$ be the character of $M\Gamma_{0}(4)$
defined in \cite[Th.\,5.4]{Borcherds00}. 
Let $k\in{\bf Z}_{>0}$ be such that $\sigma(\Lambda)+8k\geq0$.
Since $|A_{\Lambda}|\cdot2^{\sigma(\Lambda)+8k}
=2^{2\{2+4k+(l(\Lambda)-r(\Lambda))/2\}}$, we get
$\binom{-1}{|A|}=1$ and
$\chi_{|A_{\Lambda}|\cdot2^{\sigma(\Lambda)+8k}}\equiv1$
by the definitions of the character $\chi_{n}$ and 
the symbol $\binom{c}{d}$ in \cite[p.328]{Borcherds00}.
Hence we get $\chi_{A_{\Lambda}}=\chi_{\theta}^{\sigma(\Lambda)}$
by \cite[Th.\,5.4]{Borcherds00}.
\par
Set 
$M\Gamma(4)':=\{(\binom{a\,b}{c\,d},\sqrt{c\tau+d})\in Mp_2({\bf Z});\,b\equiv c\equiv0\mod 4\}$.
By \cite[Th.\,5.4]{Borcherds00}, we get
$\rho_{\Lambda}(g){\bf e}_{0}=\chi_{\theta}(g)^{\sigma(\Lambda)}{\bf e}_{0}$
for all $g\in M\Gamma(4)'$, because $\chi_{A_{\Lambda}}=\chi_{\theta}^{\sigma(\Lambda)}$.
Since the coset $M\Gamma_{0}(4)/M\Gamma(4)'$ is represented by $\{1,T,T^{2},T^{3}\}$, 
any $g\in M\Gamma_{0}(4)$can be expressed as $g=T^{a}g_{0}$, where $a\in{\bf Z}$ and 
$g_{0}\in M\Gamma(4)'$. 
Since $\rho_{\Lambda}(T){\bf e}_{0}={\bf e}_{0}$ by \eqref{eqn:(6.1)} 
and $\chi_{\theta}(T)=1$ by \cite[Lemma 5.2]{Borcherds00}, we get
$$
\rho_{\Lambda}(g){\bf e}_{0}
=
\rho(T)^{a}\rho(g_{0}){\bf e}_{0}
=
\chi_{\theta}(g_{0})^{\sigma(\Lambda)}{\bf e}_{0}
=
\chi_{\theta}(T^{a}g_{0})^{\sigma(\Lambda)}{\bf e}_{0}
=
\chi_{\theta}(g)^{\sigma(\Lambda)}{\bf e}_{0}.
$$
Since $g\in M\Gamma_{0}(4)$ is an arbitrary element, we get (1).
\newline{\bf (2) }
By (1), the result follows from \cite[Th.\,6.2]{Scheithauer06}. 
See also \cite[Lemma 2.6]{Borcherds98} and
\cite[Proof of Lemma 11.1]{Borcherds00}.
\end{pf}

\begin{lemma}\label{Lemma7.2}
The function $f_{k}^{(0)}(\tau)$ is a modular form for $M\Gamma_{0}(4)$ of weight $-4+\frac{k}{2}$
with character $\chi_{\theta}^{k}$.
\end{lemma}

\begin{pf}
The result follows from
\cite[Lemma\,5.2 and Th.\,6.2]{Borcherds00}.
\end{pf}

Set $Z:=S^2=(-\binom{1\,0}{0\,1},i)$ and
$V:=S^{-1}T^{2}S=(\binom{\,\,1\,\,\,\,0}{-2\,\,1},\sqrt{-2\tau+1})$.

\begin{lemma}\label{Lemma7.3}
The coset $M\Gamma_{0}(4)\backslash Mp_{2}({\bf Z})$ is represented by $\{1,S,ST,ST^2,ST^3,V\}$.  
\end{lemma}

\begin{pf}
Since none of two elements of $\{1,S,ST,ST^2,ST^3,V\}$
represent the same element of
$M\Gamma_{0}(4)\backslash Mp_{2}({\bf Z})$ and 
$\#M\Gamma_{0}(4)\backslash Mp_{2}({\bf Z})=6$,
we get the result.
\end{pf}

Recall that ${\bf 1}_{\Lambda}\in A_{\Lambda}$ was defined in Sect.\,\ref{subsect:1.2}. Define 
${\bf v}_{0},{\bf v}_{1},{\bf v}_{2},{\bf v}_{3}\in{\bf C}[A_{\Lambda}]$ by
$$
{\bf v}_{k}:=\sum_{\delta\in A_{\Lambda},\,\delta^{2}\equiv k/2\,\,{\rm mod}\,\,2}{\bf e}_{\delta}.
$$

\begin{lemma}\label{Lemma7.4}
The following identities hold:
$$
(1)\quad
\rho_{\Lambda}((ST^{l})^{-1})\,{\bf e}_{0}
=
i^{\frac{\sigma(\Lambda)}{2}}2^{-\frac{l(\Lambda)}{2}}\sum_{k=0}^{3}i^{-lk}\,{\bf v}_{k},
\qquad
(2)\quad
\rho_{\Lambda}(V^{-1})\,{\bf e}_{0}={\bf e}_{{\bf 1}_{\Lambda}}.
$$
\end{lemma}

\begin{pf}
{\bf (1)}
Since $S^{-1}=SZ^{3}$ and since 
$\rho_{\Lambda}(Z)\,{\bf e}_{\gamma}=i^{-\sigma(\Lambda)}\,{\bf e}_{-\gamma}$ by \eqref{eqn:(6.1)}, we get
$$
\rho_{\Lambda}(S^{-1})\,{\bf e}_{0}
=
\rho_{\Lambda}(S)\,\rho_{\Lambda}(Z^{3})\,{\bf e}_{0}
=
i^{\sigma(\Lambda)}\frac{i^{-\frac{\sigma(\Lambda)}{2}}}{|A_{\Lambda}|^{1/2}}
\sum_{\delta\in A_{\Lambda}}{\bf e}_{\delta}
=
i^{\frac{\sigma(\Lambda)}{2}}2^{-\frac{l(\Lambda)}{2}}\sum_{\delta\in A_{\Lambda}}{\bf e}_{\delta}.
$$
This, together with the first equation of \eqref{eqn:(6.1)}, yields (1).
\newline{\bf (2)}
By \cite[p.325 l.16]{Borcherds00}, we get
$$
\rho_{\Lambda}(ST^{-2}S)\,{\bf e}_{0}
=
i^{-\sigma(\Lambda)}|A_{\Lambda}|^{-1}\sum_{\gamma,\delta\in A_{\Lambda}}
e^{2\pi i\{\langle\gamma,\delta\rangle+\gamma^{2}\}}\,{\bf e}_{\delta}
=
i^{-\sigma(\Lambda)}\,{\bf e}_{{\bf 1}_{\Lambda}},
$$
where we used the identity
$\sum_{\gamma\in A_{\Lambda}}e^{2\pi i\langle\gamma,\epsilon+\gamma\rangle}
=
\sum_{\gamma\in A_{\Lambda}}e^{2\pi i\langle\gamma,\epsilon+{\bf 1}_{\Lambda}\rangle}
=
|A_{\Lambda}|\,\delta_{{\bf 1}_{\Lambda},\epsilon}$
(cf. \cite[Lemma 3.1]{Borcherds00}) to get the second equality.
Since $S^{-1}=S^{7}=Z^{3}S$, we get
$$
\rho_{\Lambda}(V^{-1})\,{\bf e}_{0}
=
\rho_{\Lambda}(Z)^{3}\rho_{\Lambda}(ST^{-2}S)\,{\bf e}_{0}
=
i^{-\sigma(\Lambda)}\rho_{\Lambda}(Z)^{3}\,{\bf e}_{{\bf 1}_{\Lambda}}
=
i^{-\sigma(\Lambda)}i^{-3\sigma(\Lambda)}\,{\bf e}_{{\bf 1}_{\Lambda}}
=
{\bf e}_{{\bf 1}_{\Lambda}}.
$$
This proves (2).
\end{pf}

\begin{lemma}\label{Lemma7.5}
The following identities hold:
$$
(1)\quad 
f_{k}^{(0)}|_{ST^{l}}(\tau)
=
2^{\frac{8-k}{2}}i^{-\frac{k}{2}}
f_{k}^{(0)}\left(\frac{\tau+l}{4}\right),
\qquad\quad
(2)\quad
f_{k}^{(0)}|_{V}(\tau)=f_{k}^{(1)}(\tau).
$$ 
\end{lemma}

\begin{pf}
We apply \cite[Th.\,5.1]{Borcherds98} to
the lattice ${\Bbb A}_{1}^{+}=\langle2\rangle$. 
Since
$A_{{\Bbb A}_{1}^{+}}=
\langle2\rangle^{\lor}/\langle2\rangle
=\{0,\frac{1}{2}\}$, the group ring
${\bf C}[A_{{\Bbb A}_{1}^{+}}]$
is equipped with the standard basis
$\{{\bf e}_{0},{\bf e}_{1/2}\}$. Set
$\Theta_{{\Bbb A}_{1}^{+}}(\tau):=
\theta_{{\Bbb A}_{1}^{+}}(\tau)\,{\bf e}_{0}+
\theta_{{\Bbb A}_{1}^{+}+1/2}(\tau)\,{\bf e}_{1/2}$.
By \cite[Th.\,5.1]{Borcherds98} applied to 
${\Bbb A}_{1}^{+}$, we get
\begin{equation}\label{eqn:(7.1)}
\Theta_{{\Bbb A}_{1}^{+}}(g\cdot\tau)=j(g,\tau)\,\rho_{{\Bbb A}_{1}^{+}}(g)\Theta_{{\Bbb A}_{1}^{+}}(\tau),
\qquad
g\in Mp_{2}({\bf Z}).
\end{equation}
By \eqref{eqn:(6.1)} and \eqref{eqn:(7.1)}, we have
$$
\begin{array}{lll}
\Theta_{{\Bbb A}_{1}^{+}}(ST^{l}\cdot\tau)
&=&
j(ST^{l},\tau)
\left\{
\frac{{\bf e}_{0}+{\bf e}_{1/2}}{\sqrt{2i}}\,
\theta_{{\Bbb A}_{1}^{+}}(\tau)+
i^{l}\frac{{\bf e}_{0}-i{\bf e}_{1/2}}{\sqrt{2i}}\,
\theta_{{\Bbb A}_{1}^{+}+1/2}(\tau)
\right\},
\\
\Theta_{{\Bbb A}_{1}^{+}}(V\cdot\tau)
&=&
j(V,\tau)
\left\{
{\bf e}_{0}\,\theta_{{\Bbb A}_{1}^{+}+1/2}(\tau)+
{\bf e}_{1/2}\,\theta_{{\Bbb A}_{1}^{+}}(\tau)
\right\}.
\end{array}
$$
Comparing the coefficients of ${\bf e}_{0}$, we get
\begin{equation}\label{eqn:(7.2)}
\theta_{{\Bbb A}_{1}^{+}}|_{ST^{l}}(\tau)
=
(2i)^{-\frac{1}{2}}\{\theta_{{\Bbb A}_{1}^{+}}(\tau)+i^{l}\theta_{{\Bbb A}_{1}^{+}+1/2}(\tau)\}
=
(2i)^{-\frac{1}{2}}\theta_{{\Bbb A}_{1}^{+}}\left(\frac{\tau+l}{4}\right),
\end{equation}
\begin{equation}\label{eqn:(7.3)}
\theta_{{\Bbb A}_{1}^{+}}|_{V}(\tau)=\theta_{{\Bbb A}_{1}^{+}+1/2}(\tau).
\end{equation}
Here we get the second equality of \eqref{eqn:(7.2)} as follows:
$$
\theta_{{\Bbb A}_{1}^{+}}\left(\frac{\tau+l}{4}\right)
=
\sum_{n\,{\rm even}}e^{2\pi in^{2}(\tau+l)/4}
+
\sum_{n\,{\rm odd}}e^{2\pi in^{2}(\tau+l)/4}
=
\theta_{{\Bbb A}_{1}^{+}}(\tau)+i^{l}\theta_{{\Bbb A}_{1}^{+}+1/2}(\tau).
$$
\par
Set $\eta_{1^{-8}2^{8}4^{-8}}(\tau):=\eta(\tau)^{-8}\eta(2\tau)^{8}\eta(4\tau)^{-8}$,
which is a modular form for $M\Gamma_{0}(4)$ by Lemma~\ref{Lemma7.2}.
Since $ST^{l}=(\binom{0\,\,-1}{1\quad l},\sqrt{\tau+l})$
and since $\eta(-\tau^{-1})^{8}=\tau^{4}\eta(\tau)^{8}$ by \cite[Lemma 6.1]{Borcherds00}, we get
$$
\begin{aligned}
\eta_{1^{-8}2^{8}4^{-8}}|_{ST^{l}}(\tau)
&=
(\tau+l)^{\frac{8}{2}}\,
\eta_{1^{-8}2^{8}4^{-8}}\left(-\frac{1}{\tau+l}\right)
\\
&=
(\tau+l)^{4}\,
\eta\left(-\frac{1}{\tau+l}\right)^{-8}
\eta\left(-\frac{2}{\tau+l}\right)^{8}
\eta\left(-\frac{4}{\tau+l}\right)^{-8}
\\
&=
(\tau+l)^{4}
\left(\tau+l\right)^{-4}
\left(\frac{\tau+l}{2}\right)^{4}
\left(\frac{\tau+l}{4}\right)^{-4}
\\
&\quad
\times\,
\eta(\tau+l)^{-8}
\eta\left(\frac{\tau+l}{2}\right)^{8}
\eta\left(\frac{\tau+l}{4}\right)^{-8}
=
2^{4}\eta_{1^{-8}2^{8}4^{-8}}
\left(\frac{\tau+l}{4}\right),
\end{aligned}
$$
which, together with \eqref{eqn:(7.2)}, yields (1).
\par
Since $V=(\binom{\,\,\,1\quad0}{-2\quad1},\sqrt{-2\tau+1})$ and since $\eta_{1^{-8}2^{8}4^{-8}}(\tau)$ 
has weight $-4$, we get
$$
\begin{aligned}
\eta_{1^{-8}2^{8}4^{-8}}|_{V}(\tau)
&=
(-2\tau+1)^{4}\,
\eta\left(\frac{\tau}{-2\tau+1}\right)^{-8}
\eta\left(\frac{2\tau}{-2\tau+1}\right)^{8}
\eta\left(\frac{4\tau}{-2\tau+1}\right)^{-8}
\\
&=
(-2\tau+1)^{4}\left(2-\frac{1}{\tau}\right)^{-4}
\left(1-\frac{1}{2\tau}\right)^{4}
\left(\frac{1}{2}-\frac{1}{4\tau}\right)^{-4}
\\
&\quad
\times\,
\eta\left(2-\frac{1}{\tau}\right)^{-8}
\eta\left(1-\frac{1}{2\tau}\right)^{8}
\eta\left(\frac{1}{2}-\frac{1}{4\tau}\right)^{-8}
\\
&=
2^{4}\tau^{4}\,
\eta\left(2-\frac{1}{\tau}\right)^{-8}
\eta\left(1-\frac{1}{2\tau}\right)^{8}
\eta\left(\frac{1}{2}-\frac{1}{4\tau}\right)^{-8}.
\end{aligned}
$$
We define
$h(\tau):=\eta(\tau+\frac{1}{2})^{-8}
\eta(2\tau+1)^{8}\eta(4\tau+2)^{-8}$ 
for $\tau\in{\frak H}$. Then
\begin{equation}\label{eqn:(7.4)}
\eta_{1^{-8}2^{8}4^{-8}}|_{V}(\tau)=16\tau^{4}\,h\left(-\frac{1}{4\tau}\right).
\end{equation}
Set $\zeta:=\exp(2\pi i/48)$.
Since $h(\tau)$ is equal to
$$
\begin{aligned}
\,&
\zeta^{-8+16-32}
\{q^{-\frac{8}{24}}\prod_{n=1}^{\infty}(1-(-q)^{n})^{-8}\}
\{q^{\frac{16}{24}}\prod_{n=1}^{\infty}(1-q^{2n})^{8}\}
\{q^{-\frac{32}{24}}\prod_{n=1}^{\infty}(1-q^{4n})^{-8}\}
\\
&=
-q^{-1}\prod_{n=1}^{\infty}
\{(1-q^{2n})^{-8}(1+q^{2n-1})^{-8}\}
\cdot(1-q^{2n})^{8}\cdot
\{(1-q^{2n})^{-8}(1+q^{2n})^{-8}\}
\\
&=
-q^{-1}\prod_{n=1}^{\infty}
(1-q^{2n})^{-8}(1+q^{2n})^{-8}(1+q^{2n-1})^{-8}
\end{aligned}
$$
and since we have the identities
$\vartheta_{2}(2\tau)=2q^{1/4}\prod_{n=1}^{\infty}
(1-q^{2n})(1+q^{2n})^{2}$ and
\begin{equation}
\vartheta_{3}(2\tau)=\prod_{n=1}^{\infty}(1-q^{2n})(1+q^{2n-1})^{2},
\qquad
\vartheta_{4}(2\tau)=\prod_{n=1}^{\infty}(1-q^{2n})(1-q^{2n-1})^{2}
\label{eqn:(7.5)}
\end{equation}
by \cite[p.105, Eqs.(32-36)\,]{ConwaySloan99}, 
we get
\begin{equation}\label{eqn:(7.6)}
\vartheta_{2}(2\tau)^{4}\vartheta_{3}(2\tau)^{4}
=
2^{4}q\prod_{n=1}^{\infty}(1-q^{2n})^{8}(1+q^{2n})^{8}(1+q^{2n-1})^{8}
=
-2^{4}h(\tau)^{-1}.
\end{equation}
By \cite[p.104, Eq.(20)\,]{ConwaySloan99}, we have
$$
\vartheta_{2}(-\tau^{-1})^{4}=-\tau^{2}\vartheta_{4}(\tau)^{4},
\qquad
\vartheta_{3}(-\tau^{-1})^{4}=-\tau^{2}\vartheta_{3}(\tau)^{4},
$$
which, together with \eqref{eqn:(7.6)}, yield the equality
\begin{equation}\label{eqn:(7.7)}
\begin{aligned}
h\left(-\frac{1}{4\tau}\right)
&=
-2^{4}\vartheta_{2}\left(-\frac{1}{2\tau}\right)^{-4}\vartheta_{3}\left(-\frac{1}{2\tau}\right)^{-4}
\\
&=
-\tau^{-4}\vartheta_{3}(2\tau)^{-4}\vartheta_{4}(2\tau)^{-4}
\\
&=
-\tau^{-4}\left\{\prod_{n=1}^{\infty}(1-q^{2n})^{2}(1+q^{2n-1})^{2}(1-q^{2n-1})^{2}\right\}^{-4}
\\
&=
-\tau^{-4}\left\{\prod_{n=1}^{\infty}\frac{(1-q^{2n})(1-q^{4n})(1-q^{4n-2})}{(1-q^{4n})}\right\}^{-8}
\\
&=
-\tau^{-4}\left\{\frac{\prod_{n=1}^{\infty}(1-q^{2n})^{2}}{\prod_{n=1}^{\infty}(1-q^{4n})}\right\}^{-8}
=
-\tau^{-4}\eta(2\tau)^{-16}\eta(4\tau)^{8}.
\end{aligned}
\end{equation}
Here we used \eqref{eqn:(7.5)} to get the third equality.
We deduce from \eqref{eqn:(7.4)}, \eqref{eqn:(7.7)} that
\begin{equation}\label{eqn:(7.8)}
\eta_{1^{-8}2^{8}4^{-8}}|_{V}(\tau)=-16\,\eta(2\tau)^{-16}\eta(4\tau)^{8}.
\end{equation}
We get (2) from \eqref{eqn:(7.3)} and \eqref{eqn:(7.8)}.
\end{pf}

\begin{definition}\label{Definition7.6}
For a $2$-elementary lattice $\Lambda$, define a ${\bf C}[A_{\Lambda}]$-valued holomorphic function 
$F_{\Lambda}(\tau)$ on ${\frak H}$ by
$$
\begin{aligned}
F_{\Lambda}(\tau):
&=
f^{(0)}_{8+\sigma(\Lambda)}(\tau)\,{\bf e}_{0}
+
2^{\frac{4-\sigma(\Lambda)-l(\Lambda)}{2}}\,\sum_{l=0}^{3}g^{(l)}_{8+\sigma(\Lambda)}(\tau)\,{\bf v}_{l}
+
f^{(1)}_{8+\sigma(\Lambda)}(\tau)\,{\bf e}_{{\bf 1}_{\Lambda}}
\\
&=
f^{(0)}_{8+\sigma(\Lambda)}(\tau)\,{\bf e}_{0}
+
2^{\frac{4-\sigma(\Lambda)-l(\Lambda)}{2}}\,\sum_{\gamma\in A_{\Lambda}}
g^{(2\gamma^{2})}_{8+\sigma(\Lambda)}(\tau)\,{\bf e}_{\gamma}
+
f^{(1)}_{8+\sigma(\Lambda)}(\tau)\,{\bf e}_{{\bf 1}_{\Lambda}}.
\end{aligned}
$$
\end{definition}

By the Fourier expansions of $f^{(0)}_{k}(\tau)$ and $f^{(1)}_{k}(\tau)$ at $q=0$,
we get the following Fourier expansion of $F_{\Lambda}(\tau)$ at $q=0$:
\begin{equation}\label{eqn:(7.9)}
\begin{aligned}
F_{\Lambda}(\tau)
&=
\left\{q^{-1}+24+2\,\sigma(\Lambda)+O(q)\right\}\,{\bf e}_{0}
+
2^{\frac{4-\sigma(\Lambda)-l(\Lambda)}{2}}\left\{24+2\,\sigma(\Lambda)+O(q)\right\}\,{\bf v}_{0}
\\
&\quad
+O(q^{1/4})\,{\bf v}_{1}+O(q^{1/2})\,{\bf v}_{2}
+
2^{\frac{4-\sigma(\Lambda)-l(\Lambda)}{2}}\left\{q^{-1/4}+O(q^{3/4})\right\}\,{\bf v}_{3}
\\
&\quad
-2^{12+\sigma(\Lambda)}\,q^{\frac{8+\sigma(\Lambda)}{4}}
\{1+(24+\sigma(\Lambda))\,q^{2}+O(q^{4})\}\,{\bf e}_{{\bf 1}_{\Lambda}}.
\end{aligned}
\end{equation}

\begin{theorem}\label{Theorem7.7}
\begin{itemize}
\item[(1)]
$F_{\Lambda}(\tau)=
{\mathcal B}_{\Lambda}[\eta_{1^{-8}2^{8}4^{-8}}\theta_{{\Bbb A}_{1}^{+}}^{8+\sigma(\Lambda)}](\tau)$.
In particular,
$F_{\Lambda}(\tau)$ is a modular form for $Mp_{2}({\bf Z})$ of type $\rho_{\Lambda}$ with weight 
$\sigma(\Lambda)/2$.
\item[(2)]
The group $O(\Lambda)$ preserves $F_{\Lambda}$, i.e., ${\rm Aut}(F_{\Lambda},\Lambda)=O(\Lambda)$.
\item[(3)]
If $b^{+}(\Lambda)\leq2$ and $\sigma(\Lambda)\geq-12$,
$F_{\Lambda}(\tau)$ has integral Fourier coefficients.
\end{itemize}
\end{theorem}

\begin{pf}
{\bf (1)}
Set $k=8+\sigma(\Lambda)$ and $\phi(\tau)=f_{k}^{(0)}(\tau)$ in Proposition~\ref{Proposition7.1}. 
Since  $f_{k}^{(0)}(\tau)$ is a modular form for $M\Gamma_{0}(4)$ of weight $(k-8)/2=\sigma(\Lambda)/2$ 
with character $\chi_{\theta}^{k}=\chi_{\theta}^{\sigma(\Lambda)}$ by Lemma~\ref{Lemma7.2},
${\mathcal B}_{\Lambda}[f_{k}^{(0)}](\tau)$ is a modular form for $Mp_{2}({\bf Z})$ of type $\rho_{\Lambda}$ 
with weight $\sigma(\Lambda)/2$ by Proposition~\ref{Proposition7.1}. 
We prove that $F_{\Lambda}={\mathcal B}_{\Lambda}[f_{k}^{(0)}]$.
Since $k=8+\sigma(\Lambda)$ and $|A_{\Lambda}|=2^{l(\Lambda)}$, 
we deduce from Lemmas~\ref{Lemma7.4} (1) and \ref{Lemma7.5} (1) that
\begin{equation}\label{eqn:(7.10)}
\begin{aligned}
\sum_{l=0}^{3}f_{k}^{(0)}|_{ST^{l}}(\tau)\,\rho_{\Lambda}\left((ST^{l})^{-1}\right)\,{\bf e}_{0}
&=
\sum_{l=0}^{3}2^{\frac{8-k}{2}}i^{-\frac{k}{2}}\,i^{\frac{\sigma(\Lambda)}{2}}|A_{\Lambda}|^{-\frac{1}{2}}\,
\sum_{j=0}^{3}f_{k}^{(0)}\left(\frac{\tau+l}{4}\right)\,i^{-lj}\,{\bf v}_{j}
\\
&=
2^{\frac{-\sigma(\Lambda)-l(\Lambda)}{2}}\sum_{j=0}^{3}\sum_{l=0}^{3}
f_{k}^{(0)}\left(\frac{\tau+l}{4}\right)\,i^{-lj}\,{\bf v}_{j}
\\
&=
2^{-\frac{\sigma(\Lambda)+l(\Lambda)}{2}}
\sum_{j=0}^{3}\sum_{l=0}^{3}\sum_{s\in{\bf Z}/4{\bf Z}}g_{k}^{(s)}(\tau+l)\,i^{-lj}\,{\bf v}_{j}.
\end{aligned}
\end{equation}
Recall that $f^{(0)}_{k}(\tau)=
\sum_{n=-1}^{\infty}c_{k}^{(0)}(n)\,q^{n}$.
Since $g_{k}^{(s)}(\tau)=
\sum_{n\equiv s\,\,{\rm mod}\,\,4}
c_{k}^{(0)}(n)\,q^{n/4}$,
we get
$$
g_{k}^{(s)}(\tau+l)
=
\sum_{n\equiv s\,\,{\rm mod}\,\,4}
c_{k}^{(0)}(n)\,e^{2\pi in(\tau+l)/4}
=
\sum_{n\equiv s\,\,{\rm mod}\,\,4}
c_{k}^{(0)}(n)\,i^{sl}\,q^{n/4},
$$
which yields that
$$
\sum_{l=0}^{3}i^{-jl}g_{k}^{(s)}(\tau+l)
=
\sum_{n\equiv s\,\,{\rm mod}\,\,4}
c_{k}^{(0)}(n)\,
\sum_{l=0}^{3}i^{(s-j)l}\,q^{n/4}
=
4\delta_{js}\,g_{k}^{(s)}(\tau).
$$
Hence we get
$$
\sum_{l=0}^{3}\sum_{s\in{\bf Z}/4{\bf Z}}
i^{-jl}g_{k}^{(s)}(\tau+l)
=
\sum_{s\in{\bf Z}/4{\bf Z}}4\,\delta_{sj}\,
g_{k}^{(s)}(\tau)
=
4\,g_{k}^{(j)}(\tau),
$$
which, together with \eqref{eqn:(7.10)}, yields that
\begin{equation}\label{eqn:(7.11)}
\sum_{l=0}^{3}f_{k}^{(0)}|_{ST^{l}}(\tau)\cdot\rho_{\Lambda}\left((ST^{l})^{-1}\right)\,{\bf e}_{0}
=
2^{\frac{4-\sigma(\Lambda)-l(\Lambda)}{2}}\sum_{j=0}^{3}g_{k}^{(j)}(\tau)\,{\bf v}_{j}.
\end{equation}
Similarly, we get by Lemmas~\ref{Lemma7.4} (2) and \ref{Lemma7.5} (2)
\begin{equation}\label{eqn:(7.12)}
f_{k}^{(0)}|_{V}(\tau)\,\rho_{\Lambda}(V^{-1})\,{\bf e}_{0}=f_{k}^{(1)}(\tau)\,{\bf e}_{{\bf 1}_{\Lambda}}.
\end{equation}
By \eqref{eqn:(7.11)} and \eqref{eqn:(7.12)}, we get $F_{\Lambda}={\mathcal B}_{\Lambda}[f_{k}^{(0)}]$. 
\newline{\bf (2)}
Since $g({\bf e}_{\gamma})={\bf e}_{\bar{g}(\gamma)}$ for $g\in O(\Lambda)$ and $\gamma\in A_{\Lambda}$, 
we get $g({\bf e}_{0})={\bf e}_{0}$ and $g({\bf v}_{i})={\bf v}_{i}$ for all $g\in O(\Lambda)$ 
by the definition of ${\bf v}_{i}$. 
Since ${\bf 1}_{\Lambda}$ is $O(q_{\Lambda})$-invariant by its uniqueness, 
we get $\bar{g}({\bf 1}_{\Lambda})={\bf 1}_{\Lambda}$ for all $g\in O(\Lambda)$. 
This proves ${\rm Aut}(F_{\Lambda},\Lambda)=O(\Lambda)$.
\newline{\bf (3) }
Since $f^{(0)}_{k}(\tau)$, $g^{(j)}_{k}(\tau)$, $f^{(1)}_{k}(\tau)$ have integral Fourier coefficients
for $k\geq-4$, i.e., $\sigma(\Lambda)\geq-12$,
it suffices to prove by Definition~\ref{Definition7.6} that
$2^{\frac{4-\sigma(\Lambda)-l(\Lambda)}{2}}\in{\bf Z}$ when $b^{+}(\Lambda)\leq2$.
Since $\sigma(\Lambda)=2b^{+}(\Lambda)-r(\Lambda)$, $r(\Lambda)\geq l(\Lambda)$ and
$r(\Lambda)\equiv l(\Lambda)\mod2$, we get
$4-\sigma(\Lambda)-l(\Lambda)=2(2-b^{+}(\Lambda))+r(\Lambda)-l(\Lambda)\geq0$ and
$4-\sigma(\Lambda)-l(\Lambda)\equiv0\mod2$.
\end{pf}

\par
Recall that $F_{\Lambda}$ induces a modular form $F_{\Lambda}|_{L}$ of type $\rho_{L}$
when $\Lambda={\Bbb U}(N)\oplus L$ (cf. Sect.\,\ref{subsect:6.2}).
Since $\Lambda$ is $2$-elementary, $N\in\{1,2\}$
and $L$ is $2$-elementary in this case.

\begin{lemma}\label{Lemma7.8}
If $\Lambda={\Bbb U}(N)\oplus L$, then $F_{\Lambda}|_{L}=F_{L}$.
\end{lemma}

\begin{pf}
Write $F_{\Lambda}|_{L}(\tau)=\sum_{\gamma\in A_{L}}(F_{\Lambda}|_{L})_{\gamma}(\tau)\,{\bf e}_{\gamma}$.
Since ${\bf 1}_{{\Bbb U}(N)}=(0,0)$ for $N=1,2$, we get ${\bf 1}_{\Lambda}=((0,0),{\bf 1}_{L})$.
Since $((n/N,0),\gamma)^{2}=\gamma^{2}\mod2$ for $\gamma\in A_{L}$,
it follows from Definition~\ref{Definition7.6} and the definition of $(F_{\Lambda}|_{L})(\tau)$ (cf. \eqref{eqn:(6.3)}) 
that
\begin{equation}\label{eqn:(7.13)}
(F_{\Lambda}|_{L})_{\gamma}(\tau)
=
\begin{cases}
\begin{array}{ll}
N\,2^{\frac{4-\sigma(\Lambda)-l(\Lambda)}{2}}\,g^{(2\gamma^{2})}_{8+\sigma(\Lambda)}(\tau)
&
(\gamma\not=0,{\bf 1}_{L})
\\
f^{(0)}_{8+\sigma(\Lambda)}(\tau)
+
N\,2^{\frac{4-\sigma(\Lambda)-l(\Lambda)}{2}}\,g^{(0)}_{8+\sigma(\Lambda)}(\tau)
&
(\gamma=0)
\\
f^{(1)}_{8+\sigma(\Lambda)}(\tau)
+
N\,2^{\frac{4-\sigma(\Lambda)-l(\Lambda)}{2}}\,g^{(\sigma(\Lambda))}_{8+\sigma(\Lambda)}(\tau)
&
(\gamma={\bf 1}_{L}).
\end{array}
\end{cases}
\end{equation}
In the last equality, we used the formula
${\bf 1}_{\Lambda}^{2}\equiv\frac{\sigma(\Lambda)}{2}\mod2$,
which follows from \eqref{eqn:(6.2)}, \eqref{eqn:(7.9)}.
If $N=1$, $A_{\Lambda}=A_{L}$ and hence $F_{\Lambda}|_{L}=F_{\Lambda}=F_{L}$ 
by Definition~\ref{Definition7.6} and \eqref{eqn:(7.13)}. Assume $N=2$. 
Since $\sigma(\Lambda)=\sigma(L)$ and $l(\Lambda)=l(L)+2$, we get $F_{\Lambda}|_{L}=F_{L}$
by comparing the definition of $F_{L}$ with \eqref{eqn:(7.13)}.
This proves the lemma.
\end{pf}

\begin{lemma}\label{Lemma7.9}
Let $L$ be a $2$-elementary Lorentzian lattice.
If $r(L)\leq10$, then a subset of ${\mathcal C}_{L}^{+}$ is a Weyl chamber of $L$ 
if and only if it is a Weyl chamber of $F_{L}$.
\end{lemma}

\begin{pf}
Write
$F_{L}(\tau)=
\sum_{\gamma\in A_{L}}{\bf e}_{\gamma}\sum_{k\in\frac{\gamma^{2}}{2}+{\bf Z}}c_{L,\gamma}(k)\,q^{k}$.
By \eqref{eqn:(6.4)}, it suffices to prove that if $\lambda\in L^{\lor}$, $\lambda^{2}<0$ and
$c_{L,\bar{\lambda}}(\lambda^{2}/2)\not=0$, then $h_{\lambda}=h_{d}$ for some $d\in\Delta_{L}$.
Since $8+\sigma(L)\geq0$, this follows from \eqref{eqn:(7.9)}.
\end{pf}

\section{Borcherds products for $2$-elementary lattices}
\label{sect:8}
\par
Throughout this section, we assume that $\Lambda$ is an even $2$-elementary lattice with
${\rm sign}(\Lambda)=(2,r(\Lambda)-2)$. 
Recall that ${\mathcal D}'_{\Lambda}$ and ${\mathcal D}''_{\Lambda}$ were defined in Sect.\,\ref{subsect:1.4}.

\begin{theorem}\label{Theorem8.1}
If $r(\Lambda)\leq12$, then the Borcherds lift $\Psi_{\Lambda}(\cdot,F_{\Lambda})$ 
is a holomorphic automorphic form on $\Omega_{\Lambda}^{+}$ for $O^{+}(\Lambda)$ with zero divisor
$$
{\rm div}(\Psi_{\Lambda}(\cdot,F_{\Lambda}))
=
{\mathcal D}'_{\Lambda}+(2^{(r(\Lambda)-l(\Lambda))/2}+1)\,{\mathcal D}''_{\Lambda}.
$$
The weight $w(\Lambda)$ of $\Psi_{\Lambda}(\cdot,F_{\Lambda})$ is given by the following formula:
$$
w(\Lambda)
=
\begin{cases}
\begin{array}{ll}
(16-r(\Lambda))(2^{(r(\Lambda)-l(\Lambda))/2}+1)-8(1-\delta(\Lambda))
&(r(\Lambda)=12),
\\
(16-r(\Lambda))(2^{(r(\Lambda)-l(\Lambda))/2}+1)
&(r(\Lambda)<12).
\end{array}
\end{cases}
$$
\end{theorem}

\begin{pf}
Since $r(\Lambda)\leq12$ and
${\rm sign}(\Lambda)=(2,r(\Lambda)-2)$,
we get $\sigma(\Lambda)=4-r(\Lambda)$ and $8+\sigma(\Lambda)\geq0$.
By Theorem~\ref{Theorem7.7} (2), we get 
${\rm Aut}(\Lambda,F_{\Lambda})=O(\Lambda)$.
Write
$F_{\Lambda}(\tau)=
\sum_{\gamma\in A_{\Lambda}}{\bf e}_{\gamma}
\sum_{k\in\frac{\gamma^{2}}{2}+{\bf Z}}
c_{\Lambda,\gamma}(k)\,q^{k}$.
By \eqref{eqn:(7.9)}, we see that $c_{\Lambda,\gamma}(k)\geq0$ if $k<0$
and that the coefficient of ${\bf e}_{{\bf 1}_{\Lambda}}$, i.e., $f_{8+\sigma(\Lambda)}^{(1)}(\tau)$, 
is regular at $q=0$.  
By Theorem~\ref{Theorem6.1} and \eqref{eqn:(7.9)}, 
$\Psi_{\Lambda}(\cdot,F_{\Lambda})$ is an automorphic form for $O^{+}(\Lambda)$ such that
\begin{equation}\label{eqn:(8.1)}
\begin{aligned}
{\rm div}(\Psi_{\Lambda}(\cdot,F_{\Lambda}))
&=
\sum_{\lambda\in\Lambda^{\lor}/\pm1,\,\lambda^2<0}c_{\Lambda,\bar{\lambda}}(\lambda^{2}/2)\,H_{\lambda}
\\
&=
\sum_{\lambda\in\Lambda/\pm1,\,\lambda^2/2=-1}c_{\Lambda,\bar{0}}(\lambda^{2}/2)\,H_{\lambda}
+
\sum_{\lambda\in\Lambda^{\lor}/\pm1,\,\lambda^2/2=-1/4}c_{\Lambda,\bar{\lambda}}(\lambda^{2}/2)\,H_{\lambda}
\\
&=
\sum_{\lambda\in\Delta_{\Lambda}/\pm1}H_{\lambda}
+
2^{\frac{4-\sigma(\Lambda)-l(\Lambda)}{2}}\sum_{\lambda\in\Delta''_{\Lambda}/\pm1}H_{\lambda}
=
{\mathcal D}'_{\Lambda}+(2^{\frac{r(\Lambda)-l(\Lambda)}{2}}+1)\,{\mathcal D}''_{\Lambda}.
\end{aligned}
\end{equation}
\par
By Theorem~\ref{Theorem6.1}, $w(\Lambda)=c_{\Lambda,\bar0}(0)/2$. 
If $r(\Lambda)=12$ and $\delta(\Lambda)=0$, then ${\bf 1}_{\Lambda}=0$,
which, substituted into \eqref{eqn:(7.9)}, implies that
\begin{equation}\label{eqn:(8.2)}
\begin{aligned}
F_{\Lambda}(\tau)
&=
\left\{q^{-1}+24+2\sigma(\Lambda)+O(q)\right\}\,{\bf e}_{0}
+
2^{\frac{4-\sigma(\Lambda)-l(\Lambda)}{2}}\left\{24+2\sigma(\Lambda)+O(q)\right\}\,{\bf v}_{0}
\\
&\quad
+O(q^{1/4})\,{\bf v}_{1}+O(q^{1/2})\,{\bf v}_{2}
+
2^{\frac{4-\sigma(\Lambda)-l(\Lambda)}{2}}\left\{q^{-1/4}+O(q^{3/4})\right\}\,{\bf v}_{3}
\\
&\quad
+\{-16+O(q)\}\,{\bf e}_{0}.
\end{aligned}
\end{equation}
Since ${\bf v}_{0}$ contains ${\bf e}_{0}$ with multiplicity one and since $\sigma(\Lambda)=4-r(\Lambda)$, 
we deduce from \eqref{eqn:(8.2)} that
$$
w(\Lambda)
=
\frac{c_{\Lambda,0}(0)}{2}
=
12+\sigma(\Lambda)+2^{\frac{4-\sigma(\Lambda)-l(\Lambda)}{2}}(12+\sigma(\Lambda))-8
=
(16-r(\Lambda))(2^{\frac{r(\Lambda)-l(\Lambda)}{2}}+1)-8.
$$
This proves the formula for $w(\Lambda)$ when $r(\Lambda)=12$ and $\delta(\Lambda)=0$.
\par
If $r(\Lambda)<12$ or $(r(\Lambda),\delta(\Lambda))=(12,1)$, the coefficient of ${\bf e}_{{\bf 1}_{\Lambda}}$
does not contribute to $c_{\Lambda,0}(0)$ by \eqref{eqn:(7.9)}, so that
$$
w(\Lambda)
=
\frac{c_{\Lambda,0}(0)}{2}
=
12+\sigma(\Lambda)+2^{\frac{4-\sigma(\Lambda)-l(\Lambda)}{2}}(12+\sigma(\Lambda))
=
(16-r(\Lambda))(2^{\frac{r(\Lambda)-l(\Lambda)}{2}}+1)
$$
in this case. This completes the proof of Theorem~\ref{Theorem8.1}.
\end{pf}

\begin{corollary}\label{Corollary8.2}
If $r(\Lambda)\leq12$ and $\Delta''_{\Lambda}=\emptyset$,
then ${\rm div}(\Psi_{\Lambda}(\cdot,F_{\Lambda}))={\mathcal D}_{\Lambda}$.
\end{corollary}

\begin{pf}
Since $\Delta''_{\Lambda}=\emptyset$,
the result follows from Theorem~\ref{Theorem8.1}.
\end{pf}

\begin{corollary}\label{Corollary8.3}
The moduli space of $2$-elementary $K3$ surfaces of type $M$ is quasi-affine if $r(M)\geq10$. 
\end{corollary}

\begin{pf}
Set $\Lambda:=M^{\perp}$. Then $r(\Lambda)\leq 12$.
An automorphic form on $\Omega_{\Lambda}$ is identified with 
a holomorphic section of an ample line bundle over 
${\mathcal M}_{\Lambda}^{*}$ by Baily--Borel \cite{BailyBorel66}. 
Hence ${\mathcal M}_{\Lambda}\setminus
{\rm div}(\Psi_{\Lambda}(\cdot,F_{\Lambda}))$ is 
quasi-affine. Since
${\rm supp}\,{\rm div}
(\Psi_{\Lambda}(\cdot,F_{\Lambda}))=
{\mathcal D}_{\Lambda}$ by Theorem~\ref{Theorem8.1} and hence
${\mathcal M}_{\Lambda}^{o}=
{\mathcal M}_{\Lambda}\setminus
{\rm div}(\Psi_{\Lambda}(\cdot,F_{\Lambda}))$,
we get the result.
\end{pf}

In \cite[Sect.\,2]{Nikulin80b},
\cite[Sect.\,2.2]{AlexeevNikulin06},
\cite[Sects.\,1--3]{Dolgachev96},
the notion of lattice polarized $K3$ surface was introduced. 
We follow the definition in \cite{Dolgachev96}.

\begin{corollary}\label{Corollary8.4}
If $M$ is a primitive $2$-elementary Lorentzian sublattice of ${\Bbb L}_{K3}$ with $r(M)\geq10$, 
then the moduli space of ample $M$-polarized $K3$ surfaces is quasi-affine.
\end{corollary}

\begin{pf}
Set 
$O^{+}_{0}(M^{\perp}):=
\ker\{O^{+}(M^{\perp})\to O(q_{M^{\perp}})\}$,
where $O^{+}(M^{\perp})\to O(q_{M^{\perp}})$ denotes
the natural homomorphism. 
By \cite[p.2607]{Dolgachev96},
the coarse moduli space of ample $M$-polarized 
$K3$ surfaces is isomorphic to the analytic space
$\Omega_{M^{\perp}}^{o}/O^{+}_{0}(M^{\perp})$. 
By this description, the proof of the corollary is 
similar to that of Corollary~\ref{Corollary8.3}.
\end{pf}

For the table of isometry classes of primitive 
$2$-elementary Lorentzian sublattices 
$M\subset{\Bbb L}_{K3}$ with $r(M)\geq10$, see 
\cite[Appendix, Tables 1,2,3]{FinashinKharlamov08};
there are $49$ isometry classes.
There are some examples of lattices $\Lambda$ with 
$b^{+}(\Lambda)=2$ admitting an automorphic form 
on $\Omega_{\Lambda}^{+}$ with zero divisor
${\mathcal D}_{\Lambda}$.
See \cite[Sect.\,16 Examples 1,2,3]{Borcherds95}, 
\cite{Borcherds96}, \cite[Sect.\,12]{Borcherds00}, 
\cite[Examples 2.1, 2.2]{BKPSB98}, 
\cite[II, Th.\,5.2.1]{GritsenkoNikulin98},
\cite[Th.\,6.4]{Kondo07a},
\cite[Sect.\,10]{Scheithauer06} etc.
\par
As a related result, we mention the following:

\begin{theorem}\label{Theorem8.5}
The moduli space of $2$-elementary $K3$ surfaces of type $M$ contains no complete curves if $r(M)\geq7$. 
The same is true for the moduli space of ample $M$-polarized $K3$ surfaces if $M$ is $2$-elementary 
and $r(M)\geq7$.
\end{theorem}

\begin{pf}
By \cite[Th.\,5.9]{Yoshikawa04}, $\tau_{M}$ is a strongly
pluri-subharmonic function on ${\mathcal M}_{M^{\perp}}^{o}$
if $r(M)\geq7$. Hence ${\mathcal M}_{M^{\perp}}^{o}$ 
contains no complete curves when $r(M)\geq7$.
Since the moduli space of ample $M$-polarized $K3$ surfaces
$\Omega_{M^{\perp}}^{o}/O^{+}_{0}(M^{\perp})$
is a finite covering of ${\mathcal M}_{M^{\perp}}^{o}$, 
the second assertion follows from the first one. 
\end{pf}

\begin{question}\label{Question8.6}
The existence of a strongly pluri-subharmonic
function on a quasi-projective variety $X$ does {\em not} 
necessarily imply the quasi-affineness of $X$. 
(See \cite[p.\,232 Example 3.2]{Hartshorne70} 
for a counter example.) 
If $r(M)\geq7$, is ${\mathcal M}_{M^{\perp}}^{o}$ quasi-affine?
\end{question}

Assume $\Lambda={\Bbb U}(N)\oplus L$, where $L$ is a $2$-elementary Lorentzian lattice 
with $r(L)\leq10$ and $N\in\{1,2\}$. 
Hence $r(\Lambda)\leq12$, and $F_{\Lambda}|_{L}=F_{L}$ by Lemma~\ref{Lemma7.8}. 
By \cite[Th.\,13.3]{Borcherds98}, Definition~\ref{Definition7.6} and the definitions of 
$f_{k}^{(0)}(\tau)$, $f_{k}^{(1)}(\tau)$ and $g_{k}^{(i)}(\tau)$, 
the infinite product for $\Psi_{\Lambda}(\cdot,F_{\Lambda})$ is given explicitly as follows:
\begin{equation}\label{eqn:(8.3)}
\begin{aligned}
\Psi_{\Lambda}(z,F_{\Lambda})
&=
e^{2\pi i\langle\varrho,z\rangle}
\prod_{\lambda\in L,\,\lambda\cdot{\mathcal W}>0,\,\lambda^{2}\geq-2}
(1-e^{2\pi i\langle\lambda,z\rangle})^{c_{8+\sigma(\Lambda)}^{(0)}(\lambda^2/2)}
\\
&
\quad\times
\prod_{\lambda\in 2L^{\lor},\,\lambda\cdot{\mathcal W}>0,\,\lambda^{2}\geq-2}
(1-e^{\pi iN\langle\lambda,z\rangle}
)^{2^{\frac{r(\Lambda)-l(\Lambda)}{2}}c_{8+\sigma(\Lambda)}^{(0)}(\lambda^2/2)}
\\
&
\quad\times
\prod_{\lambda\in({\bf 1}_{L}+L),\,\lambda\cdot{\mathcal W}>0,\,\lambda^{2}\geq0} 
(1-e^{2\pi i\langle\lambda,z\rangle})^{2c_{8+\sigma(\Lambda)}^{(1)}(\lambda^2/2)},
\end{aligned}
\end{equation}
where $z\in L\otimes{\bf R}+i\,{\mathcal W}$ with $({\rm Im}\,z)^{2}\gg0$, 
${\mathcal W}\subset L\otimes{\bf R}$ is a Weyl chamber of $L$ by Lemma~\ref{Lemma7.9} and 
$\varrho=\varrho(L,F_{L},{\mathcal W})\in L\otimes{\bf Q}$ is the Weyl vector of $(L,F_{L},{\mathcal W})$.

\begin{example}\label{Example8.7}
Let $\Lambda={\Bbb U}(2)\oplus{\Bbb A}_{1}^{+}\oplus{\Bbb A}_{1}^{\oplus k}$ with $0\leq k\leq8$.
By \cite[Th.\,1.1]{Yoshikawa09a}, $\Psi_{\Lambda}(\cdot,F_{\Lambda})$ is regarded as an automorphic form 
on the K\"ahler moduli of a Del Pezzo surface of degree $9-k$, which appears in the formula for 
the BCOV invariant \cite{FLY08} of certain Borcea--Voisin threefolds.
By \cite[Props.\,4.1 and 4.3]{Yoshikawa09a} and \cite[Proof of Th.\,2.3 (a) and Sect.\,3]{GritsenkoNikulin97},
there is a Borcherds--Kac--Moody superalgebra with denominator function $\Psi_{\Lambda}(\cdot,F_{\Lambda})$.
In \cite[Cors.\,3.4 and 3.5]{Gritsenko10}, Gritsenko gave a very explicit Fourier series expansion of 
$\Psi_{\Lambda}(\cdot,F_{\Lambda})$ under an appropriate identification of the domains
$\Omega_{\Lambda}^{+}$ and $\Omega_{{\Bbb U}^{\oplus2}\oplus{\Bbb D}_{k-1}}$.
\end{example}

\begin{example}\label{Example8.8}
Let $\Lambda={\Bbb U}(2)\oplus{\Bbb U}(2)\oplus{\Bbb E}_{8}(2)$. 
We have $l(\Lambda)=12$ and $w(\Lambda)=0$.
This $\Lambda$ admits no primitive embedding into ${\Bbb L}_{K3}$ by \cite[Th.\,1.12.1]{Nikulin80a}.
Since $\Delta_{\Lambda}=\emptyset$, we get ${\mathcal D}_{\Lambda}=\emptyset$, 
so that $\Psi_{\Lambda}(\cdot,F_{\Lambda})$ is a constant function. 
This $F_{\Lambda}(\tau)$ gives an example of non-trivial elliptic modular form for $Mp_{2}({\bf Z})$ 
whose Borcherds lift becomes trivial.
\end{example}

\begin{example}\label{Example8.9}
Let $\Lambda={\Bbb U}\oplus{\Bbb U}(2)\oplus{\Bbb E}_{8}(2)$.
We have $l(\Lambda)=10$ and $w(\Lambda)=4$.
Then $\Psi_{\Lambda}(\cdot,F_{\Lambda})$ is the Borcherds $\Phi$-function of dimension $10$.
See \cite[Sect.\,15, Example 4]{Borcherds95}, 
\cite{Borcherds96}, \cite[Example 13.7]{Borcherds98}, 
\cite[Sect.\,13]{FLY08},
\cite[Sect.\,11]{FreitagSalvatiManni07},
\cite[Remark 4.7, Th.\,7.1]{Kondo02}, 
\cite{Scheithauer00}, \cite[Sect.\,8.1]{Yoshikawa04}
for more about this example and related results.
\end{example}

\begin{example}\label{Example8.10}
Let $\Lambda={\Bbb U}^{2}\oplus{\Bbb E}_{8}(2)$.
We have $l(\Lambda)=8$ and $w(\Lambda)=12$. Then 
$\Psi_{\Lambda}(\cdot,F_{\Lambda})=\Psi_{\Lambda}(\cdot,\Theta_{\Lambda_{16}^{+}}(\tau)/\eta(\tau)^{24})$
is the restriction of the Borcherds $\Phi$-function of dimension $26$ to $\Omega_{\Lambda}$, 
where $\Theta_{\Lambda_{16}^{+}}(\tau)$ is the theta function \cite[Sect.\,4]{Borcherds98}
for the positive-definite $16$-dimensional Barnes--Wall lattice $\Lambda_{16}^{+}$. 
See \cite[Sect.\,8.2]{Yoshikawa04}.
\end{example}

\begin{example}\label{Example8.11}
Let $\Lambda={\Bbb U}\oplus{\Bbb U}(2)\oplus{\Bbb D}_{4}^{2}$. 
We have $l(\Lambda)=6$ and $w(\Lambda)=28$. 
Kond\=o \cite[Th.\,6.4]{Kondo07a} used $\Psi_{\Lambda}(\cdot,F_{\Lambda})$ to study the projective model 
of the moduli space of $8$ points on ${\bf P}^{1}$. 
By \cite[Th.\,6.7 and its proof]{Kondo07a}, $\Psi_{\Lambda}(\cdot,F_{\Lambda})^{15}$ is expressed
as the product of certain $105$ {\it additive} Borcherds lifts \cite[Sect.\,14]{Borcherds98}.
See also \cite[Sect.\,12]{FreitagSalvatiManni07}.
\end{example}

\begin{example}\label{Example8.12}
Let $\Lambda={\Bbb U}\oplus{\Bbb U}\oplus{\Bbb E}_{8}$.
Then $l(\Lambda)=0$ and $w(\Lambda)=252$.
We get $F_{\Lambda}(\tau)=E_{4}(\tau)^{2}/\eta(\tau)^{24}$,
where $E_{4}(\tau)$ is the Eisenstein series of weight $4$. The corresponding Borcherds lift
$\Psi_{\Lambda}(\cdot,F_{\Lambda})=\Psi_{\Lambda}(\cdot,E_{4}(\tau)^{2}/\eta(\tau)^{24})$
was introduced by Borcherds \cite[Th.\,10.1, Sect.\,16 Example 1]{Borcherds95}.
By Harvey--Moore \cite[Sects.\,4 and 5]{HarveyMoore96},
$\Psi_{\Lambda}(\cdot,E_{4}(\tau)^{2}/\eta(\tau)^{24})$
appears in the formula for the one-loop coupling renormalization 
\cite[Eqs.\,(4.1), (4.5), (4.16), (4.27)]{HarveyMoore96}. 
\end{example}

\begin{example}\label{Example8.13}
When $\Lambda={\Bbb U}^{2}\oplus{\Bbb D}_{4}$, $\Psi_{\Lambda}(\cdot,F_{\Lambda})$ coincides with
the automorphic form $\Delta$ of Freitag--Hermann \cite[Th.\,11.6]{FreitagHermann00}.
Notice that the weight of $\Delta$ is $72$ in our definition (cf. \cite[p.250 l.21--l.23]{FreitagHermann00}).
By \cite[Proof of Th.\,11.5]{FreitagHermann00},
$\Psi_{\Lambda}(\cdot,F_{\Lambda})$ is expressed as the product of certain $36$ theta functions.
\end{example}

\begin{example}\label{Example8.14}
When $\Lambda=({\Bbb A}_{1}^{+})^{\oplus2}\oplus{\Bbb A}_{1}^{\oplus4}$, 
$\Psi_{\Lambda}(\cdot,F_{\Lambda})$ is the product of all even Freitag theta functions 
\cite{Yoshida97}, \cite[Th.\,7.9]{Yoshikawa07a},
so that the structure of $\Psi_{\Lambda}(\cdot,F_{\Lambda})$ is similar to that of
$\Psi_{{\Bbb U}\oplus{\Bbb U}(2)\oplus{\Bbb D}_{4}^{2}}
(\cdot,F_{{\Bbb U}\oplus{\Bbb U}(2)\oplus{\Bbb D}_{4}^{2}})$,
$\Psi_{{\Bbb U}^{2}\oplus{\Bbb D}_{4}}(\cdot,F_{{\Bbb U}^{2}\oplus{\Bbb D}_{4}})$.
For the corresponding $2$-elementary $K3$ surfaces, see \cite{Yoshida97}.
\end{example}

\begin{example}\label{Example8.15}
When $\Lambda=({\Bbb A}_{1}^{+})^{\oplus2}\oplus{\Bbb A}_{1}^{\oplus3}$, 
$\Psi_{\Lambda}(\cdot,F_{\Lambda})$ coincides with the automorphic form $\Delta_{11}$ of Gritsenko--Nikulin
\cite[II, Example 3.4 and Th.\,5.2.1]{GritsenkoNikulin98}.
When $\Lambda={\Bbb U}^{2}\oplus{\Bbb A}_{1}$, $\Psi_{\Lambda}(\cdot,F_{\Lambda})$ coincides with
the automorphic form $\Delta_{5}^{4}\Delta_{35}$ of Gritsenko--Nikulin
\cite[II, Examples 2.4 and 3.9, Th.\,5.2.1]{GritsenkoNikulin98}.
\end{example}

\begin{theorem}\label{Theorem8.16}
Let $\Lambda={\Bbb U}\oplus{\Bbb U}\oplus{\Bbb E}_{8}(2)\oplus{\Bbb A}_{1}$. 
Then the Borcherds lift $\Psi_{\Lambda}(\cdot,F_{\Lambda})$ is a meromorphic automorphic form for 
$O^{+}(\Lambda)$ of weight $15$ with zero divisor
$$
{\mathcal D}'_{\Lambda}+5{\mathcal D}''_{\Lambda}-8\,{\mathcal H}_{\Lambda}({\bf 1}_{\Lambda},-1/2).
$$
\end{theorem}

\begin{pf}
We have $r(\Lambda)=13$, $l(\Lambda)=9$, $\sigma(\Lambda)=-9$ and $\delta(\Lambda)=1$.
By Theorem~\ref{Theorem6.1} and \eqref{eqn:(7.9)}, 
the weight of $\Psi_{\Lambda}(\cdot,F_{\Lambda})$ is given by 
$(12+\sigma(\Lambda))(2^{(4-\sigma(\Lambda)-l(\Lambda))/2}+1)=15$
and the divisor of $\Psi_{\Lambda}(\cdot,F_{\Lambda})$ is given by
$$
{\mathcal D}_{\Lambda}
+
2^{\frac{4-\sigma(\Lambda)-l(\Lambda)}{2}}\,{\mathcal D}''_{\Lambda}
-
2^{12+\sigma(\Lambda)}\,{\mathcal H}_{\Lambda}({\bf 1}_{\Lambda},-1/2)
=
{\mathcal D}'_{\Lambda}+5{\mathcal D}''_{\Lambda}
-
8\,{\mathcal H}_{\Lambda}({\bf 1}_{\Lambda},-1/2),
$$
where $-2^{12+\sigma(\Lambda)}{\mathcal H}_{\Lambda}({\bf 1}_{\Lambda},-\frac{1}{2})$
comes from the negative coefficient of $q^{\frac{8+\sigma(\Lambda)}{4}}
{\bf e}_{{\bf 1}_{\Lambda}}$ in \eqref{eqn:(7.9)}.
This proves the theorem.
\end{pf}

\section{An explicit formula for $\tau_{M}$}
\label{sect:9}
\par

\begin{theorem}\label{Theorem9.1}
Let $M$ be a primitive $2$-elementary Lorentzian sublattice of ${\Bbb L}_{K3}$.
If $r(M)>10$ or $(r(M),\delta(M))=(10,1)$, then there is a constant $C_{M}>0$ depending only on $M$ 
such that for every $2$-elementary $K3$ surface $(X,\iota)$ of type $M$,
$$
\tau_{M}(X,\iota)^{-2^{g(M)+1}(2^{g(M)}+1)}
=
C_{M}\,
\|\Psi_{M^{\perp}}(\overline{\varpi}_{M}(X,\iota),F_{M^{\perp}})\|^{2^{g(M)}}
\|\chi_{g(M)}(\varOmega(X^{\iota})\|^{16}.
$$
In particular, if $\ell\in{\bf Z}_{>0}$ is an integer such that 
${\mathcal F}_{g(M)}^{2^{g(M)+1}(2^{g(M)}+1)\ell}$ extends to a very ample line bundle on 
${\mathcal A}_{g(M)}^{*}$, then the following equality holds in Theorem~\ref{Theorem5.1}:
$$
\Phi_{M}
=
C_{M}^{\ell/2}\,
\Psi_{M^{\perp}}(\cdot,F_{M^{\perp}})^{2^{g(M)-1}\ell}
\otimes J_{M}^{*}\chi_{g(M)}^{8\ell}.
$$
\end{theorem}

\begin{pf}
By our assumption $r(M)\geq10$, we get $r(M^{\perp})\leq12$. 
If the equality $r(M^{\perp})=12$ holds, then $\delta(M)=1$.
We set $\Lambda=M^{\perp}$ in Theorem~\ref{Theorem8.1}.
Then we have 
$$
16-r(\Lambda)=r(M)-6,
\qquad
\frac{r(\Lambda)-l(\Lambda)}{2}=
11-\frac{r(M)+l(M)}{2}=g(M).
$$
\par
Recall that the Bergman kernel $K_{M^{\perp}}\in C^{\infty}(\Omega_{M^{\perp}}^{+})$ 
was defined in Sect.\,\ref{subsect:4.2}.
Let $\omega_{M^{\perp}}$ be the K\"ahler form of the Bergman metric on $\Omega_{M^{\perp}}^{+}$, 
i.e., 
$$
\omega_{M^{\perp}}:=-dd^{c}\log K_{M^{\perp}}.
$$
By \cite[Eq.\,(7.1)]{Yoshikawa04},
\cite[Th.\,4.1]{Yoshikawa09b}, we have the following equation 
of currents on $\Omega_{M^{\perp}}$:
\begin{equation}\label{eqn:(9.1)}
dd^{c}\log\tau_{M}
=
\frac{r(M)-6}{4}\,\omega_{M^{\perp}}
+
J_{M}^{*}\omega_{{\mathcal A}_{g(M)}}
-
\frac{1}{4}\delta_{{\mathcal D}_{M^{\perp}}}.
\end{equation}
By Theorem~\ref{Theorem8.1}, \eqref{eqn:(4.17)} and the Poincar\'e-Lelong formula, 
we get
\begin{equation}\label{eqn:(9.2)}
\begin{aligned}
\,&
-2^{g(M)-1}\,dd^{c}\log\|\Psi_{M^{\perp}}(\cdot,F_{M^{\perp}})\|^{2}
\\
&=
2^{g(M)-1}(2^{g(M)}+1)(r(M)-6)\,\omega_{M^{\perp}}
-
2^{g(M)-1}\delta_{{\mathcal D}'_{M^{\perp}}}
-
2^{g(M)-1}(2^{g(M)}+1)\delta_{{\mathcal D}''_{M^{\perp}}}.
\end{aligned}
\end{equation}
By Proposition~\ref{Proposition4.2} (2), there exist
$a\in{\bf Z}_{\geq0}$ and an $O^{+}(M^{\perp})$-invariant
effective divisor $E$ on $\Omega_{M^{\perp}}^{+}$ such that
\begin{equation}\label{eqn:(9.3)}
-dd^{c}\log\|J_{M}^{*}\chi_{g(M)}^{8\ell}\|^{2}
=
2^{g(M)+1}(2^{g(M)}+1)\ell\,J_{M}^{*}\omega_{{\mathcal A}_{g(M)}}
-
2(2^{2g(M)-2}+a)\ell\,\delta_{{\mathcal D}'_{M^{\perp}}}
-
\delta_{E}.
\end{equation}
By \eqref{eqn:(9.1)}, \eqref{eqn:(9.2)}, \eqref{eqn:(9.3)}, 
we get the following equation of currents on $\Omega_{M^{\perp}}^{+}$:
\begin{equation}\label{eqn:(9.4)}
\begin{aligned}
\,&
-dd^{c}\log\left[
\tau_{M}^{2^{g(M)+1}(2^{g(M)}+1)\ell}
\|\Psi_{M^{\perp}}(\cdot,F_{M^{\perp}})^{2^{g(M)-1}\ell}
\otimes J_{M}^{*}\chi_{g(M)}^{8\ell}\|^{2}
\right]
\\
&=
-2a\ell\,\delta_{{\mathcal D}'_{M^{\perp}}}-\delta_{E}.
\end{aligned}
\end{equation}
Since $\log\tau_{M}$, $\log\|\Psi_{M^{\perp}}(\cdot,F_{M^{\perp}})\|$ and $\log\|J_{M}^{*}\chi_{g(M)}^{8\ell}\|$ 
are $O^{+}(M^{\perp})$-invariant $L^{1}_{\rm loc}$-function on $\Omega_{M^{\perp}}^{+}$, 
we deduce from \eqref{eqn:(9.4)} and \cite[Th.\,3.17]{Yoshikawa04}, \cite[Eq.\,(4.8)]{Yoshikawa09b}
the existences of an integer $m$ and an $O^{+}(M^{\perp})$-invariant meromorphic function $\varphi$ 
on $\Omega_{M^{\perp}}^{+}$ with divisor $m(2a\ell\,{\mathcal D}'_{M^{\perp}}+E)$ such that
\begin{equation}\label{eqn:(9.5)}
\tau_{M}^{2^{g(M)+1}(2^{g(M)}+1)\ell}
\|\Psi_{M^{\perp}}(\cdot,F_{M^{\perp}})^{2^{g(M)-1}\ell}\otimes J_{M}^{*}\chi_{g(M)}^{8\ell}\|^{2}
=
|\varphi|^{2/m}.
\end{equation}
Since $a\geq0$, $\ell>0$ and $E$ is effective, $\varphi$ is holomorphic. 
By the $O^{+}(M^{\perp})$-invariance of $\varphi$,
there exists a holomorphic function $\widetilde{\varphi}$ on ${\mathcal M}_{M^{\perp}}$ such that 
$$
\varPi_{M^{\perp}}^{*}\widetilde{\varphi}=\varphi,
$$
where $\varPi_{M^{\perp}}\colon\Omega_{M^{\perp}}^{+}\to{\mathcal M}_{M^{\perp}}$ is the projection.
Recall that ${\mathcal M}_{M^{\perp}}^{*}$ is the Baily--Borel--Satake compactification of 
${\mathcal M}_{M^{\perp}}$. 
We define $B_{M^{\perp}}:={\mathcal M}_{M^{\perp}}^{*}\setminus{\mathcal M}_{M^{\perp}}$.
\newline{\bf Case 1 }
Assume that $r(M)\leq17$.
Since ${\mathcal M}_{M^{\perp}}^{*}$ is an irreducible normal projective variety and since
$\dim B_{M^{\perp}}\leq\dim{\mathcal M}_{M^{\perp}}^{*}-2$ by the condition $r(M)\leq17$, 
$\widetilde{\varphi}$ extends to a holomorphic function on ${\mathcal M}_{M^{\perp}}^{*}$.
Since ${\mathcal M}_{M^{\perp}}^{*}$ is compact,
$\widetilde{\varphi}$ must be a constant function on ${\mathcal M}_{M^{\perp}}^{*}$. 
Hence $a=0$, $E=0$ and $\varphi$ is a constant.
Setting $C_{M}:=|\varphi|^{-2/m}$ in \eqref{eqn:(9.5)}, we get the result.
\newline{\bf Case 2 }
Assume that $r(M)\geq18$. Then $g(M)=11-\frac{r(M)+l(M)}{2}\leq2$. 
By Proposition~\ref{Proposition4.2} (3), we get $a=0$ and $E=0$. 
Hence $\widetilde{\varphi}$ is a {\em nowhere vanishing} holomorphic function on ${\mathcal M}_{M^{\perp}}$.
By \cite[Th.\,1.1]{Yoshikawa09b}, $\widetilde{\varphi}$ has at most zeros or poles on $B_{M^{\perp}}$. 
In particular, $\widetilde{\varphi}$ extends to a meromorphic function on ${\mathcal M}_{M^{\perp}}^{*}$ 
such that ${\rm div}(\widetilde{\varphi})\subset B_{M^{\perp}}$.
Since $B_{M^{\perp}}$ is irreducible when $r(M)\geq18$ by Proposition~\ref{Proposition11.7} below, 
either ${\rm div}(\widetilde{\varphi})$ or $-{\rm div}(\widetilde{\varphi})$ is effective.
In both cases, $\widetilde{\varphi}$ must be a constant. 
This completes the proof.
\end{pf}

The following is the table of $M^{\perp}$ such that
$M$ is a primitive $2$-elementary Lorentzian sublattices 
$M\subset{\Bbb L}_{K3}$ with $r(M)>10$ or 
$(r(M),\delta(M))=(10,1)$:
\begin{table}[H]
\caption{List of $M^{\perp}$ with $r(M)>10$ or
$(r(M),\delta(M))=(10,1)$}
\begin{center}
\begin{tabular}{|c|cc|c|} \hline
$g(M)$&$M^{\perp}$ with $\delta(M^{\perp})=1$&\,
&$M^{\perp}$ with $\delta(M^{\perp})=0$
\\ \hline
$0$
&
$({\Bbb A}_{1}^{+})^{\oplus 2}\oplus{\Bbb A}_{1}^{\oplus k}$
&
($0\leq k\leq 9$)
&${\Bbb U}(2)^{\oplus 2}$
\\ \hline
$1$
&
${\Bbb U}\oplus{\Bbb A}_{1}^{+}\oplus{\Bbb A}_{1}^{\oplus k}$
&
($0\leq k\leq 9$)
&
${\Bbb U}(2)^{\oplus 2}\oplus{\Bbb D}_{4}$,
${\Bbb U}\oplus{\Bbb U}(2)$
\\ \hline
$2$
&
${\Bbb U}^{\oplus 2}\oplus{\Bbb A}_{1}^{\oplus k}$
&
($1\leq k\leq 8$)
&
${\Bbb U}\oplus{\Bbb U}(2)\oplus{\Bbb D}_{4}$,
${\Bbb U}^{\oplus 2}$
\\ \hline
$3$
&
${\Bbb U}^{\oplus 2}\oplus{\Bbb D}_{4}
\oplus{\Bbb A}_{1}^{\oplus k}$
&
($1\leq k\leq 4$)
&${\Bbb U}^{\oplus 2}\oplus{\Bbb D}_{4}$
\\ \hline
$4$
&
$({\Bbb A}_{1}^{+})^{\oplus2}\oplus{\Bbb E}_{8}
\oplus{\Bbb A}_{1}^{\oplus k}$
&
($0\leq k\leq 2$)
&\,
\\ \hline
$5$
&
${\Bbb U}\oplus{\Bbb A}_{1}^{+}\oplus
{\Bbb E}_{8}\oplus{\Bbb A}_{1}^{\oplus k}$
&
($0\leq k\leq 1$)
&\,
\\ \hline
\end{tabular}
\end{center}
\end{table}

When $(r,\delta)=(10,0)$, the same formula for $\tau_{M}$ as in Theorem~\ref{Theorem9.1} does not hold.

\begin{proposition}\label{Proposition9.2}
If $(r(M),\delta(M))=(10,0)$ and $M\not\cong{\Bbb U}(2)\oplus{\Bbb E}_{8}(2)$, then
$$
J_{M}^{o}(\Omega_{M^{\perp}}^{o})\subset\theta_{{\rm null},g(M)}.
$$
\end{proposition}

\begin{pf}
We prove that $J_{M}^{o}(\Omega_{M^{\perp}}^{o})\not\subset\theta_{{\rm null},g(M)}$ yields a contradiction.
In what follows, assume $J_{M}^{o}(\Omega_{M^{\perp}}^{o})\not\subset\theta_{{\rm null},g(M)}$. 
Since $\delta(M^{\perp})=0$ and $r(M^{\perp})=12$,
$$
\varphi
:=
\Psi_{M^{\perp}}(\cdot,F_{M^{\perp}})^{2^{g(M)-1}(2^{g(M)}+1)\ell}
\otimes
(J_{M}^{*}\chi_{g(M)}^{8\ell})^{2^{g(M)}-1}
$$
is an automorphic form on $\Omega_{M^{\perp}}^{+}$ for $O^{+}(M^{\perp})$
of weight $2^{g(M)-1}(2^{2g(M)}-1)\ell(4,4)$ by Theorem~\ref{Theorem8.1}.
Since $J_{M}^{o}(\Omega_{M^{\perp}}^{o})\not\subset\theta_{{\rm null},g(M)}$, 
we get $\varphi\not\equiv0$.
We can put $\nu=2^{g(M)-1}(2^{g(M)}+1)\ell$ in Theorem~\ref{Theorem5.1}. 
Set $\psi:=\varphi/\Phi_{M}^{2^{g(M)}-1}$. 
Since $\psi$ is an $O^{+}(M^{\perp})$-invariant meromorphic function on $\Omega_{M^{\perp}}^{+}$, 
we identify $\psi$ with the corresponding meromorphic function on ${\mathcal M}_{M^{\perp}}^{*}$. 
We compute the divisor of $\psi$.
\par
Since $\delta(M^{\perp})=0$ implies $\Delta''_{M^{\perp}}=\emptyset$, 
we get ${\mathcal D}''_{M^{\perp}}=\emptyset$. 
Since $r(M)=10$ and $M\not\cong{\Bbb U}(2)\oplus{\Bbb E}_{8}(2)$,
we get $g(M)>0$ by Proposition~\ref{Proposition2.1}.
By Proposition~\ref{Proposition4.2} (2) and Theorem~\ref{Theorem8.1}, we get
\begin{equation}\label{eqn:(9.6)}
\begin{aligned}
{\rm div}(\varphi)
&=
2^{g(M)-1}(2^{g(M)}+1)\ell\,{\mathcal D}'_{M^{\perp}}
+
(2^{g(M)}-1)\{2(2^{2g(M)-2}+a)\ell\,{\mathcal D}'_{M^{\perp}}+E\}
\\
&=
\{2^{g(M)-1}(2^{2g(M)}+1)+2a\,(2^{g(M)}-1)\}\ell\,{\mathcal D}'_{M^{\perp}}+(2^{g(M)}-1)\,E.
\end{aligned}
\end{equation}
By Theorem~\ref{Theorem5.1},
${\rm div}(\Phi_{M})=\nu\,{\mathcal D}'_{M^{\perp}}=2^{g(M)-1}(2^{g(M)}+1)\ell\,{\mathcal D}'_{M^{\perp}}$,
which, together with \eqref{eqn:(9.6)}, yields that
\begin{equation}\label{eqn:(9.7)}
\begin{aligned}
{\rm div}(\psi)
&=
{\rm div}(\varphi)-(2^{g(M)}-1)\,{\rm div}(\Phi_{M})
\\
&=
\{2^{g(M)}+2a\,(2^{g(M)}-1)\}\ell\,{\mathcal D}'_{M^{\perp}}+(2^{g(M)}-1)\,E.
\end{aligned}
\end{equation}
Since $\ell\geq1$, $a\geq0$ and since $E$ is an effective divisor, ${\rm div}(\psi)$ is a {\it non-zero} and
{\it effective} divisor on $\Omega_{M^{\perp}}^{+}$ by \eqref{eqn:(9.7)}. 
This contradicts the fact that $\psi$ descends to a meromorphic function on ${\mathcal M}_{M^{\perp}}^{*}$.
\end{pf}

When $2\leq g(M)\leq5$, one can verify Proposition~\ref{Proposition9.2} 
by using the explicit equations defining the corresponding log del Pezzo surfaces
\cite[pp.494-495, Table 14]{Nakayama07}.

\begin{theorem}\label{Theorem9.3}
If $M\cong{\Bbb A}_{1}^{+}\oplus{\Bbb A}_{1}^{\oplus8}$, 
then there exists a constant $C_{M}>0$ depending only on $M$ 
such that for every $2$-elementary $K3$ surface $(X,\iota)$ of type $M$,
$$
\tau_{M}(X,\iota)^{-40}
=
C_{M}\,\|\Psi_{M^{\perp}}(\varpi_{M}(X,\iota),F_{M^{\perp}})\|^{4}\|\chi_{g(M)}(\varOmega(X^{\iota})\|^{16}.
$$
\end{theorem}

\begin{pf}
Since $M\cong{\Bbb A}_{1}^{+}\oplus{\Bbb A}_{1}^{\oplus8}$, 
we get $M^{\perp}\cong{\Bbb U}^{\oplus2}\oplus{\Bbb E}_{8}(2)\oplus{\Bbb A}_{1}$ by e.g.
\cite[Appendix, Table 2]{FinashinKharlamov08}.
By \eqref{eqn:(9.1)} and Proposition~\ref{Proposition4.4}, we get
\begin{equation}\label{eqn:(9.8)}
\begin{aligned}
\,
&
dd^{c}\{-40\ell\,\log\tau_{M}-\log\|J_{M}^{*}\chi_{2}^{8\ell}\|^{2}\}
\\
&=
\ell\left\{
-30\,\omega_{M^{\perp}}
+
10\,\delta_{{\mathcal D}_{M^{\perp}}}
-
8\,\delta_{{\mathcal D}'_{M^{\perp}}}
-
16\,\delta_{{\mathcal H}_{M^{\perp}}({\bf 1}_{M^{\perp}},
-\frac{1}{2})}
\right\}
\\
&=
\ell
\left\{
-30\,\omega_{M^{\perp}}
+
2\delta_{{\mathcal D}'_{M^{\perp}}}
+
10\,\delta_{{\mathcal D}''_{M^{\perp}}}
-
16\,\delta_{{\mathcal H}_{M^{\perp}}({\bf 1}_{M^{\perp}},
-\frac{1}{2})}
\right\}.
\end{aligned}
\end{equation}
By \eqref{eqn:(9.8)} and \cite[Th.\,3.17]{Yoshikawa04}, 
there is a meromorphic automorphic form $\varphi_{M}$ on $\Omega_{M^{\perp}}^{+}$ 
for $O^{+}(M^{\perp})$ of weight $30\ell$ with
\begin{equation}\label{eqn:(9.9)}
{\rm div}\,\varphi_{M}
=
\ell
\left\{
2{\mathcal D}'_{M^{\perp}}+10{\mathcal D}''_{M^{\perp}}-16\,{\mathcal H}_{M^{\perp}}({\bf 1}_{M^{\perp}},-1/2)
\right\}
\end{equation}
such that
\begin{equation}\label{eqn:(9.10)}
40\ell\,\log\tau_{M}+\log\|J_{M}^{*}\chi_{2}^{8\ell}\|^{2}=-\log\|\varphi_{M}\|^{2}.
\end{equation}
Since $O^{+}(M^{\perp})/[O^{+}(M^{\perp}),O^{+}(M^{\perp})]$ is a finite Abelian group, there exists 
$\nu\in{\bf Z}_{>0}$ such that $\varphi_{M}^{\nu}$ and $\Psi_{M^{\perp}}(\cdot,F_{M^{\perp}})^{2\nu}$ 
are automorphic forms with trivial character.
By Theorem~\ref{Theorem8.16} and \eqref{eqn:(9.9)}, 
$(\Psi_{M^{\perp}}(\cdot,F_{M^{\perp}})^{2\ell}/\varphi_{M})^{\nu}$ 
is an $O^{+}(M^{\perp})$-invariant meromorphic function on $\Omega_{M^{\perp}}^{+}$ with
\begin{equation}\label{eqn:(9.11)}
\begin{aligned}
{\rm div}(\Psi_{M^{\perp}}(\cdot,F_{M^{\perp}})^{2\ell}/\varphi_{M})^{\nu}
&=
\nu\ell\{
2{\mathcal D}'_{M^{\perp}}
+
10{\mathcal D}''_{M^{\perp}}
-
16{\mathcal H}_{M^{\perp}}({\bf 1}_{M^{\perp}},-1/2)
\}
\\
&\quad
-\nu\ell
\{
2{\mathcal D}'_{M^{\perp}}
+
10\,{\mathcal D}''_{M^{\perp}}
-
16\,{\mathcal H}_{M^{\perp}}({\bf 1}_{M^{\perp}},-1/2)
\}
=0.
\end{aligned}
\end{equation}
Since ${\rm div}(\Psi_{M^{\perp}}(\cdot,F_{M^{\perp}})^{2\ell}/\varphi_{M})^{\nu}$ is empty, 
there exists a non-zero constant $C_{M}$ with
$$
\varphi_{M}=C_{M}^{\ell/2}\,\Psi_{M^{\perp}}(\cdot,F_{M^{\perp}})^{2\ell}.
$$
By \eqref{eqn:(9.10)}, \eqref{eqn:(9.11)}, we get the result.
\end{pf}

\begin{question}\label{Question9.4}
Is ${\rm div}(J_{M}^{*}\chi_{g(M)}^{8\ell})$ a linear combination of Heegner divisors on 
$\Omega_{M^{\perp}}^{+}$?
If it is the case and if $M^{\perp}\cong{\Bbb U}^{\oplus 2}\oplus K$ for some lattice $K$,
$\Phi_{M}/J_{M}^{*}\chi_{g(M)}^{8\ell}$ will be expressed as a Borcherds product by \cite[Th.\,0.8]{Bruinier02}. 
Is there a Siegel modular form $\psi$ on ${\frak S}_{g(M)}$ 
such that ${\rm div}(J_{M}^{*}\psi)$ is a linear combination of Heegner divisors on $\Omega_{M^{\perp}}^{+}$?
\end{question}

\section{Equivariant determinant of the Laplacian on real $K3$ surfaces}
\label{sect:10}
\par
In this section, we give an explicit formula for
the equivariant determinant of real $K3$ surfaces. 
We refer to \cite{DIK00}, \cite{Yoshikawa08} 
for more details about real $K3$ surfaces.
\par
The pair of a $K3$ surface and an anti-holomorphic
involution is called a real $K3$ surface.
Let $(Y,\sigma)$ be a real $K3$ surface.
There exists a primitive $2$-elementary Lorentzian 
sublattice $M\subset{\Bbb L}_{K3}$ and a marking
$\alpha$ of $Y$ such that 
$\alpha\sigma^{*}\alpha^{-1}=I_{M}$.
A holomorphic $2$-form $\eta$ on $Y$ is said to be defined
over ${\bf R}$ if $\sigma^{*}\eta=\bar{\eta}$.
Let $\gamma$ be a $\sigma$-invariant Ricci-flat
K\"ahler metric on $Y$ with volume $1$. 
Let $\Delta_{(Y,\gamma)}$ be the Laplacian of 
$(Y,\gamma)$.
Since $\sigma$ preserves $\gamma$, 
$\Delta_{(Y,\gamma)}$ commutes with the $\sigma$-action 
on $C^{\infty}(Y)$. We define
$C^{\infty}_{\pm}(Y):=
\{f\in C^{\infty}(Y);\,\sigma^{*}f=\pm f\}$,
which are preserved by $\Delta_{(Y,\gamma)}$. We set 
$\Delta_{(Y,\gamma),\pm}:=
\Delta_{(Y,\gamma)}|_{C^{\infty}_{\pm}(Y)}$.
Let $\zeta_{\pm}(Y,\gamma)(s)$ denote
the spectral zeta function of $\Delta_{(Y,\gamma),\pm}$.
Then it converges absolutely for ${\rm Re}\,s\gg0$
and extends meromorphically to
the complex plane ${\bf C}$, 
and it is holomorphic at $s=0$.
We define
$$
{\rm det}_{{\bf Z}_{2}}^{*}\Delta_{(Y,\gamma)}(\sigma)
:=\exp[-\zeta'_{+}(Y,\gamma)(0)+\zeta'_{-}(Y,\gamma)(0)].
$$
\par
Let $Y({\bf R}):=\{y\in Y;\,\sigma(y)=y\}$
be the set of real points of $(Y,\sigma)$ and 
let $Y({\bf R})=\amalg_{i}C_{i}$ be the decomposition 
into the connected components. Then $Y({\bf R})$
is the disjoint union of oriented two-dimensional
manifolds. The Riemannian metric $g|_{Y({\bf R})}$ induces
a complex structure on $Y({\bf R})$. The period of 
$Y({\bf R})$ with respect to this complex structure
is denoted by 
$\varOmega(Y({\bf R}),\gamma|_{Y({\bf R})})$.
Let $\Delta_{(C_{i},\gamma|_{C_{i}})}$ 
be the Laplacian of the Riemannian manifold
$(C_{i},\gamma|_{C_{i}})$ and let 
$\zeta(C_{i},\gamma|_{C_{i}})(s)$ 
denote the spectral zeta function of 
$\Delta_{(C_{i},\gamma|_{C_{i}})}$.
The regularized determinant of 
$\Delta_{(C_{i},\gamma|_{C_{i}})}$ is defined as 
$$
{\rm det}^{*}\Delta_{(C_{i},\gamma|_{C_{i}})}
:=
\exp\left[-\zeta(C_{i},\gamma|_{C_{i}})'(0)\right].
$$
After \cite[Def.\,4.4]{Yoshikawa08}, we define
$$
\tau(Y,\sigma,\gamma)
:=
\left\{{\rm det}_{{\bf Z}_{2}}^{*}
\Delta_{(Y,\gamma)}(\sigma)\right\}^{-2}
\prod_{i}{\rm Vol}(C_{i},\gamma|_{C_{i}})\,
({\rm det}^{*}\Delta_{(C_{i},\gamma|_{C_{i}})})^{-1}.
$$

\begin{theorem}\label{Theorem10.1}
Let $(Y,\sigma)$ be a real $K3$ surface and let $\alpha$ be a marking of $Y$ such that
$\alpha\sigma^{*}\alpha^{-1}=I_{M}$.
Let $\gamma$ be a $\sigma$-invariant Ricci-flat K\"ahler metric on $Y$ with volume $1$. 
Let $\omega_{\gamma}$ be the K\"ahler form of $\gamma$, 
and let $\eta_{\gamma}$ be a holomorphic $2$-form on $Y$ defined over ${\bf R}$ such that 
$\eta_{\gamma}\wedge\bar{\eta}_{\gamma}=2\omega_{\gamma}^{2}$. 
If $r(M)>10$ or $(r(M),\delta(M))=(10,1)$, then the following identity holds:
$$
\begin{aligned}
-4(2^{g(M)}+1)\,\log\tau(Y,\sigma,\gamma)
&=
\log\|\Psi_{M^{\perp}}(\alpha(\omega_{\gamma}+i\,{\rm Im}\,\eta_{\gamma}),F_{M^{\perp}})\|^{2}
\\
&\quad
+2^{(4-g(M))}\,\log\|\chi_{g(M)}(\varOmega(Y({\bf R}),\gamma|_{Y({\bf R})}))\|^{2}+C'_{M},
\end{aligned}
$$
where $C'_{M}=2\log C_{M}$ and $\omega_{\gamma}$, $\eta_{\gamma}$ are identified with 
their cohomology classes.
\end{theorem}

\begin{pf}
The result follows from Theorem~\ref{Theorem9.1} and
\cite[Lemma 4.5 Eq.\,(4.6)]{Yoshikawa08}.
\end{pf}

\section{Appendix}
\label{sect:11}
\par
In this section, we prove some technical results about lattices.

\subsection{A proof of the equality $\Gamma_{M}=O(M^{\perp})$}
\label{subsect:11.1}
\par
Let $M$ be a primitive sublattice of ${\Bbb L}_{K3}$ 
and set $H_{M}:={\Bbb L}_{K3}/(M\oplus M^{\perp})$.
Since ${\Bbb L}_{K3}$ is unimodular, we get
$M\oplus M^{\perp}\subset{\Bbb L}_{K3}=
{\Bbb L}_{K3}^{\lor}\subset M^{\lor}\oplus
(M^{\perp})^{\lor}$, so that
$H_{M}\subset A_{M}\oplus A_{M^{\perp}}$.
Let $p_{1}\colon H_{M}\to A_{M}$ and
$p_{2}\colon H_{M}\to A_{M^{\perp}}$ be 
the homomorphism induced by the projections
$A_{M}\oplus A_{M^{\perp}}\to A_{M}$ and
$A_{M}\oplus A_{M^{\perp}}\to A_{M^{\perp}}$, respectively. 
By \cite[Props.\,1.5.1 and 1.6.1]{Nikulin80a}, 
the following hold:

\begin{itemize}
\item[(a)]
$p_{1}$ and $p_{2}$ are isomorphisms.
\item[(b)]
$A_{M}\cong A_{M^{\perp}}$ via the isomorphism
$\gamma^{{\Bbb L}_{K3}}_{M,M^{\perp}}:=
p_{2}\circ p_{1}^{-1}$.
\item[(c)]
$q_{M^{\perp}}\circ\gamma^{{\Bbb L}_{K3}}_{M,M^{\perp}}
=-q_{M}$.
\end{itemize}

\par
Recall that $g\in O(M^{\perp})$ induces 
$\overline{g}\in O(q_{M^{\perp}})$. 
For $g\in O(M^{\perp})$, we set
$\psi_{g}:=(\gamma_{M,M^{\perp}}^{{\Bbb L}_{K3}})^{-1}
\circ\overline{g}\circ
\gamma_{M,M^{\perp}}^{{\Bbb L}_{K3}}$.
Then $\psi_{g}\in{\rm Aut}(A_{M})$.

\begin{lemma}\label{Lemma11.1}
The automorphism $\psi_{g}$ preserves $q_{M}$, i.e., $\psi_{g}\in O(q_{M})$.
\end{lemma}

\begin{pf}
The result follows from Condition (c) and the fact 
$\overline{g}\in O(q_{M^{\perp}})$.
\end{pf}

Assume that $M\subset{\Bbb L}_{K3}$ is a primitive $2$-elementary Lorentzian sublattice.
Recall that the isometry $I_{M}\in O({\Bbb L}_{K3})$ was defined in Sect.\,\ref{subsect:1.2}.
In \cite[Sect.\,1.4 (c)]{Yoshikawa04}, 
we introduced the following subgroup $\Gamma_{M}\subset O(M^{\perp})$:
$$
\Gamma_{M}:=\{g|_{M^{\perp}}\in O(M^{\perp});\,g\in O({\Bbb L}_{K3}),\,g\circ I_{M}=I_{M}\circ g\}.
$$

\begin{proposition}\label{Proposition11.2}
The following equality holds: 
$$
\Gamma_{M}=O(M^{\perp}).
$$
\end{proposition}

\begin{pf}
By the definition of $\Gamma_{M}$, it suffices to prove $O(M^{\perp})\subset\Gamma_{M}$.
Let $g\in O(M^{\perp})$ be an arbitrary element.
Since $M$ is $2$-elementary and indefinite, the natural homomorphism $O(M)\to O(q_{M})$ 
is surjective by \cite[Th.\,3.6.3]{Nikulin80a}, which implies the existence of $\varPsi_{g}\in O(M)$ 
with $\psi_{g}=\overline{\varPsi_{g}}$. 
Define $\widetilde{g}:=\varPsi_{g}\oplus g\in O(M\oplus M^{\perp})$. Then
\begin{equation}\label{eqn:(11.1)}
\gamma_{M,M^{\perp}}^{{\Bbb L}_{K3}}\circ\overline{\varPsi_{g}}
=
\gamma_{M,M^{\perp}}^{{\Bbb L}_{K3}}\circ\psi_{g}
=
\overline{g}\circ\gamma_{M,M^{\perp}}^{{\Bbb L}_{K3}}.
\end{equation}
By \eqref{eqn:(11.1)} and the criterion of Nikulin \cite[Cor.\,1.5.2]{Nikulin80a},
we get $\widetilde{g}\in O({\Bbb L}_{K3})$.
We have $\widetilde{g}\circ I_{M}=I_{M}\circ\widetilde{g}$ on $M\oplus M^{\perp}$ 
because for all $(m,n)\in M\oplus M^{\perp}$,
$$
\widetilde{g}\circ I_{M}(m,n)
=
\widetilde{g}(m,-n)
=
(\varPsi_{g}(m),-g(n))
=
I_{M}(\varPsi_{g}(m),g(n))
=
I_{M}\circ\widetilde{g}(m,n).
$$
Since $M\oplus M^{\perp}$ linearly spans ${\Bbb L}_{K3}\otimes{\bf Q}$,
we have $\widetilde{g}\circ I_{M}=I_{M}\circ\widetilde{g}$ in $O({\Bbb L}_{K3})$. 
Hence $\widetilde{g}\in\Gamma_{M}$. This proves the inclusion $O(M^{\perp})\subset\Gamma_{M}$.
\end{pf}

\subsection{A formula for $g([M\perp d])$}
\label{subsect:11.2}

\begin{lemma}\label{Lemma11.3}
Let $d\in\Delta_{M^{\perp}}$.
The smallest primitive $2$-elementary Lorentzian sublattice of ${\Bbb L}_{K3}$ containing $M\oplus{\bf Z}d$
is given by $[M\perp d]=(M^{\perp}\cap d^{\perp})^{\perp}$. 
\end{lemma}

\begin{pf}
Set $L:={\bf Z}d\cong{\Bbb A}_{1}$.
Then $[M\perp d]$ is the smallest primitive 
Lorentzian sublattice of ${\Bbb L}_{K3}$ containing 
$M\oplus L$. Since 
$M\oplus L\subset[M\perp d]\subset
[M\perp d]^{\lor}\subset M^{\lor}\oplus L^{\lor}$ 
and hence
$[M\perp d]/(M\oplus L)\subset
[M\perp d]^{\lor}/(M\oplus L)\subset A_{M}\oplus A_{L}
\cong{\bf Z}_{2}^{l(M)+1}$,
$A_{[M\perp d]}=[M\perp d]^{\lor}/[M\perp d]$ 
is a vector space over ${\bf Z}_{2}$. 
Hence $[M\perp d]$ is $2$-elementary.
\end{pf}

\begin{lemma}\label{Lemma11.4}
Let $d\in\Delta_{M^{\perp}}$. Then
$$
l([M\perp d])
=
l(M^{\perp}\cap d^{\perp})
=
\begin{cases}
\begin{array}{lll}
l(M^{\perp})+1&{\rm if}&d\in\Delta'_{M^{\perp}},
\\
l(M^{\perp})-1&{\rm if}&d\in\Delta''_{M^{\perp}}.
\end{array}
\end{cases}
$$
\end{lemma}

\begin{pf}
See \cite[Prop.\,3.1]{FinashinKharlamov08}.
\end{pf}

\begin{lemma}\label{Lemma11.5}
Let $d\in\Delta_{M^{\perp}}$. Then
$$
g([M\perp d])
=
\begin{cases}
\begin{array}{lll}
g(M)-1&{\rm if}&d\in\Delta'_{M^{\perp}},
\\
g(M)&{\rm if}&d\in\Delta''_{M^{\perp}}.
\end{array}
\end{cases}
$$
\end{lemma}

\begin{pf}
Since $r(M^{\perp}\cap d^{\perp})=r(M^{\perp})-1$
and
$$
g(M)=\{r(M^{\perp})-l(M^{\perp})\}/2,
\qquad
g([M\perp d])=
\{r(M^{\perp}\cap d^{\perp})-
l(M^{\perp}\cap d^{\perp})\}/2,
$$
the result follows from Lemma~\ref{Lemma11.4}.
\end{pf}

\subsection{The $K3$-graph}
\label{subsect:11.3}
\par
In \cite{FinashinKharlamov08}, Finashin and Kharlamov introduced
the notion of lattice graph $\Gamma_{L}$ for an even unimodular 
lattice $L$. When $L={\Bbb L}_{K3}$, the $K3$-graph 
$\Gamma_{K3}=\Gamma_{{\Bbb L}_{K3}}$ is defined as follows
(cf. \cite[Sect.\,3]{FinashinKharlamov08}):
\par
The set of vertices of $\Gamma_{K3}$, denoted by $V_{K3}$, consists of the isometry classes of 
primitive $2$-elementary Lorentzian sublattices of ${\Bbb L}_{K3}$. 
For a primitive $2$-elementary Lorentzian sublattice $M\subset{\Bbb L}_{K3}$, 
write $[M]\in V_{K3}$ for its isometry class. We identify $[M]$ with the triplet $(r(M),l(M),\delta(M))$. 
The vertex $[M]\in V_{K3}$ is even (resp. odd) if $\delta(M)=0$ (resp. $\delta(M)=1$).
In \cite{FinashinKharlamov08}, an even (resp. odd) vertex is said to be of type $I$  (resp. type $II$).
The set $V_{K3}$ was determined by Nikulin \cite{Nikulin80a, Nikulin83}.
\par
The set of oriented edges of $\Gamma_{K3}$, denoted by $E_{K3}$,
consists of the $O({\Bbb L}_{K3})$-orbits of the pairs
$(M,[d])$, where $M$ is a primitive $2$-elementary Lorentzian
sublattice of ${\Bbb L}_{K3}$ and 
$[d]\in\Delta_{M^{\perp}}/O(M^{\perp})$.
The oriented edge represented by $(M,[d])$ is denoted by $[(M,[d])]$.
Then $[(M,[d])]$ connects the vertices $[M]$ and $[M\perp d]$
with arrow starting from $[M]$ to $[M\perp d]$. 
By identifying $(M,[d])$ with the divisor
$\overline{H}_{d}=O(M^{\perp})\cdot H_{d}\subset
\overline{\mathcal D}_{M^{\perp}}$,
there is a bijection between the following sets (i), (ii):
\begin{itemize}
\item[(i)]
The edges of $\Gamma_{K3}$ starting from $[M]$.
\item[(ii)]
The irreducible components of 
$\overline{\mathcal D}_{M^{\perp}}$.
\end{itemize}
By the equivalence of (i) and (ii), 
two vertices $[M],[M']\in\Gamma_{K3}$ are connected by
an oriented edge of $\Gamma_{K3}$ going from $[M]$ to $[M']$ 
if and only if there exist $\gamma\in O({\Bbb L}_{K3})$ and 
$d\in\Delta_{M^{\perp}}$ such that
${\mathcal M}_{\gamma(M')^{\perp}}$ is an irreducible component 
of $\overline{\mathcal D}_{M^{\perp}}$.
\par
An edge $[(M,[d])]$ with $[d]\in\Delta'_{M^{\perp}}/O(M^{\perp})$
is called an odd edge.
An edge $[(M,[d])]$ with $[d]\in\Delta''_{M^{\perp}}/O(M^{\perp})$
is called an even edge. If an even edge $[(M,[d])]$ satisfies
$\delta(d^{\perp}\cap M^{\perp})=0$, then $[(M,[d])]$ is called
an even Wu edge. If $\delta(d^{\perp}\cap M^{\perp})=1$, 
$[(M,[d])]$ is called an even non-Wu edge.
The set $E_{K3}$ was determined by Finashin--Kharlamov
\cite{FinashinKharlamov08}.
See \cite[p.694 Figure 1]{FinashinKharlamov08} 
for the $K3$-graph $\Gamma_{K3}$.

\begin{proposition}\label{Proposition11.6}
The following hold:
\begin{itemize}
\item[(1)]
$(M,[d])=(M,[d'])$ in $E_{K3}$ if and only if $[M\perp d]=[M\perp d']$ in $V_{K3}$.
\item[(2)]
If $[(M,[d])]\in E_{K3}$ is odd, then
$$
(r([M\perp d]),l([M\perp d]),\delta([M\perp d]))=
(r(M)+1,l(M)+1,1).
$$
\item[(3)]
If $[(M,[d])]\in E_{K3}$ is even Wu, then
$$
(r([M\perp d]),l([M\perp d]),\delta([M\perp d]))=
(r(M)+1,l(M)-1,0).
$$ 
\item[(4)]
If $[(M,[d])]\in E_{K3}$ is even non-Wu, then
$$
(r([M\perp d]),l([M\perp d]),\delta([M\perp d]))=
(r(M)+1,l(M)-1,1).
$$ 
\item[(5)]
$\Gamma_{K3}$ contains no multiple edges. In particular,
$\#[\Delta'_{M^{\perp}}/O(M^{\perp})]\leq1$ and $\#[\Delta''_{M^{\perp}}/O(M^{\perp})]\leq2$.
\end{itemize}
\end{proposition}

\begin{pf}
The proof can be found in \cite[Sect.\,3]{FinashinKharlamov08}.
For the completeness reason, we give it here.
We get (1) by \cite[Prop.\,3.3]{FinashinKharlamov08}.
When $[(M,[d])]\in E_{K3}$ is odd, the equality
$r([M\perp d])=r(M)+1$ is trivial, the equality
$l([M\perp d])=l(M)+1$ follows from 
\cite[Prop.\,3.1]{FinashinKharlamov08}, and
the equality $\delta([M\perp d])=1$ follows from
\cite[Proof of Prop.\,3.3]{FinashinKharlamov08},
because $A_{\langle2\rangle}$ is a direct summand of
$A_{[M\perp d]}$. This proves (2).
When $[(M,[d])]\in E_{K3}$ is even Wu (resp non-Wu), 
the equality $r([M\perp d])=r(M)+1$ is trivial, the equality
$l([M\perp d])=l(M)-1$ follows from 
\cite[Prop.\,3.1]{FinashinKharlamov08}, and
the equality $\delta([M\perp d])=0$ 
(resp. $\delta([M\perp d])=1$) follows from the definition of 
Wu (resp. non-Wu) edge. 
This proves (3) and (4). We get (5) by (1), (2), (3), (4).
\end{pf}

\subsection{The irreducibility of the boundary locus: the case $r(\Lambda)\leq4$}
\label{subsect:11.4}
\par
The following was used in the proof of Theorem~\ref{Theorem9.1}.

\begin{proposition}\label{Proposition11.7}
If $r(\Lambda)\leq4$, then $B_{\Lambda}={\mathcal M}_{\Lambda}^{*}\setminus{\mathcal M}_{\Lambda}$ 
is irreducible.
\end{proposition}

\begin{pf}
When $r(\Lambda)=2$, we get $\Lambda\cong({\Bbb A}_{1}^{+})^{\oplus2}$ and 
$B_{\Lambda}=\emptyset$. Assume $r(\Lambda)\geq3$.
Let $I_{\max}(\Lambda)$ be the set of maximal primitive isotropic sublattices of $\Lambda$. 
The number of the irreducible components of $B_{\Lambda}$
of maximal dimension is given by $\#[I_{\max}(\Lambda)/O(\Lambda)]$ (cf. \cite[Sect.\,2.1]{Scattone87}). 
We must prove that when $3\leq r(\Lambda)\leq4$,
\begin{equation}\label{eqn:(11.2)}
\#[I_{\max}(\Lambda)/O(\Lambda)]\leq 1.
\end{equation}
Since $3\leq r(\Lambda)\leq4$, $\Lambda$ is one of the following $7$ lattices (cf. Table 1):
$$
{\Bbb U}(k)^{\oplus2}\,(k=1,2),
\quad
({\Bbb A}_{1}^{+})^{\oplus 2}\oplus{\Bbb A}_{1}^{\oplus l}
\,(l=1,2),
\quad
{\Bbb U}\oplus{\Bbb U}(2),
\quad
{\Bbb U}\oplus{\Bbb A}_{1}^{+}\oplus{\Bbb A}_{1},
\quad
{\Bbb U}\oplus{\Bbb A}_{1}^{+}.
$$
{\bf Case 1 }
Assume that $\Lambda={\Bbb U}(k)\oplus{\Bbb U}(k)$ $(k\leq2)$
or $\Lambda=({\Bbb A}_{1}^{+})^{\oplus 2}\oplus{\Bbb A}_{1}^{\oplus l}$ $(l=1,2)$. 
Since there exist an indefinite unimodular lattice 
$\Lambda'$ and $k\in{\bf Z}_{>0}$ with $\Lambda=\Lambda'(k)$,
we get \eqref{eqn:(11.2)} by \cite[Prop.\,1.17.1]{Nikulin80a}.
\newline{\bf Case 2 }
Assume that $\Lambda={\Bbb U}\oplus{\Bbb A}_{1}^{+}$.
There exist isomorphisms $\Omega_{\Lambda}^{+}\cong{\frak H}$ and
$O^{+}(\Lambda)\cong SL_{2}({\bf Z})$ such that 
the $O^{+}(\Lambda)$-action on $\Omega_{\Lambda}^{+}$ is 
identified with the projective action of $SL_{2}({\bf Z})$ on ${\frak H}$ via these isomorphisms 
(cf. \cite[Th.\,7.1]{Dolgachev96}). Hence
${\mathcal M}_{\Lambda}\cong SL_{2}({\bf Z})\backslash{\frak H}\cong{\bf C}$ and 
${\mathcal M}_{\Lambda}^{*}\setminus{\mathcal M}_{\Lambda}=\{+i\infty\}$, 
which implies \eqref{eqn:(11.2)} in this case.
\newline{\bf Case 3 }
Assume that $\Lambda={\Bbb U}\oplus{\Bbb U}(2)$ or
$\Lambda={\Bbb U}\oplus{\Bbb A}_{1}^{+}\oplus{\Bbb A}_{1}$.
Let $L\subset\Lambda$ be a primitive isotropic sublattice of  rank $2$. 
Let $\{{\bf e}_{1},{\bf e}_{2}\}$ be a basis of $L$.
Extending this basis of $L$, we get a basis $\{{\bf e}_{1},{\bf e}_{2},{\bf e}_{3},{\bf e}_{4}\}$
of $\Lambda$ with Gram matrix $G$ as follows:
$$
G
=
(\langle{\bf e}_{i},{\bf e}_{j}\rangle)_{1\leq i,j\leq4}
=
\begin{pmatrix}
O&A
\\
A&B
\end{pmatrix},
\qquad
A=\begin{pmatrix}
1&0
\\
0&2
\end{pmatrix},
\quad
B\in M_{2}({\bf Z}),
\quad
{}^{t}B=B.
$$
Set ${\bf e}'_{3}:={\bf e}_{3}+\alpha{\bf e}_{1}+\beta{\bf e}_{2}$ and 
${\bf e}'_{4}:={\bf e}_{4}+\gamma{\bf e}_{1}+\delta{\bf e}_{2}$,
where $\alpha,\beta,\gamma,\delta\in{\bf Z}$. The Gram matrix of $\Lambda$ with respect to
$\{{\bf e}_{1},{\bf e}_{2},{\bf e}'_{3},{\bf e}'_{4}\}$ is given by
$$
G'=
\begin{pmatrix}
O&A
\\
A&B+C
\end{pmatrix},
\qquad
C=
\begin{pmatrix}
2\alpha&\gamma+2\beta
\\
\gamma+2\beta&4\delta
\end{pmatrix}.
$$
Since $\Lambda$ is even, we can write $\langle{\bf e}_{4},{\bf e}_{4}\rangle=4k$ or $4k+2$.
We set $\alpha:=-\langle{\bf e}_{3},{\bf e}_{3}\rangle/2$, $\beta:=0$,
$\gamma:=-\langle{\bf e}_{3},{\bf e}_{4}\rangle$ and
$\delta:=-2k$. Then we get $B+C=O$ if $\delta(\Lambda)=0$
(i.e., $\langle{\bf e}_{4},{\bf e}_{4}\rangle\equiv0\mod4$)
and $B+C=\binom{0\,\,0}{0\,\,1}$ if $\delta(\Lambda)=1$
(i.e., $\langle{\bf e}_{4},{\bf e}_{4}\rangle\equiv2\mod4$).
This proves the existence of a basis $\{{\bf e}_{1,L},{\bf e}_{2,L},{\bf e}_{3,L},{\bf e}_{4,L}\}$
of $\Lambda$ with $L={\bf Z}{\bf e}_{1,L}+{\bf Z}{\bf e}_{2,L}$, 
such that the Gram matrix of $\Lambda$ with respect to this basis is of the form $\binom{O\,\,A}{A\,\,B}$.
Here $A=\binom{1\,\,0}{0\,\,2}$, $B=O$ if $\delta(\Lambda)=0$ and
$B=\binom{0\,\,0}{0\,\,2}$ if $\delta(\Lambda)=1$.
If $L'\subset\Lambda$ is another primitive isotropic sublattice of rank $2$, 
then we get an isometry of $\Lambda$ sending $L$ to $L'$ by identifying the basis 
$\{{\bf e}_{1,L},{\bf e}_{2,L},{\bf e}_{3,L},{\bf e}_{4,L}\}$ and
$\{{\bf e}_{1,L'},{\bf e}_{2,L'},{\bf e}_{3,L'},{\bf e}_{4,L'}\}$ via the map ${\bf e}_{i,L}\mapsto{\bf e}_{i,L'}$.
This proves \eqref{eqn:(11.2)}.
\end{pf}


\end{document}